\documentclass[11pt]{amsart}

\usepackage{amssymb,amsmath,amsthm,amsfonts,amscd,latexsym, dsfont, color}
\usepackage[dvipsnames]{xcolor}

\newcommand{\NAH}{\mathrm{NAH}}
\usepackage{multirow}
        
\renewcommand{\ln}{\mathsf{link}}
\newcommand{\Nest}{\mathrm{Nest}}

\newcommand{\image}{\mathsf{image}} 
 
\newcommand{\graph}{\mathsf{graph}}
\newcommand{\p}{\mathfrak{p}}
\newcommand{\Fiber}{\mathscr{M}}

\usepackage{dynkin-diagrams}
\usepackage[title]{appendix}
\usepackage[american]{babel} 

\usepackage{array}
\usepackage[citecolor=ForestGreen,linkcolor=blue,urlcolor=blue]{hyperref} 
 \usepackage{graphics}
\usepackage{wrapfig}
\usepackage{pst-node}

\usepackage{tikz}
\usetikzlibrary{shapes.geometric}
\usepackage{tikz-cd}

\hypersetup{colorlinks}
\usepackage[normalem]{ulem}

\usepackage{setspace} 

\usepackage{tcolorbox}

\usepackage{footnote}

\usepackage{cancel}

\newcommand{\dev}{\mathrm{dev}}
\newcommand{\hol}{\mathsf{hol}} 
\newcommand{\T}{\mathrm{T}}

\usepackage{MnSymbol}
\usepackage{enumitem}
\newcommand{\Fix}{\mathsf{Fix}}

\newcommand{\spann}{\mathsf{span}}

\newcommand{\V}{\mathcal{E}} 
\definecolor{bluegreen}{RGB}{5, 170, 204}

\definecolor{defnyellow}{RGB}{255,224,102} 
\definecolor{lightblue}{RGB}{170, 220, 255} 
\definecolor{darkblue}{RGB}{0, 125, 230}

\newcommand{\Gtwo}{\mathsf{G}_2}

\newcommand{\Gtwosplit}{\mathsf{G}_2'}

\newtcolorbox{bluebox}[1]{colback=lightblue,colframe=darkblue,fonttitle=\bfseries,title=#1}
\newtcolorbox{redbox}[1]{colback=red!5!white,colframe=red!75!black,fonttitle=\bfseries,title=#1}

\usepackage{pgfplots}

\usepgfplotslibrary{colormaps}
\pgfplotsset{
    compat=newest,
    colormap={mycolormap}{color=(lightgray) color=(white) color=(lightgray) } }

\definecolor{mydarkblue}{RGB}{37, 42, 200}

\newcommand{\sff}{\mathrm{I\!I}}

\definecolor{mygreen}{RGB}{0, 150, 50}

\newcommand{\vis}{\partial_{\mathrm{vis}}}

\newcommand{\Ein}{\mathsf{Ein}}

\newcommand{\imoct}{\mathrm{Im}(\Oct')}

\newcommand{\SO}{\mathrm{SO}}
\newcommand{\Gr}{\mathsf{Gr}} 
\newcommand{\Stab}{\mathsf{Stab}} 
\usepackage{mathrsfs}
\usepackage[utf8]{inputenc}
\usepackage[T1]{fontenc}
\newcommand{\R}{\mathbb R}
\newcommand{\RP}{\mathbb R \mathbb P}

\newcommand{\C}{\mathbb C}
\newcommand{\Z}{\mathbb Z}

\newcommand{\Ha}{\mathbb{H}} 
\newcommand{\F}{\mathbb{F}}
\newcommand{\g}{\mathfrak{g}}
\newcommand{\frakk}{\mathfrak{k}} 
\newcommand{\frakp}{\mathfrak{p}} 

\newcommand{\X}{\mathbb{X}}
\newcommand{\K}{\mathcal{K}}
\newcommand{\End}{\mathsf{End}} 

\newcommand{\sphere}{\mathbb{S}} 
 
\newcommand{\quadric}{\hat{ \mathbb{S}}^{2,4}}
\newcommand{\sllie}{\mathfrak{sl}}
\newcommand{\SL}{\mathrm{SL}}
\newcommand{\PSL}{\mathsf{PSL}}
\newcommand{\Sp}{\mathsf{Sp}}
\newcommand{\GL}{\mathsf{GL}}
\newcommand{\gl}{\mathfrak{gl}}

  \newcommand{\Aut}{\mathsf{Aut}}

 \newcommand{\Hom}{\mathsf{Hom}}

 \newcommand{\Ad}{\text{Ad}} 
 \newcommand{\ad}{\text{ad}} 
 
\newcommand{\delbar}{\overline{\partial}}
 \newcommand{\Der}{\mathsf{Der}}

\newcommand{\Ann}{\mathsf{Ann}}
 
\newcommand{\Iso}{\mathsf{Iso}} 
\usepackage{mathtools}

\DeclarePairedDelimiter\abs{\lvert}{\rvert}
\newcommand{\Isom}{\mathsf{Isom}} 
\newcommand{\Pho}{\mathsf{Pho}}

\newcommand{\zbar}{\bar{z}} 

\renewcommand{\Nest}{\mathrm{Nest}}

\newcommand{\im}{\text{im}}

\newcommand{\Oct}{\mathbb{O}} 
 
\newcommand{\diag}{\mathsf{diag}}

\newcommand{\Hit}{\mathsf{Hit}}

\usepackage{thmtools, thm-restate}

\newtheorem{mainthm}{Theorem}

\newtheorem{maincor}[mainthm]{Corollary}

\theoremstyle{plain}

\newtheorem{theorem}{Theorem}[section]
\newtheorem{proposition}[theorem]{Proposition}
\newtheorem{lemma}[theorem]{Lemma}
\newtheorem{corollary}[theorem]{Corollary}

\newtheorem{definition}[theorem]{Definition}
\newtheorem{remark}[theorem]{Remark}

\newcommand{\alignL}{\begin{flushleft}}
\newcommand{\alignLend}{\end{flushleft}}

\usepackage{soul}

\usepackage{bbm}

\title{$\Gtwosplit$-Slodowy Slices II: Anosov Representations}
\author{Colin Davalo and Parker Evans}
\address{Colin Davalo: Dipartimento di Matematica, Università degli Studi di Torino, Torino, Italy}
\email{colinmarcthierry.davalo@unito.it}
\address{Parker Evans: Department of Mathematics, Washington University in St. Louis, St. Louis, MO}
\email{evans.p@wustl.edu}

\usepackage[margin= 1in]{geometry}

\numberwithin{equation}{section}

\begin{document}

\begin{abstract}
We find many \emph{exotic} connected components of Anosov representations of surface groups in the split real group  $\Gtwosplit$ that are neither $\iota$-Fuchsian deformations nor higher rank Teichm\"uller spaces.
These representations are found through the non-abelian Hodge correspondence via $\Gtwosplit$-Higgs bundles that we call \emph{$\alpha$-bundles} and \emph{$\beta$-bundles}, which are particular cyclic Higgs bundles in the Slodowy slices of the $\sllie_2$-triples associated to the simple roots $\alpha$ and $\beta$ of $\g_2^\C$. We first prove that (most) $\beta$-bundles that are also Hodge bundles are $P_{\beta}$-Anosov and $\alpha$-bundles are $P_{\alpha}$-Anosov. The strategy employed is new,  geometric, and general. 
For $\sigma \in \{\alpha, \beta\}$, we then show the locally homogeneous $(\Gtwosplit, \Gtwosplit/P_{\sigma})$-manifold ${\Fiber}_{\sigma}(\rho)$ recently described in \cite{DE26b}, for the representation $\rho:\pi_1S\rightarrow \Gtwosplit$ associated to a $\sigma$-bundle, coincides up to finite cover with the corresponding geometric structure for $\rho$ defined via domains of discontinuity by Guichard-Wienhard and Kapovich-Leeb-Porti. Consequently, we find that $\pi_1({\Fiber}_{\sigma}(\rho))$ is an invariant for the connected component of $\rho$ in the space of $P_{\sigma}$-Anosov representations. After verifying the folklore result that the $\Gtwosplit$-character variety has exactly three connected components, we conclude lower bounds for the number of connected components of $P_{\sigma}$-Anosov representations in the two non-Hitchin components, most of which are exotic. These are the first higher rank examples of connected components of a character variety containing multiple (in fact, arbitrarily many) exotic components of Anosov representations. 
\end{abstract}
\maketitle

\vspace{-3ex} 

\tableofcontents

 \newpage
\section{Introduction} 
\subsection{Background}
Let $G$ be a semisimple Lie group and $S=S_g$ be a closed surface of genus $g\geq 2$. For any choice of Riemann surface $\Sigma = (S,J)$ on $S$, the  \emph{non-abelian Hodge correspondence} $\mathrm{NAH}_{\Sigma}$ provides a transcendental identification between the moduli space $\mathcal{M}_{\Sigma}(G)$ of polystable \emph{$G$-Higgs bundles} on $\Sigma$ up to gauge transformations, and the $G$-character variety $\chi(S,G) = \Hom^{red}(\pi_1S, G)/G$, the space of reductive representations $\pi_1S \to G$ up to conjugation \cite{Hit87,Sim88, Cor88, Don87}. 
We thus obtain an equivalence between two very different kinds of objects: one analytic and one purely algebraic. 
Recall that a \emph{Higgs bundle} on $\Sigma$ is a pair $(\V, \Phi)$ consisting of a holomorphic vector bundle $\V$ on $\Sigma$ and a \emph{Higgs field} $\Phi$, a holomorphic endomorphism-valued one-form $\Phi \in H^0(\End(\V) \otimes \K)$, where $\K  = (\T^{1,0}\Sigma)^*$ is the holomorphic cotangent line bundle.

To apply the map $\NAH_{\Sigma}$ and convert Higgs bundles to representations, one must solve a non-linear system of elliptic PDE called \emph{Hitchin's equations}. 
Hence, it is typically quite difficult to identify geometric features of the representation $\pi_1S \rightarrow G$ associated to a particular Higgs bundle $(\V,\Phi)$. In particular, it is difficult to see if the representation associated to a Higgs bundle is discrete. Even further, one can ask for the representation to satisfy \emph{Anosov properties}, generalizations of convex cocompactness in rank one Lie groups that imply the representation is \emph{discrete and faithful} \cite{Lab06, GW12, GGKW17,KLP17,BPS}. A bit more precisely, given a parabolic subgroup $P$ of $G$, a \emph{$P$-Anosov} representation $\pi_1S \rightarrow G$ induces an action of $\pi_1S$ on the flag manifold $G/P$ with strong contraction properties that are analogous to the source-sink dynamics exhibited by Fuchsian representations $\pi_1S\rightarrow \PSL(2,\R)$ acting on the visual boundary $\sphere^1 \cong\partial_{\infty}\Ha^2$ of hyperbolic space $\Ha^2$. 

Large families of Anosov representations are given by the $G$-\emph{Hitchin components} \cite{Hit92} and more generally the \emph{$\Theta$-positive} components of \cite{GW25}. These connected components of character varieties are the image via $\NAH_\Sigma$ of certain slices of the moduli space of Higgs bundles: the \emph{Slodowy slices} for \emph{magical} $\mathfrak{sl}_2$-triples \cite{BCGBO24}. In fact, Slodowy slices can be defined for any $\mathfrak{sl}_2$-triple in $\g$ \cite{CS21}, and provide a family of Higgs bundles that does not always fill connected components of $\mathcal{M}_{\Sigma}(G)$. Recently, some representations associated with \emph{cyclic} Higgs bundles in Slodowy slices have been shown to have remarkable properties including being Anosov \cite{Fil23,Zha25,DB25,Vir25, YaoZhang26}.

One of the defining qualities of the Slodowy slice of an $\sllie_2$-triple $\{e,f,h\}$ is that it consists of $G$-Higgs bundles $(\V,\Phi)$ whose holomorphic vector bundle $\V$ admits a `special' non-vanishing holomorphic endomorphism valued one-form $\phi_f \in H^0(\End(\V) \otimes \K)$. In particular, for any holomorphic tangent vector $v \in \T^{1,0}_p\Sigma$, the endomorphism $\phi_f(v) \in \End(\V_p)$ is in the $G$-orbit of the nilpotent $f$. 

As a model example of Slodowy slices, one can consider Hitchin's parametrization of Teichm\"uller space $\mathscr{T}(S)$ via Higgs bundles in \cite{Hit87}. 
In this case, an arbitrary Riemann surface $\Sigma_0 = (S,J)$ on $S$ is chosen, and Fuchsian representations $\pi_1S \rightarrow \PSL(2,\R)$ are realized via $\NAH_{\Sigma_0}$ through the following slice of Higgs bundles $(\V, \Phi)$ on $\Sigma_0$:  
\[ \V=\K^{1/2}\oplus \K^{-1/2}\ , \qquad \Phi=\begin{pmatrix} 0 & q_2 \\ 1 & 0 \end{pmatrix} \in H^0(\End(\V)\otimes \K)\ .\]
Here, $\K^{1/2}$ is a square root of the canonical bundle $\K$, $q_2\in H^0(\K^2)$, and $\phi_f=\begin{pmatrix} 0& 0\\ 1&0\end{pmatrix}$ is modeled on the standard $(2\times 2)$ nilpotent. Note that `1' here means 
$1 \in H^0(\Hom( \K^{1/2}, \K^{-1/2})\otimes \K) \simeq H^0(\mathcal{O})$. 

In fact, a solution to Hitchin's equations gives a Hermitian metric $h$ on $\V$ such that the connection $\nabla = \nabla^h+\Phi+\Phi^{*h}$ is flat, where $\nabla^h$ is the Chern connection of $(\V,h)$. In this case, the Fuchsian representation associated to the Higgs bundle $(\V, \Phi)$ is the holonomy of this flat connection $\nabla$.

In standard fashion, we may schematically represent the above Higgs bundle as follows: 
\begin{equation*}
    \begin{tikzcd}
     \mathcal{K}^{1/2}\arrow[r,"1"]&\mathcal{K}^{-1/2}\arrow[l, bend left,"q"]
    \end{tikzcd}\ .
\end{equation*}

We shall now transition to discuss the Higgs bundles of interest in our setting. 

\subsection{Our focus} In this paper, we continue our study from \cite{DE26b} of two families of cyclic $\Gtwosplit$-Higgs bundles, where $\Gtwosplit$ denotes the split real form of the complex exceptional Lie group $\Gtwo^{\C}$. 
We call these two families \emph{$\alpha$-bundles} and \emph{$\beta$-bundles}; the latter were introduced (not by this name) by Collier and Toulisse in \cite{CT24}. Here, $\beta$ and $\alpha$ denote short and long simple roots, respectively, in the $G_2$ root system. 
These $\sigma$-bundles, for $\sigma \in \{\alpha, \beta\}$, are $\Gtwosplit$-Higgs bundles in the $\Gtwo^{\C}$-Slodowy slices associated to the $\mathfrak{sl}_2$-triple $\{e, f, h\}$ of the simple root $\sigma$, i.e., the nilpotent $f$ is the root vector $e_{-\sigma}$. 

The holomorphic vector bundle $\V$ underlying a $\sigma$-bundle $(\V,\Phi)$ has the following form, for some holomorphic line bundles $\mathcal{L}_1,\mathcal{L}_2$:
\[ \mathcal{E} =\mathcal{L}_1\mathcal{L}_2 \oplus \mathcal{L}_2 \oplus \mathcal{L}_1 \oplus \mathcal{O} \oplus \mathcal{L}_{1}^* \oplus \mathcal{L}_{2}^* \oplus \mathcal{L}_1^*\mathcal{L}_2^*,\]
which allows us to reduce the structure group of $\V$ from $\GL(7,\C)$ to $\Gtwo^{\C}$. The $\sigma$-Slodowy condition further constrains the shape of $\V$, forcing it to admit the non-vanishing one-form $\phi_f$. For concreteness, we recall the shape of the $\sigma$-bundles and the one-form $\phi_f$ in more detail in each case. 

Fix a Riemann surface $\Sigma= (S,J)$ on $S$. 
The $\alpha$-bundles on $\Sigma$ obtain the following form: 
\begin{equation}\label{AlphaCyclic}
    \begin{tikzcd}
     \mathcal{T}^2\mathcal{K} \arrow[r, "\beta"] & \mathcal{T}\mathcal{K}\arrow[r, "1"]&\mathcal{T}\arrow[r, "-i\sqrt{2}\beta"]&\mathcal{O} \arrow[r, "-i\sqrt{2}\beta"]&\mathcal{T}^{-1} \arrow[r, "1"]&\mathcal{T}^{-1}\mathcal{K}^{-1} \arrow[r, "\beta"]\arrow[bend left, lllll, "\delta"]&\mathcal{T}^{-2}\mathcal{K}^{-1} \arrow[bend left, lllll, "\delta"]
    \end{tikzcd} , 
\end{equation} 
where $\mathcal{T}$ is a holomorphic line bundle, $\beta \in H^0(\mathcal{T}^{-1}\K)$, $\delta \in H^0(\mathcal{T}^{3}\mathcal{K}^3)$.  
Here, $\phi_f$ consists of the sum of the pair of $1$'s and is modeled on the root vector $e_{-\alpha}\in \g_2^{\C}$.

The $\beta$-bundles obtain the following similar form: 
\begin{equation}\label{BetaCyclic}
\begin{tikzcd}
\mathcal{B} \arrow[r, "1"] & \mathcal{B}\mathcal{K}^{-1}\arrow[r, "\alpha"]&\mathcal{K}\arrow[r, "-i\sqrt{2}"]&\mathcal{O} \arrow[r, "-i\sqrt{2}"]&\mathcal{K}^{-1} \arrow[r, "\alpha"]&\mathcal{B}^{-1}\mathcal{K} \arrow[r, "1"]\arrow[bend left, lllll, "\delta"]&\mathcal{B}^{-1} \arrow[bend left, lllll, "\delta"]
\end{tikzcd}.
\end{equation}
Here, $\mathcal{B}$ is a holomorphic line bundle on $\Sigma$, $\alpha\in H^0(\mathcal{B}^{-1} \mathcal{K}^3)$ and $\delta \in H^0(\mathcal{B}^2)$. 
In this case, $\phi_f$ consists of the sum of the two $1$'s and the two $-i\sqrt{2}$'s, and is modeled on the root vector $e_{-\beta} \in \g_2^{\C}$. 

\subsection{Main results}
Before stating the results, we very briefly review some $\Gtwosplit$-geometry. 

The Lie group $\Gtwosplit$ admits the following characterization: $\Gtwosplit =\Aut(\R^{3,4}, \times_{3,4})$. That is, $\Gtwosplit$ is the stabilizer in $\GL(7,\R)$ of the \emph{cross-product} $\times_{3,4}:\R^{3,4} \times \R^{3,4} \rightarrow \R^{3,4}$. The partial flag manifolds $\Gtwosplit/P_{\sigma}$ associated to short and long roots $\beta,\alpha$ are as follows: 
$\Gtwosplit/P_{\beta} \cong \Ein^{2,3}$ is the \emph{Einstein $(2,3)$-universe} of isotropic lines in $\R^{3,4}$ and $\Gtwosplit/P_{\alpha} \cong \Pho^\times$ is the space of \emph{annihilator photons} in $\R^{3,4}$. The latter space $\Pho^{\times} \subset \Gr_2(\R^{3,4})$ consists of all 2-planes for which the cross-product vanishes identically. 

\subsubsection{The Anosov property} 
We first verify new Anosov representations in $\Gtwosplit$ using $\sigma$-bundles. 

\begin{mainthm}[Anosov properties for $\sigma$-bundles]\label{MainThm:A}
Suppose that $\rho: \pi_1S \rightarrow \Gtwosplit$ is the representation associated to a polystable $\sigma$-Higgs bundle $(\V, \Phi)$. 
\begin{enumerate}[label=(\roman*), noitemsep]
    \item Let $\sigma =\beta$. If $\delta=0$, or $\alpha=0$ and $\delta$ has at least one zero, then $\rho$ is $P_{\beta}$-Anosov. 
    \item Let $\sigma =\alpha$. Then $\rho$ is $P_{\alpha}$-Anosov. 
\end{enumerate}
\end{mainthm}

Theorem \ref{MainThm:A}(i) partially answers a question of \cite{CT24} on the Anosov properties of representations associated to $\beta$-bundles. 
An $\alpha$ (resp. $\beta$)-bundle is Hitchin exactly when $\mathcal{T} \cong \K$ (resp. $\mathcal{B} \cong \K^3)$ \cite{DE26b, CT24}. 
Thus, when $\alpha,\beta \neq 0$, Theorem \ref{MainThm:A}(i),(ii) certifies the first irreducible and non-Hitchin Anosov surface group representations in $\Gtwosplit$. Note that Theorem \ref{MainThm:A}(i) also includes the $\beta$-Higgs bundles with $\alpha=0$ and $\delta$ with at least one zero, which is not new, as such representations factor through a copy of $\SO(2,2) < \Gtwosplit$. We include this case for its use in Theorems \ref{MainThm:D}, \ref{MainThm:E}.

\subsubsection{Comparison of Geometric Structures} 
In \cite{DE26b}, using the geometry of $\sigma$-bundles, we constructed locally homogeneous $(\Gtwosplit,\Gtwosplit/P_{\sigma})$-manifolds. We briefly recall the results. In particular, if $\rho: \pi_1S\rightarrow \Gtwosplit$ is the representation associated to a $\sigma$-bundle $(\V,\Phi)$, the following hold:
\begin{itemize}[itemsep = 1ex]
    \item When $\sigma = \beta$, if $(\V, \Phi)$ is polystable and $\delta$ has at least one zero, there is a canonically associated $(\Gtwosplit, \Ein^{2,3})$-manifold ${\Fiber}_\beta\rightarrow S$ with fiber $\Ein^{2,1}$ and descended holonomy $\rho$. 
    \item When $\sigma = \alpha$, if $(\V, \Phi)$ is polystable and $\delta = 0$ or $\beta=0$, there is canonically associated $(\Gtwosplit, \Pho^\times)$-manifold ${\Fiber}_\alpha\rightarrow S$ with fiber $(\RP^2\times \sphere^1)$ and descended holonomy $\rho$. 
\end{itemize}
Here, by `descended holonomy', we mean the holonomy map $\hol: \pi_1{\Fiber}_\sigma\rightarrow \Gtwosplit$ factors through the bundle projection $\pi_1{\Fiber}_\sigma\rightarrow \pi_1S$ to obtain the form $\overline{\hol}: \pi_1S\rightarrow \Gtwosplit$. 
The fiber bundles ${\Fiber}_\sigma$ are explicitly built using \emph{moving bases of pencils}, which involves the differential geometry of Busemann functions and equivariant harmonic maps to the $\Gtwosplit$-symmetric space. 

On the other hand, Anosov representations are also intimately linked with locally homogeneous $(G,X)$-structures from a different perspective.
For the group $G=\Gtwosplit$, 
by work of Guichard-Wienhard and Kapovich-Leeb-Porti, a $P_{\sigma}$-Anosov representation $\rho:\pi_1S \rightarrow \Gtwosplit$ admits a cocompact domain of discontinuity $\Omega^\rho_{\sigma}$ in the flag manifold $\mathcal{F}_{\sigma}=\Gtwosplit/P_{\sigma}$ \cite{GW12, KLP18}. The quotient 
$M_{\sigma}^{\rho} = \rho(\pi_1S) \backslash \Omega^{\rho}_{\sigma}$ is a closed manifold equipped with a locally homogeneous $(\Gtwosplit, \mathcal{F}_{\sigma})$-structure. Note that under this abstract construction, the diffeomorphism type of this manifold $M_{\sigma}^{\rho}$ is not clear a priori. 

Now, when $\rho:\pi_1S\rightarrow \Gtwosplit$ is \emph{Hitchin}, it is known a priori to be both $P_{\alpha}$ and $P_{\beta}$-Anosov by \cite{FG06}.
In the $\Gtwosplit$-Hitchin case, in \cite{DE26b}, we studied the relation between the base of pencil manifold ${\Fiber}_{\sigma}^{\rho}$ and the quotient $M_{\sigma}^{\rho}$ of the domain of discontinuity, which are constructed by completely different methods.
We now have many newly certified Anosov representations
$\rho:\pi_1S\rightarrow \Gtwosplit$ associated to $\sigma$-bundles via Theorem \ref{MainThm:A}. For these representations, a natural problem is then to determine the relationship between the $(\Gtwosplit,\mathcal{F}_{\sigma})$-manifolds ${\Fiber}_{\sigma}$ and $M_{\sigma}$. 

To describe the result, we denote the pullback of ${\Fiber}_\sigma$ under the universal covering $\pi:\tilde{S} \rightarrow S$ by $\overline{{\Fiber}_\sigma}\rightarrow \tilde{S}$. Due to the `descended holonomy', the developing map of ${\Fiber}_\sigma$ can be viewed as of type $\dev:\overline{{\Fiber}_\sigma} \rightarrow \mathcal{F}_{\sigma}$. Now, the relation between the geometric structures is as follows.

\begin{mainthm}[Quotients of domains via bases of pencils]\label{MainThm:B}
Let $\rho: \pi_1S\rightarrow \Gtwosplit$ be the representation associated to a polystable $\sigma$-Hodge bundle, satisfying $\deg(\mathcal{B})\neq 0$ in the case $\sigma=\beta$. Then the associated developing map $\dev: \overline{{\Fiber}_{\sigma}} \rightarrow \mathcal{F}_{\sigma}$ has image the domain $\Omega^{\rho}_{\sigma}$ and moreover is a finite covering. 
Consequently, ${\Fiber}_{\sigma}$ is a (possibly trivial) finite cover of $M^{\rho}_{\sigma}$. 
\end{mainthm}
Here, \emph{Hodge bundle} means that $(\V,\Phi)$ is a fixed point of the $\C^*$-action. Concretely, a $\beta$-Hodge bundle is a $\beta$-bundle with $\delta=0$ or $\alpha=0$, and an $\alpha$-Hodge bundle is an $\alpha$-bundle with $\beta=0$ or $\delta=0$.

There has been much attention given recently to the topology of the quotients $M_{\rho}= \rho(\pi_1S)\backslash \Omega_{\rho}$ for $\rho:\pi_1S\rightarrow G$ an Anosov representation and $\Omega_{\rho}$ a \cite{GW12, KLP18}-domain of discontinuity.
See \cite{DE26} or \cite{ET26} for descriptions of recent results. Due to work of Alessandrini-Maloni-Tholozan-Wienhard in \cite{AMTW25}, when $\rho$ is an $\iota$-Fuchsian deformation\footnote{Recall that an \emph{$\iota$-Fuchsian representation} $\rho:\pi_1S\rightarrow G$ is an Anosov representation that preserves a totally geodesic hyperbolic plane $\mathbb{H}^2_{\iota}$ in the Riemannian symmetric space $\mathbb{X}_G$ of $G$ associated to a homomorphism $\iota:\SL(2,\R)\rightarrow G$.}, the quotient manifold $M_{\rho}$ fibers over $S$. A difficult question already is what the fiber $\mathfrak{F}_{\rho}$ of $M_{\rho} \rightarrow S$ is. When $\rho$ is not an $\iota$-Fuchsian, it is not known generally whether $M_{\rho}$ fibers over $S$. Positive examples occur in \cite{CTT19}, and negative examples are due to \cite{GLT88} (as explained in \cite{AMTW25}). Thus, questions regarding the topology of $M_{\rho}$ are even more subtle for non-Fuchsian deformations. The value of Theorem \ref{MainThm:B} is that it describes (almost) entirely the abstract manifold $M_{\sigma}^{\rho}$ via the explicit manifold $\mathscr{M}_{\sigma}$, including the case of many non-$\iota$-Fuchsian deformations $\rho$ (see Theorem \ref{MainThm:E} and the surrounding discussion). 

\subsubsection{Verifying Anosov Properties: Nestedness of Cones}
The proofs of Theorems \ref{MainThm:A} and \ref{MainThm:B} are related to an auxiliary construction of interest, which we now describe.  

To prove the Anosov property of a representation $\rho:\pi_1S\rightarrow \Gtwosplit$ associated to a $\sigma$-bundle $(\V,\Phi)$, we make essential use of the non-vanishing 1-form $\phi_f$. Let us describe the key construction. 
Pull back the flat vector bundle $\V\rightarrow S$ under the universal covering $\pi:\tilde{S}\rightarrow S$ to obtain $\pi^*\V \rightarrow \tilde{S}$. The bundle $\pi^*\mathcal{E}$ is trivialized via $\pi^*\nabla$-parallel transport, where $\nabla = \nabla^h+\Phi+\Phi^{*h}$ is the flat connection. We identify $\pi^*\mathcal{E} \cong \tilde{\Sigma}\times \C^7$, and $\pi^*\mathcal{E}^{\R} \cong \tilde{\Sigma}\times \R^{7}$ going forward, where $\V^{\R}$ is the $\nabla$-parallel real locus.  

Now, for $\sigma = \beta$, let us set $k=1$ and for $\sigma = \alpha$, we set $k=2$. In each respective case, we use $\phi_f$ to define a $\rho$-equivariant map
$B: \T^1 \widetilde{\Sigma}\rightarrow \mathscr{B}_{k}$ to the space $\mathscr{B}_{k}$ of signature $(k,7-k)$-bilinear forms on $\R^7$. The construction is as follows: for $v \in \T^1_p\widetilde{\Sigma}$, the endomorphism $\phi_f(v) \in \End(\V^{\R}_{\mid {p}})\cong \End(\R^7)$ has an attracting $k$-plane $E(v)$ and a repelling $(7-k)$-plane $F(v)$ since $\phi_f(v)$ is $\Gtwosplit$-equivalent to the root vector $f=e_{-\sigma}$. The harmonic metric $h$ on $(\V,\Phi)$ provides a Euclidean metric on the real locus $\V^{\R}$, which we alter to obtain a bilinear form $B_v$ of signature $(k,7-k)$ as follows: 
\[ B_{v} = h|_{E(v)} - h|_{F(v)}.\]
We then obtain a \emph{cone} $C(v)$ in projective space by defining $C(v) =\{\ell \in \RP^6 \mid B_v|_{\ell} > 0\}$, namely taking all lines in $\R^7$ for which $B_v$ restricts positively. We show these cones are nested under an appropriate geodesic flow, 
and that their limiting behavior determines the corresponding Anosov limit maps. For the latter, keep in mind that $\mathcal{F}_{\beta} \subset \Gr_1(\R^7)$ and $\mathcal{F}_{\alpha} \subset \Gr_2(\R^7)$. 

\begin{mainthm}[Nestedness of cones]\label{MainThm:C}
Let $(\V,\Phi)$ be a polystable $\sigma$-bundle, with the additional demand if $\sigma =\beta$ that $\delta=0$, or $\alpha=0$ and $\delta$ has at least one zero. There is a ($\pi_1S$-invariant) negatively curved Riemannian metric $g$ on $\tilde{\Sigma}$ such that the $g$-geodesic flow $\mathcal{G}_t: \T^1\widetilde{\Sigma} \rightarrow \T^1\widetilde{\Sigma}$ satisfies the following nestedness for all $t> 0$ and $v \in \T^1\widetilde{\Sigma}$: 
\[ \overline{C({\mathcal{G}_tv}}) \subset C(v). \]
Moreover, if $v$ points towards $x \in \partial \pi_1S$, then the Anosov limit map $\xi^{\sigma}: \partial \pi_1S\rightarrow \mathcal{F}_{\sigma}$ satisfies 
\[ \mathbb{P}(\xi^{\sigma}(x))= \bigcap_{t >0} C(\mathcal{G}_t(v)).\]
\end{mainthm}

We now take a moment to clarify the relationship between the theorems mentioned so far. 
First, the nestedness of cones in Theorem \ref{MainThm:C} is used to verify the Anosov properties in Theorem \ref{MainThm:A}, analogous to the approach first innovated in \cite{BPS}. Additionally, the latter result of Theorem \ref{MainThm:C} provides partial information about the location of the limit set $\Lambda_{\sigma} = \im(\xi^{\sigma})$. Namely, take any point $x \in \partial \pi_1S$ in the Gromov boundary, and select some tangent vector $v \in \T^1\widetilde{\Sigma}$ pointing towards $x$. Then Theorem \ref{MainThm:C} guarantees that the projectivized $k$-plane $\mathbb{P}(\xi^{\sigma}(x))$ is trapped inside the cone $C(v)$. We use this information vitally in the proof of Theorem \ref{MainThm:B} to show that the developed image of $\overline{{\Fiber}_\sigma}$ avoids the thickening $K_{\Lambda_{\sigma}}$ of the limit set, namely, that $\dev$ has image in $\Omega_{\sigma}^{\rho} = \mathcal{F}_{\sigma} \setminus K_{\Lambda_{\sigma}}$. 

\subsubsection{Components of Anosov Representations} 
We now discuss consequences of Theorems \ref{MainThm:A} and \ref{MainThm:B} for the connected components of Anosov representations in $\Gtwosplit$. Let us denote $\chi^{Anosov}_{\sigma}(S,\Gtwosplit)$ as the space of $P_{\sigma}$-Anosov representations in the character variety $\chi(S,\Gtwosplit)$. In general, for $G$ a non-compact semisimple Lie group and $P_{\Theta} <G$ a parabolic subgroup, it is non-trivial task to find invariants for the connected components of $\chi^{Anosov}_{\Theta}(S,G)$. By the work of \cite{KLP18}, there is some flag manifold $\mathcal{F}_{D} = G/Q$ for which $P_{\Theta}$-Anosov representations admit a non-empty cocompact domain of discontinuity $\Omega_{\rho}^{I} \subset \mathcal{F}_{D}$ associated to some balanced ideal $I$. 
One invariant for $\chi^{Anosov}_{\Theta}(S,G)$, which is difficult to compute, is the topology of the $(G,\mathcal{F}_{D})$-manifold $M^{\rho}_{I} = \rho(\pi_1S)\backslash \Omega^{\rho}_{I}$. Indeed, it is known, by a variant of the Ehresmann-Thurston principle, that for $\rho, \rho'$ in the same component of $\chi_{\Theta}^{Anosov}(S,G)$ that $M^{\rho}_{I}$ and $M^{\rho'}_I$ are homeomorphic \cite{AMTW25}.

In our case, we use double covers $\hat{M}_\sigma^\rho$ and $\hat{\Fiber}_\sigma$, corresponding to lifting the geometric structure to a cover $\widehat{\mathcal{F}}_{D}$. 
Though we do not know $\hat{M}_{\sigma}^\rho$ and $\hat{\Fiber}_\sigma$ are homeomorphic, we can verify using geometric topology that at worst we have $\pi_1(\hat{M}^\rho_{\sigma}) \cong N\rtimes \pi_1(\hat{\Fiber}_\sigma)$ for some finite group $N$. From this fact, we show that if $\pi_1(\hat{M}^\rho_{\sigma})\cong \pi_1(\hat{M}^{\rho'}_{\sigma})$, then $\pi_1(\hat{\Fiber}_{\sigma}^{\rho})\cong \pi_1(\hat{\mathscr{M}}_{\sigma}^{\rho'})$, for $\sigma$-bundles $(\V,\Phi), (\V,\Phi')$ with associated representations $\rho,\rho'$.   Consequently, the fundamental group of the base of pencil manifold $\hat{\Fiber}_{\sigma}$ is an invariant for the component of the corresponding representation $\rho$ in $\chi_{\Theta}^{Anosov}(S,\Gtwosplit)$. 
The fiber bundle $\hat{\Fiber}_\sigma \rightarrow S$ 
has a circle bundle factor in the fiber product, which causes $\pi_1(\hat{\Fiber}_\sigma)$ to be a central extension of $\pi_1S$ by $\Z$.
In the end, $\pi_1(\hat{\Fiber}_\sigma)$ is controlled by the degree $d=\deg(\mathcal{T})$ or $d=\deg(\mathcal{B})$, which can take many values for polystable $\sigma$-bundles.
Consequently, we obtain many connected components of $P_{\sigma}$-Anosov representations. We can state the results two ways. 

First, we describe the relation between connected components of $\sigma$-bundles and the connected components of $P_{\sigma}$-Anosov representations. 
To this end, we recall the following results on the connected components of the space $\mathcal{M}_{\sigma}(\Sigma)$ of polystable $\sigma$-bundles on a fixed Riemann surface $\Sigma$: 

\begin{itemize}[itemsep = 1ex]
    \item The space $\mathcal{M}_{\alpha}(\Sigma)$ has $4g-3$ connected components $\mathcal{M}_{\alpha}^d(\Sigma)$, which are labeled by the integer $d=\deg(\mathcal{T})$ satisfying $-2g+2\leq d \leq 2g-2$ \cite{DE26b}. 
    \item The space $\mathcal{M}_{\beta}(\Sigma)$ has $6g-5$ connected components $\mathcal{M}^d_{\beta}(\Sigma)$, which are labeled by the integer $d= \deg(\mathcal{B})$ satisfying $0\leq d \leq 6g-6$ \cite{CT24}.
\end{itemize}
We find that the $\NAH$ correspondence respects the stratification of $\sigma$-bundles. 
\begin{mainthm}[Correspondence of strata]\label{MainThm:D}
Let $S$ be a closed surface of genus $g \geq2$ and $\Sigma =(S,J)$ be a Riemann surface on $S$. Let $\mathcal{M}_{\sigma}(\Sigma) \subset \mathcal{M}_{\Gtwosplit}(\Sigma)$ be the moduli space of polystable $\sigma$-bundles. 
\begin{enumerate}[label=(\roman*), itemsep = 1ex]
    \item The non-abelian Hodge correspondence $\NAH_{\Sigma}$ defines a map $\mathcal{M}_{\alpha}(\Sigma) \rightarrow \chi_{\alpha}^{Anosov}(S,\Gtwosplit)$ that is injective at the level of connected components.
    \item Let $\mathcal{M}_{\beta}^{*}(\Sigma)$ be the sub-moduli space of polystable $\beta$-Hodge bundles such that $d\neq0$. Then $\NAH_{\Sigma}$ defines a map $\mathcal{M}_{\beta}^{*}(\Sigma) \rightarrow \chi_{\beta}^{Anosov}(S,\Gtwosplit)$ that is injective at the level of connected components. 
\end{enumerate}
\end{mainthm}
The locus $\mathcal{M}^{*}_{\beta}(\Sigma)$ intersects every connected component $\mathcal{M}^d_{\beta}(\Sigma)$ of $\mathcal{M}_{\beta}(\Sigma)$, excluding, by definition, the case $d=0$. However, the associated representations when $d=0$ factor through a $\Z_2$-cover of a particular $\SL(3,\R)$-subgroup of $\Gtwosplit$ and hence can never be $\beta$-Anosov, see Remark \ref{Remk:SL3}. 
\medskip

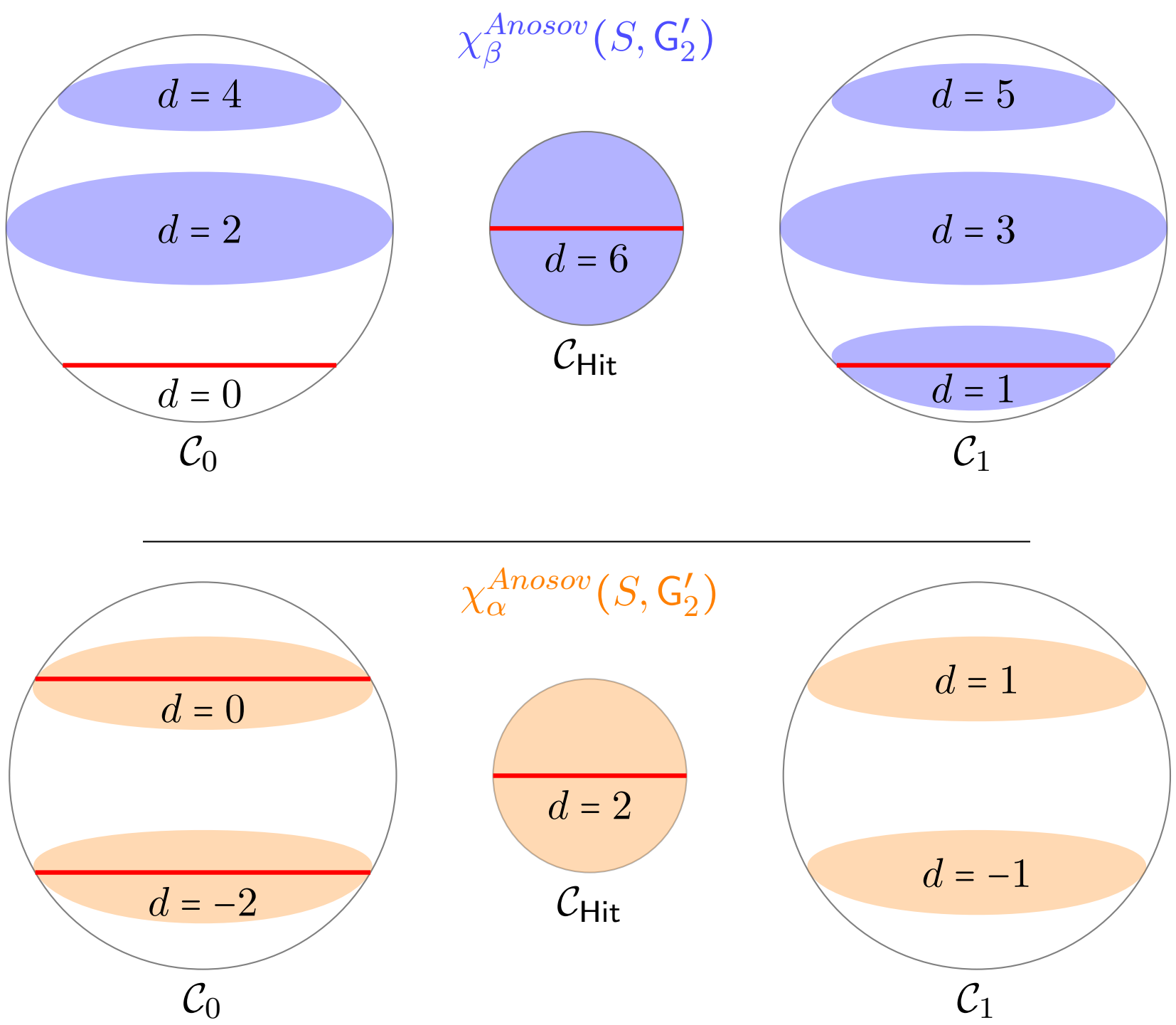
\begin{figure}[ht!]
\begin{center}
\begin{tikzpicture}[scale=.9]
\begin{scope}[
    shift={(-4,0)}]
\draw[opacity=.5] (0,0) circle (2);

\fill[blue, opacity=.3] (0:2) to[bend left=90, looseness=.5] (180:2) to[bend left=90,looseness=.5] (0:2);
\fill[blue, opacity=.3] (45:2) to[bend left=135,looseness=.7] (135:2) to[bend left=45,looseness=.5] (45:2);
\draw[red, very thick] (-45:2) -- (225:2);
\node[below] at (0,-2) {$\mathcal{C}_0$};
\node at (0,0) {$d=2$};
\node[below] at (0,-1.4) {$d=0$};
\node at (0,1.4) {$d=4$};
\end{scope}

\begin{scope}[
    shift={(4,0)}]
    \draw[opacity=.5] (0,0) circle (2);

\fill[blue, opacity=.3] (0:2) to[bend left=90, looseness=.5] (180:2) to[bend left=90,looseness=.5] (0:2);
\fill[blue, opacity=.3] (45:2) to[bend left=135,looseness=.7] (135:2) to[bend left=45,looseness=.5] (45:2);
\fill[blue, opacity=.3] (-45:2) to[bend left=-135,looseness=.7] (-135:2) to[bend left=-45,looseness=.8] (-45:2);
\draw[red, very thick] (-45:2) -- (225:2);
\node[below] at (0,-2) {$\mathcal{C}_1$};
\node at (0,0) {$d=3$};
\node[below] at (0,-1.35) {$d=1$};
\node at (0,1.4) {$d=5$};
\end{scope}

\begin{scope}[
    shift={(0,-1)}]
\draw[opacity=.5,fill=blue, fill opacity=.3] (0,1) circle (1);
\draw[red, very thick, shift={(0,1)}] (0:1) -- (180:1);
\node[below] at (0,0) {$\mathcal{C}_\Hit$};
\node[below] at (0,1) {$d=6$};
\node[above, blue, opacity =.7] at (0,2.5) {$\chi_{\beta}^{Anosov}(S,\Gtwosplit)$};
\end{scope}
\end{tikzpicture}
\end{center}

\rule{0.5\textwidth}{.4pt}

\begin{center}
\begin{tikzpicture}[scale=.9]
\begin{scope}[
    shift={(-4,0)}]
 \draw[opacity=.5] (0,0) circle (2);

\fill[orange, opacity=.3] (30:2) to[bend left=120,looseness=.6] (150:2) to[bend left=60,looseness=.5] (30:2);
 \draw[red, very thick] (30:2) -- (150:2);

\fill[orange, opacity=.3] (-30:2) to[bend left=60,looseness=.6] (210:2) to[bend left=120,looseness=.5] (-30:2);
 \draw[red, very thick] (-30:2) -- (210:2);
 \node[below] at (0,-2) {$\mathcal{C}_0$};
\node[below] at (0,-1) {$d=-2$};
 \node[below] at (0,1) {$d=0$};
\end{scope}

\begin{scope}[
    shift={(4,0)}]

 \draw[opacity=.5] (0,0) circle (2);
 \fill[orange, opacity=.3] (30:2) to[bend left=120,looseness=.5] (150:2) to[bend left=60,looseness=.5] (30:2);
 \fill[orange, opacity=.3] (-30:2) to[bend left=60,looseness=.5] (210:2) to[bend left=120,looseness=.5] (-30:2);
 \node[below] at (0,-2) {$\mathcal{C}_1$};
\node at (0,-1) {$d=-1$};
 \node at (0,1) {$d=1$};
\end{scope}

\begin{scope}[
    shift={(0,-1)}]
\draw[opacity=.3] (0,1) circle (1);
\fill[orange, opacity=.3] (0,1) circle (1);
\draw[red, very thick, shift={(0,1)}] (0:1) -- (180:1);
\node[below] at (0,0) {$\mathcal{C}_\Hit$};
\node[below] at (0,1) {$d=2$};
\node[above, orange, opacity =1] at (0,2.5) {$\chi_{\alpha}^{Anosov}(S,\Gtwosplit)$};
\end{scope}
\end{tikzpicture}
\end{center}
\caption{\small{Illustration of the connected components of $P_\sigma$-Anosov representations in the character variety $\chi(S,\Gtwosplit)$ from Theorem \ref{MainThm:E} for a surface $S$ of genus two. The red lines represent the locus of $\sigma$-bundles whose representations factor through $\SL(2,\R)$. Only three of four $\alpha$-Anosov Fuchsians appear via $\alpha$-bundles; the fourth occurs for $\beta$-bundles with $d=0$.}}
\label{fig:ThmE}
\end{figure}

The second consequence involves the relation between $\chi(S,\Gtwosplit)$ and $\chi^{Anosov}_{\sigma}(S,\Gtwosplit)$. 
Now, there are exactly three connected components in $\chi(S,\Gtwosplit)$, which we denote by $\mathcal{C}_{0}$, $\mathcal{C}_{1}$, $\mathcal{C}_\Hit$. For the details, see Appendix \ref{Apendix:CC}. 
The Hitchin component $\mathcal{C}_\Hit$ is well known to consist entirely of Borel Anosov representations. We prove that inside each of the non-Hitchin components, the open set of $P_{\sigma}$-Anosov representations decomposes into many connected components. See Figure \ref{fig:ThmE}.

\begin{mainthm}[Connected components of Anosov representations]\label{MainThm:E}
Let $S$ be a closed surface of genus $g\geq 2$ and 
write $\chi(S,\Gtwosplit) = \mathcal{C}_0 \sqcup \mathcal{C}_1 \sqcup \mathcal{C}_{\Hit}$ as above. Then we have the following: 
\begin{enumerate}[label=(\roman*),itemsep = 1ex]
    \item The set of $P_{\alpha}$-Anosov representations in $\mathcal{C}_0$ and  $\mathcal{C}_1$ has at least $2g-2$ connected components.
    \item The set of $P_{\beta}$-Anosov representations in $\mathcal{C}_0$ and $\mathcal{C}_1$ has at least $3g-4$ and $3g-3$ connected components, respectively. 
\end{enumerate}
Moreover, for $\alpha$ and $\beta$-bundles, the associated representation is in $\mathcal{C}_1$ if and only if $d$ is odd, where $d=\deg(\mathcal{T})$ or $d=\deg(\mathcal{B})$, respectively. 
\end{mainthm}

Hence, for every genus $g \geq 2$, each of the non-Hitchin components of $\mathcal{C}_0$ and $\mathcal{C}_1$ has a disconnected intersection with the space $\chi_{\sigma}^{Anosov}(S_g,\Gtwosplit)$ of $P_{\sigma}$-Anosov representations. 

Now, there are exactly $5$ distinct copies of $\mathfrak{sl}_2\R$ in the Lie algebra $\g_2'$ \cite[Appendix B]{DE26b}.  
In particular, apart from the Hitchin component, exactly one other component of $\beta$-Anosov representations contains $\iota$-Fuchsian representations and exactly three other components of $\alpha$-Anosov representations contain $\iota$-Fuchsian representations.
Consequently, Theorem \ref{MainThm:E} provides many new examples in a higher rank group $G$ of components of Anosov representations that do not contain $\iota$-Fuchsian representations, which we shall call \emph{exotic}. We can state this result quantifiably as follows. \medskip 

\begin{maincor}[Exotic components]\label{Maincor:F}
Let $S$ be a closed surface of genus $g\geq 2$. 
\begin{enumerate}[label=(\roman*),itemsep = 0ex]
    \item There are at least $3g-4$ exotic components of $P_{\beta}$-Anosov representations both in $\mathcal{C}_0$ and $\mathcal{C}_1$.
    \item There are at least $2g-5$ and $2g-2$ exotic components of $P_{\alpha}$-Anosov representations in $\mathcal{C}_0$ and $\mathcal{C}_1$, respectively.
\end{enumerate}
\end{maincor}

Unexpected components of character varieties were found by Gothen in $\Sp(4,\R)$ \cite{Got01}, which were later found to consist entirely of Anosov representations \cite{BILW05}. 
Analogous components in $\SO(p,q)$ were discovered in \cite{ABCGGO19}. 
These components of Anosov representations are exotic, 
but are also 
\emph{higher Teichm\"uller spaces} (cf. \cite{Wie18}), occupying entire components of the character variety. See \cite{BCGBO24} for a conjectural list of all such remarkable components. 

Other exotic components of Anosov representations of closed surface groups, which are not higher Teichm\"uller spaces, were detected recently by Filip \cite{Fil23} and Zhang \cite{Zha25} and by Yao-Zhang \cite{YaoZhang26}. In particular, \cite{Zha25} gives Anosov representations in $\SO_0(2,3)$ of any non-maximal Toledo invariant, many of which must be exotic. Yao-Zhang found non-Hitchin Borel-Anosov representations in $\SL(3k,\R)$ \cite{YaoZhang26}, which are necessarily exotic in $\SL(6k,\R)$. 
We note that all of these constructions detect the existence of at least one exotic connected component (though many may exist) in a given connected component of the corresponding character variety.

On the other hand, our result is the first case in higher rank in which one can distinguish several exotic connected components of Anosov representations in the \emph{same} connected component of the character variety. In fact, Corollary \ref{Maincor:F} says that the connected components $\mathcal{C}_0$, $\mathcal{C}_1$ can contain arbitrarily many exotic components of Anosov representations as the genus diverges.
The mechanism underlying the exotic nature of these components is the global topology of their associated geometric structures, which is not detected by the topological invariants of the representations. 

The first known exotic components of Anosov surface group representations in a rank one group $G$ were discovered in the striking work of Gromov-Lawson-Thurston in \cite{GLT88} for $G=\SO_0(4,1)$. Our strategy is inspired in part by their work, which verifies exotic components of Anosov representations by studying associated conformally flat geometric structures in the three-sphere that do \emph{not} fiber over $S$, meaning they cannot be $\iota$-Fuchsian deformations by \cite{AMTW25}. 
Recently, Bronstein \cite{Bro26} constructed examples of exotic Anosov representations for this same group $G$ via different analytic techniques.\footnote{In these two works, the representations are constructed into the group $G=\SO(4,1) \cong \Isom(\Ha^4)$, but one can simply take an index two subgroup that is mapped into $\SO_0(4,1)$.} Here, we recall that the $\SO_0(4,1)$-character variety of $S=S_g$ has just two connected components \cite{Bra24}. Note, however, that the constructions in these two papers apply only for surfaces of sufficiently large genus, and hence do not provide a quantitative lower bound for the components of Anosov representations of a given surface, as we obtain here.\medskip 

\subsubsection{Explicit description of limit maps}
Our final result concerns an explicit description of the Anosov limit maps associated to $\sigma$-bundles. Indeed, these limit maps are related to degenerations of 
equivariant \emph{Frenet frames} $(\mathscr{L}_p,T_p,N_p,B_p)_{p\in \widetilde{\Sigma}}$. Let us recall the $\Gtwosplit$-geometry of such objects. Given a $\sigma$-bundle, we obtain a quadruple of subspaces that is equivariant with respect to the associated representation $\rho:\pi_1S\to \Gtwosplit$, and has the following properties for all $p\in \widetilde{\Sigma}$:
\begin{itemize}
    \item $\mathscr{L}_p\oplus T_p\oplus N_p \oplus B_p$ is an orthogonal splitting of $\R^{3,4}$,
    \item The planes $T_p$, $N_p$, $B_p$ are invariant under cross product by the line $\mathscr{L}_p$,
    \item $\mathscr{L}_p$ and $N_p$ are spacelike, $T_p$ and $B_p$ are timelike. 
\end{itemize}
The space of such splittings $(\mathscr{L},T,N,B)$ is (a $\Z_2$-cover of) the homogeneous space $\Gtwosplit/H$, where $H$ is a maximal compact torus in $\Gtwosplit$.  
In our setting, the Frenet frames manifest as the Gauss lifts of equivariant $J$-holomorphic curves in the pseudosphere $\mathbb{S}^{2,4}$, as described in \cite{CT24, DE26b}.

We relate the limits of these equivariant Frenet frames to the Anosov limit maps of the associated representation as follows. For the statement, denote $\Gamma:=\pi_1S$ and $\partial \Gamma$ as the Gromov boundary. Recall that (see \S\ref{Subsec:AnosovBackground}), for $\Gamma'$ any word hyperbolic group, $P_{\sigma}$-Anosov representations $\Gamma' \rightarrow \Gtwosplit$ include as $P_{1,3}$-Anosov (resp. $P_{2}$-Anosov) in $\SL(7,\R)$ for $\sigma=\beta$ (resp. $\sigma=\alpha$). 

\begin{mainthm}[Limits of Frenet frames]
Let $(\mathscr{L}_p, T_p,N_p,B_p)_{p\in \widetilde{S}}$ be the Frenet frame associated to a $\sigma$-bundle, which we require to be a Hodge bundle such that $\delta$ has at least one zero when $\sigma=\beta$. Let $\xi^k:\partial \Gamma\to \Gr_k(\R^{3,4})$ be the Anosov limits maps of the associated representation. 
\begin{enumerate}[label=(\alph*)]
    \item  When $\sigma =\alpha$, as $p \in \tilde{S}$ converges to $x \in \partial\Gamma$,
            \begin{enumerate}[label=(\roman*)]
                \item $T_p\to \xi^2(x)$
                \item $N_p\to \xi^2(x)$
                \item $\mathscr{L}_p \oplus T_p\oplus B_p\to \xi^5(x)$
                \item $\mathscr{L}_p\oplus N_p\oplus B_p\to \xi^5(x)$
    \end{enumerate}
    \item When $\sigma = \beta$, as $p \in \tilde{S}$ converges to $x \in \partial\Gamma$, 
            \begin{enumerate}[label=(\roman*)]
                    \item $\mathscr{L}_p\to \xi^1(x)$
                    \item $\mathscr{L}_p\oplus B_p\to \xi^3(x)$
                    \item $T_p\oplus N_p\to \xi^4(x)$
                    \item $T_p\oplus N_p \oplus B_p\to \xi^6(x)$.
   \end{enumerate}
\end{enumerate}
\end{mainthm}

\subsection{Related work}\label{Sec:RelatedWork}

Maximal representations in Hermitian type Lie groups, Hitchin representation in split real Lie groups, and more generally $\Theta$-positive representations, are known to be Anosov \cite{BILW05,Lab06,GLW26}. However, the arguments to certify the Anosov condition in these cases mostly use \emph{positivity} properties, and do not directly involve the non-abelian Hodge correspondence.

The first argument directly using Higgs bundles to show that representations are Anosov is due to Filip, for slices of non-maximal Higgs bundles in $\SO(2,3)$ \cite{Fil23}. These Higgs bundles were Hodge bundles, namely cyclic Higgs bundles of further restricted shape. Filip's result was extended to a larger slice of cyclic Higgs bundles, which includes non-Hodge bundles, by Zhang \cite{Zha25} using similar techniques. Further extensions of these results can be found in \cite{YaoZhang26,YuZhang26}. However, their  approach is different from the one in this paper, and it is not clear how to apply it to more general Lie groups. 
Our technique can be used to recover Filip's result but it does not seem to work to recover the result of Zhang. We believe that this is related to the fact that in $\Gtwosplit$ our argument does not seem to prove the Anosov property for non-Hodge $\beta$-bundles. Quite different techniques were also used by Yao-Zhang in \cite{YaoZhang26} to construct Anosov representations that are close to non-Anosov $\iota$-Fuchsian representations.

Bronstein and the first named author proved the Anosov property for slices of representations in $\SL(3,\R)$ in \cite{DB25}. Both that work and the current work utilize a \emph{parallel distribution of planes} to certify the existence of a certain nested multicones. However, there are a few essential differences between the two strategies. First, the multicones in \cite{DB25} involve the geometry of the symmetric space and its visual compactification, which provides complications for computations. The cones we use in this work are simply cones in projective space. In particular, latter cones are much simpler to define and to study. As a testament to the effectiveness of the new strategy, the calculations to verify nestedness of the cones are easier presently in $\Gtwosplit$, despite the increase in complexity of the ambient Lie group from the $\SL(3,\R)$-case. Moreover, as we alluded to earlier in the introduction, the cones in projective space make sense to define in a much broader context, as they are built directly from the non-vanishing one-form $\phi_f$ that Higgs bundles in the Slodowy slice carry. 
To this end, we note that the parallel work \cite{BDL26} by the Bronstein, Li, and the first named author uses the same techniques as in Theorem \ref{MainThm:C} to prove the Anosov property for families of cyclic Higgs bundles in more general Lie groups. 
However, \cite{BDL26} does not obtain analogues of Theorems \ref{MainThm:B}, \ref{MainThm:D}, \ref{MainThm:E}, namely the use of control on the Anosov limit maps to discuss the relation with the associated geometric structures, and consequences for the space of Anosov representations.

\subsection{Acknowledgments}
We thank Brian Collier for interesting and helpful discussions. C. Davalo was funded by the European Union via the ERC 101124349 ``GENERATE.'' Views and opinions expressed are however those of the authors only and do not necessarily reflect those of the European Union or the European Research Council Executive Agency. Neither the European Union nor the granting authority can be held responsible for them. 

\subsection{Tool use}
Some tedious matrix algebra calculations, such as determinants of $(7\times 7)$-matrices, were verified with Mathematica and/or Sage. 

For this project, we used AI exclusively for the purpose of literature searches of known mathematics, checking for typos and grammatical errors after completion, and assistance with figures. In particular, the only references obtained in this fashion were \cite{PT07, Fin36} for Remark \ref{Remk:Positivity} and \cite{GP10} for standard facts on transversality. 
As this project began in 2024 alongside \cite{DE26b}, we felt it best to finish without AI proof assistance.

\section{Preliminaries}
In this section, we provide a quick overview of the essential objects needed for the paper. We first introduce some basic notions regarding $\Gtwosplit$, then discuss Anosov representations in $\Gtwosplit$, and end with a recollection of the structure of the $\Gtwosplit$-Higgs bundles of interest. 

\subsection{\texorpdfstring{$\Gtwosplit$}{G2'}-background}

To begin, let us recall some definitions and basic facts related to the Lie group $\Gtwosplit$. We shall be quite brief. The reader is referred to \cite[Section 2.1, Section 2.2]{DE26b} for a comprehensive 
review of $\Gtwosplit$-geometry. Additional background on Lie theory for groups of type $G_2$ can be found in \cite{Fon18}. 
\medskip

\subsubsection{Definitions}
In order to define $\Gtwosplit$, we first define a \emph{cross-product} $\times_{3,4} : \R^{3,4}\times \R^{3,4} \rightarrow \R^{3,4}$, which is an antisymmetric bilinear map on pseudo-Euclidean space $\R^{3,4}=(\R^7,q)$, where $q$ is a quadratic form of signature $(3,4)$. As we will need to do computations with this cross product, we recall its explicit definition. However, this definition necessitates a brief detour.  
Indeed, $\times_{3,4}$ arises most naturally from the non-associative unital $\R$-algebra $\Oct'$, the \emph{split octonions}. 

\begin{table}[ht]
    \centering
    \[
\begin{array}{|c|c|c|c|c|c|c|c|c|}
\hline
 col\odot row& 1 & \mathbf{i} & \mathbf{j} & \mathbf{k} & \mathbf{l} & \mathbf{l}\mathbf{i} & \mathbf{l}\mathbf{j} & \mathbf{l}\mathbf{k} \\ \hline
1 & 1 & \mathbf{i} & \mathbf{j} & \mathbf{k} & \mathbf{l} & \mathbf{l}\mathbf{i} & \mathbf{l}\mathbf{j} & \mathbf{l}\mathbf{k} \\[4pt] \hline
\mathbf{i} & \mathbf{i} & -1 & \mathbf{k} & -\mathbf{j} &- \mathbf{l}\mathbf{i} & \mathbf{l} & -\mathbf{l}\mathbf{k} & \mathbf{l}\mathbf{j} \\[4pt] \hline 
\mathbf{j} & \mathbf{j} & -\mathbf{k} & -1 & \mathbf{i} & -\mathbf{l}\mathbf{j} & \mathbf{l}\mathbf{k} & \mathbf{l} & -\mathbf{l}\mathbf{i} \\[4pt] \hline
\mathbf{k} & \mathbf{k} & \mathbf{j} & -\mathbf{i} & -1 & -\mathbf{l}\mathbf{k} & -\mathbf{l}\mathbf{j} & \mathbf{l}\mathbf{i} & \mathbf{l} \\[4pt] \hline
\mathbf{l} & \mathbf{l} & \mathbf{l}\mathbf{i} & \mathbf{l}\mathbf{j} & \mathbf{l}\mathbf{k} & 1 & \mathbf{i} & \mathbf{j} & \mathbf{k} \\[4pt] \hline
\mathbf{l}\mathbf{i} & \mathbf{l}\mathbf{i} & -\mathbf{l} & -\mathbf{l}\mathbf{k} & \mathbf{l}\mathbf{j} & -\mathbf{i} & 1 & \mathbf{k} & -\mathbf{j} \\[4pt] \hline 
\mathbf{l}\mathbf{j} & \mathbf{l}\mathbf{j} & \mathbf{l}\mathbf{k} & -\mathbf{l} & -\mathbf{l}\mathbf{i} & -\mathbf{j} & -\mathbf{k} & 1 & \mathbf{i} \\[4pt] \hline
\mathbf{l}\mathbf{k} & \mathbf{l}\mathbf{k} & -\mathbf{l}\mathbf{j} & \mathbf{l}\mathbf{i} & -\mathbf{l} & -\mathbf{k} & \mathbf{j} & -\mathbf{i} & 1 \\ \hline 
\end{array}
\]
\caption{The multiplication table for $\Oct'$ in the multiplication basis $\mathcal{M}$. }
\label{Table:Oct'Mult}
\end{table}
We now introduce this algebra $\Oct'$, and use it to build the cross product $\times_{3,4}$. To begin, $\Oct'$ admits a vector space basis $\mathcal{M} =(1,\mathbf{i}, \mathbf{j}, \mathbf{k}, \mathbf{l}, \mathbf{li},\mathbf{lj},\mathbf{lk})$, that satisfies Multiplication Table \ref{Table:Oct'Mult} for the algebra product $\odot$. We may refer to $\mathcal{M}$ as the \emph{multiplication basis} for $\Oct'$.  
In what follows, we may write $uv =u\odot v$ to denote the algebra product simply by juxtaposition. The algebra $\Oct'$ is equipped with a linear involution called (split octonion) \emph{conjugation} $*:\Oct'\rightarrow \Oct'$ such that $1^*=1$, and $\mathbf{e}^*=-\mathbf{e}$ for all other basis elements $\mathbf{e}\in\mathcal{M}$. The involution $*$ is anti-multiplicative: $(xy)^* = y^*x^*$. Given $u\in \Oct'$, define $\text{Re}(u)$ and $\text{Im}(u)$ be its projection to the $+1$ and $-1$ eigenspace of the conjugation. The $-1$-eigenspace of this conjugation forms the $7$-dimensional vector subspace $\imoct$, the \emph{imaginary split octonions}, thought of analogously to the imaginary complex numbers or imaginary quaternions.

Let $q$ be the bilinear form on $\Oct'$ defined by $q(u,v)=\text{Re}(uv^*)$, which is of split signature $(4,4)$. In fact, the induced quadratic form satisfies the remarkable multiplicative property 
\[ q(xy)=q(x)q(y),\] 
which promotes $(\Oct', q)$ to the elusive type of a \emph{composition algebra}. By Hurwitz's theorem, the only such $\R$-algebras are $\R,\C, \Ha, \Oct$, and their split counterparts $\C',\Ha', \Oct'$ (see \cite{Fon18}).  

Since $q(1)=1$, we can identify $(\imoct, q)$ with $\R^{3,4}$ as vector spaces with quadratic forms. Note that $\imoct$ is not a subalgebra, for example $\mathbf{i}^2=-1$. However, $\imoct$ carries a binary operation, namely 
the \emph{cross-product $\times:\imoct \times \imoct \rightarrow \imoct$}, given as follows: 
\[ u \times v := \mathrm{Im}(u\odot v) = u\odot v- q(u,v)1_{\Oct'}\ .\]
We may write $\times_{3,4}$ for the cross-product obtained on $\R^{3,4}$ constructed from $\times$ using the identification of $\R^{3,4}$ with $\imoct$.

The Lie group $\Gtwosplit\subset \GL(7,\R)$ can be characterized as the subgroup of elements preserving the cross product $\times_{3,4}$. However, in fact $\Gtwosplit$ preserves $q=q_{3,4}$ and the standard volume form on $\R^7$, so that as a connected Lie group, $\Gtwosplit <\SO_0(3,4)$, where $G_0$ denote the connected component of the identity of $G$.  
The group $\Gtwosplit$ is $14$-dimensional, and is the split real (adjoint) form of the complex group $\Gtwo^{\C}$ associated with the complex simple Lie algebra $\g_2^{\C}$.
Analogously, one can define the cross product on $\C^7\cong \imoct^{\C}$ by complex bilinear extension, and then $\Gtwo^{\C}=\Aut(\C^7, \times_{\C^7})$.

\subsubsection{Flag manifolds and parabolics} We now describe the flag manifolds of $\Gtwosplit$. To this end, we first recall a way to realize Cartan subalgebras $\mathfrak{a}$ in the Lie algebra $\g_2'$ as well as a description of the simple roots $\Delta = \{\alpha, \beta\}$ in the corresponding real root system of $(\g_2',\mathfrak{a})$.

We shall denote $\Gtwo^{\F}$ and $\g_2^{\F}$ for the objects $(\Gtwosplit, \g_2')$ or $(\Gtwo^{\C}, \g_2^{\C})$ as appropriate. Now, the Lie algebra $\g_2^{\F}$ can be described as space of \emph{derivations} of the cross-product on $\imoct ^\F = \imoct\otimes_{\R} \F$: 
\[ \g_2^{\F} = \Der(\imoct^\F) \ = \{\psi \in \mathfrak{gl}(\imoct^{\F}) \mid \psi(x\times y) = \psi(x) \times y+ x\times \psi(y) \}.\]

To describe Cartan subalgebras $\mathfrak{a}^{\F}$ in $\g_2^{\F}$, we
consider the following auxiliary notion from \cite{Eva24a}. 

\begin{definition}
Let $X= (x_k)_{k=3}^{-3}$ be a vector space basis for $\imoct^\F$. We call $X$ an $\F$-cross-product basis when $x_k \times x_l=c_{k,l}x_{k+l}$ for some constants $c_{k,l} \in \F$. 
\end{definition}

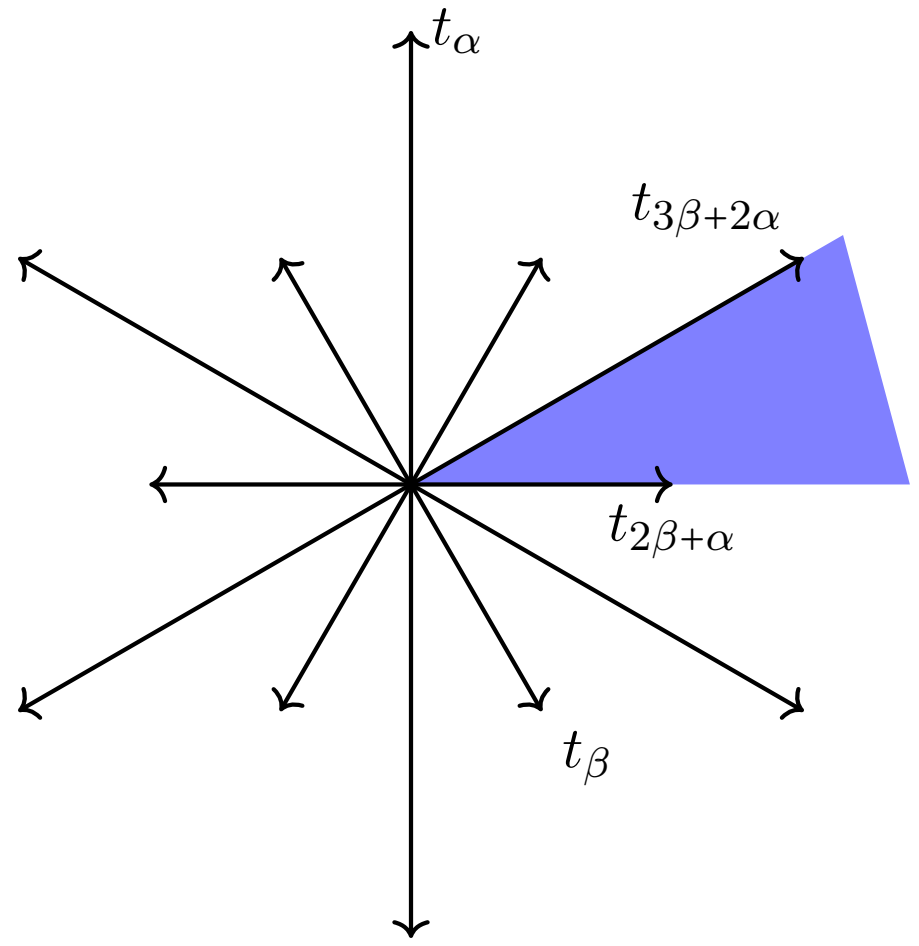
\begin{figure}[ht]
    \centering
\begin{tikzpicture}
    \fill[blue!50] (0,0) -- (0:3.3cm) -- (30:3.3cm) -- cycle;
    \foreach\ang in {0,60,120,...,300}{
     \draw[->,thick] (0,0) -- (\ang:1.731cm);
    }
    \foreach\ang in {30,90,...,330}{
     \draw[->,thick] (0,0) -- (\ang:3cm);
    }
    \node[right] at (90:3cm) {$t_{\alpha}$};
    \node[below right] at (300:1.731cm) {$t_{\beta}$};
    \node[below] at (0:1.731cm) {$t_{2\beta+\alpha}$};
    \node[above left] at (30:3cm) {$t_{3\beta+2\alpha}$};

  \end{tikzpicture}
\caption{The $G_2$ root system, including coroots $t_{\alpha},t_{\beta}$ of the simple roots, and corresponding Weyl chamber.}
\label{Fig:G2RootSystem}
\end{figure}

One can easily check that the subgroup $A^{\mathbb{F}} < \Gtwo^{\F}$ of all diagonal transformations in a given $\F$-cross-product basis $X=(x_k)_{k=3}^{-3}$ obtains the following form: 
\[ A^{\F} = \{ \diag(\lambda \mu, \mu, \lambda, 1, \lambda^{-1}, \mu^{-1}, \lambda^{-1}\mu^{-1} ) \in \GL(\imoct^\F) \mid \lambda, \mu \in \F^*\} \ .\]
The Lie algebra $\mathfrak{a}^{\F} < \g_2^{\F}$ of $A^{\F}$ is a \emph{Cartan subalgebra} given by
\[ \mathfrak{a} = \{ \diag(r+s, r, s, 0, -s, -r,-r-s) \in \gl(\imoct^{\F}) \mid r,s \in \F\}.\]
One can take the roots $\beta: = s^* \in \mathfrak{a}^*$ and $\alpha: =r^*-s^* \in \mathfrak{a}^*$ as a base for the $\F$-root system. Here, among simple roots, $\beta$ is the short root and $\alpha$ the long root. See Figure \ref{Fig:G2RootSystem}, where a corresponding Weyl chamber is shown, along with the corresponding coroots $t_{\alpha}, t_{\beta}$. 

Now, one can form parabolic subalgebras $\mathfrak{p}_{\sigma} < \g_2'$, for $\sigma \in \{\alpha, \beta, \Delta\}$, as the sum of all root spaces of non-negative $\sigma$-height. The parabolic subgroup $P_{\sigma}$ is then found as the normalizer in $\Gtwosplit$ of $\mathfrak{p}_{\sigma}$. 
The corresponding homogeneous spaces $\Gtwosplit/P_{\sigma}$ are identified with the following model spaces: 
\begin{align}
    \Gtwosplit/P_{\beta} \cong \Ein^{2,3} &= \{ [x] \in \mathbb{P}(\R^{3,4}) \mid q(x) = 0 \},\\ 
    \Gtwosplit/P_{\alpha} \cong \Pho^\times  &= \{ \omega \in \Pho(\R^{3,4}) \mid  \omega \times_{3,4} \omega = 0\}, \\ 
    \Gtwosplit/P_{\Delta}\cong \mathcal{F}_{1,2}^\times &= \{  \,(\ell, \omega) \in \Ein^{2,3} \times \Pho^\times \mid \ell \subset \omega \}.
\end{align}
Here, $\Ein^{2,3}$ is the \emph{Einstein universe} of isotropic lines in $\R^{3,4}$, $\Pho^\times$ is the space of \emph{annihilator photons}, namely 2-planes with trivial cross-product, and $\mathcal{F}_{1,2}^\times$, the full flag manifold, consisting of all \emph{pointed annihilator photons}. We recall further facts about $\Ein^{2,3}$, $\Pho^\times$ as needed later. 
\begin{remark}
Each of the $\Gtwosplit$-flag manifolds is self-opposite, and the equivariant inclusion maps $(\Gtwosplit, \Ein^{2,3}) \hookrightarrow (\SO(3,4), \Ein^{2,3})$ and $(\Gtwosplit, \Pho^\times)\hookrightarrow (\SO(3,4), \Pho(\R^{3,4}))$ are transversality-preserving.
\end{remark}

One important notion related to the $\Gtwosplit$-flag manifolds is as follows.
\begin{definition}
Let $u \in \imoct$. Define the \textbf{annihilator} of $u$ by 
\[ \Ann(u) = \{ v \in \imoct \mid u\times v=0\}.\]
\end{definition}

When $q(u)\neq0$, then $\Ann(u) = \R\{u\}$, as naively expected. However, if $u$ is isotropic, its annihilator is nontrivial. 
\begin{proposition}\label{Prop:Annihilators}
Suppose $u \in \imoct$ has $q(u)=0$. Then $\Ann(u)$ is 3-dimensional and isotropic. 
\end{proposition}

\begin{remark}[line to 3-plane]
In fact, the equivariant inclusion $\Ann:(\Gtwosplit,\Ein^{2,3})\hookrightarrow (\SO(3,4),\Iso_3(\R^{3,4}))$ 
by $\ell \mapsto \Ann(\ell)$ is transversality-preserving \cite[Section 3.1]{ER25}. 
\end{remark}

In the study of the $\Gtwosplit$-flag manifolds, we may use the following space. 
\begin{remark}
The Riemannian symmetric space $\X_{\Gtwosplit}\cong \Gtwosplit/K$ can be identified with the model space $\{ U \in \Gr_{(3,0)}^{\times}(\R^{3,4}) \mid U \times_{3,4} U = U\}$ of spacelike 3-planes closed under cross-product. 
For a proof, see \cite[Lemma 3.5]{Eva24a}.  
In this way, may write $U \in \X_{\Gtwosplit}$ to refer to such a 3-plane. 
\end{remark}

We may also use the following notation to define certain quadric hypersurfaces.
\begin{remark}[Quadrics]
For a vector space $V$ equipped with a quadratic form $q_V$, 
we write 
$Q_+(V) = \{ x \in V \mid q_V(x)=+1\}$ and $Q_-(V) = \{x \in V \mid q_V(x)=-1\}$ for the positive and negative quadrics of $V$. 
\end{remark}

Now let us formally define $\Gtwosplit$-Frenet splittings, which will appear frequently going forward. 

\begin{definition}
\label{defn:FrenetFrame}
 A \textbf{Frenet splitting} is the data of a quadruple $(\mathscr{L},T,N,B)$ of a line and three $(2-)$planes in $\R^{3,4}$ such that:
\begin{itemize}
    \item $\mathscr{L}\oplus T\oplus N \oplus B$ is an orthogonal splitting of $\R^{3,4}$,
    \item The planes $T$, $N$, $B$ are invariant under cross product by the line $\mathscr{L}$,
    \item $\mathscr{L}$ and $N$ are spacelike, $T$ and $B$ are timelike. 
\end{itemize}

An \textbf{oriented Frenet splitting} is the data of quadruple $(e,T,N,B)$ of a vector and three planes in $\R^{3,4}$ such that $e$ is spacelike of norm $1$ and $(\R\{ e\},T,N,B)$ is a Frenet frame.
\end{definition}

A \emph{Frenet frame} is a map $p\in \widetilde{S}\mapsto (\mathscr{L}_p,T_p,N_p,B_p)$ from a surface into the space of Frenet splittings, that defines a $J$-holomorphic curve $(\mathscr{L}_p)_{p\in \widetilde{S}}$ in the pseudosphere $\mathbb{S}^{2,4}\cong \mathbb{P}Q_+(\R^{3,4})$, in the sense of \cite{CT24}, with tangent planes $(T_p)$, normal planes $(N_p)$ (i.e. the image of the second fundamental form, away from the zeros of $\sff$), and binormal planes $(B_p)$. Such a map is canonically associated to $\beta$-bundles, as we recall in \S\ref{Subsec:SigmaBundles}. 

\begin{remark}
    In fact, an $\alpha$-bundle with $\beta \neq 0$ also defines an equivariant Frenet frame. However, the map $p\mapsto \mathscr{L}_p$ in this case is $J$-holomorphic with branched points exactly at the zeros of the differential $\beta$. One may call such a map a \textbf{branched Frenet frame}. In fact, these branched Frenet frames associated to $\alpha$-bundles have the property that their second fundamental forms never vanish. 
\end{remark}

The developing map of the geometric structure built in \cite{DE26b} is constructed in terms of this (branched) Frenet frame both for $\alpha$-bundles and $\beta$-bundles. We recall the relevant details in \S\ref{Sec:Comparison}.

\subsection{Anosov representations}\label{Subsec:AnosovBackground}
Next, we recall a few characterizations of Anosov representations in $\SL(n,\R)$ and then in $\Gtwosplit$.\medskip

Let $\Gamma$ be any hyperbolic group. Let us fix a word metric $|\cdot|$ on $\Gamma$, determined by a finite generating set $F\subset \Gamma$. Let $\lVert [\gamma] \rVert=\min_{\alpha\in \Gamma} |\alpha\gamma\alpha^{-1}|$, which depends only on the conjugacy class of $\gamma$. The \emph{Gromov boundary} $\partial\Gamma$ is the visual boundary of the Cayley graph of $\Gamma$ endowed with the metric $|\cdot|$; the quasi-isometric structure of $\partial\Gamma$ is independent of $|\cdot|$. Note that $\partial\Gamma$ carries a canonical $\Gamma$-action. 

Now, suppose $\Gamma=\pi_1S$, where $S$ a closed oriented surface of genus at least two. If we fix any Riemannian metric on $S$, then there exist $A,B>0$ such that the length $\ell([\gamma])$ of the shortest geodesic representative of the free homotopy class of $\gamma$ on $S$ satisfies:
\[A\lVert [\gamma] \rVert +B  \geq \ell([\gamma]) \geq \frac{1}{A}\lVert [\gamma]\rVert -B.\] 
In particular, note that $\ell([\cdot])$ is conjugation-invariant by definition.
\medskip

Let $\boldsymbol{\lambda}: \SL(n,\R) \rightarrow \sllie(n,\R)$ be the \emph{Jordan projection}, given by the logarithm of the moduli of the eigenvalues of $g$ in non-increasing order:
\begin{align}\label{eq:JordanProj}
    \boldsymbol{\lambda}(g) \coloneqq \diag(\lambda_1(g), \lambda_2(g), \cdots , \lambda_{n}(g) ). 
\end{align}
If we fix a scalar product on $\R^{n}$, one can also define the \emph{Cartan projection} $\boldsymbol{\sigma}: \SL(n,\R) \rightarrow \sllie(n,\R)$, given by the logarithm of the singular values of $g$ in non-increasing order:
\begin{align}\label{eq:CartanProj}
    \boldsymbol{\sigma}(g) \coloneqq \diag(\sigma_1(g), \sigma_2(g), \cdots , \sigma_{n}(g) ).
\end{align}
Let us now define Anosov representations into $\SL(n,\R)$, motivated by \cite{KLP17, BPS}.
To this end, we denote $\alpha_k\big( \diag(\lambda_1,\cdots, \lambda_{n}) \big)=\lambda_k-\lambda_{k+1}.$
\begin{definition}[Anosov property via singular value gaps]\label{defn:Anosov}
Let $\Gamma$ be a word hyperbolic group and $k$ be an integer satisfying $1\leq k \leq n-1$. A representation $\rho:\Gamma\to \SL(n,\R)$ is $\boldsymbol{P_{k}}$\textbf{-Anosov} if there exists $c,D>0$ such that for all $\gamma\in \Gamma$:
      \[\alpha_k(\boldsymbol{\sigma}(\rho(\gamma)))=\sigma_k(\rho(\gamma))-\sigma_{k+1}(\rho(\gamma))\geq c|\gamma|-D.\]
\end{definition}

This definition is independent of the chosen word metric, as any two word metrics are coarsely comparable. An equivalent characterization via eigenvalues from \cite{KP22} is the following: 

\begin{proposition}[Anosov property via eigenvalue gaps]\label{Prop:Anosov}
Let $\Gamma$ be a word hyperbolic group and $k$ be an integer satisfying $1\leq k \leq n-1$. A representation $\rho:\Gamma\to \GL(n,\R)$ is $P_{k}$-Anosov if and only if there exists $c,D>0$ such that for all $\gamma\in \Gamma$:
      \[ \alpha_k(\boldsymbol{\lambda}(\rho(\gamma)))= \lambda_k(\rho(\gamma))-\lambda_{k+1}(\rho(\gamma))\geq c\lVert [\gamma]\rVert -D.\]
\end{proposition}

\begin{remark}[Word length vs. geometric length]
Note that in this characterization if $\Gamma=\pi_1S$ is a closed surface group, then one can replace $\lVert [\gamma]\rVert$ by $\ell([\gamma])$ for any fixed Riemannian metric $g$ on $S$, where $\ell([\gamma])$ is the length of the shortest curve freely homotopic to $\gamma$, up to changing the constants $c,D$. Indeed, this holds for hyperbolic metrics, and every length metric $\ell([\cdot])$ on the compact surface is comparable to that of a fixed hyperbolic metric.
\end{remark}

Anosov representations in particular admit \emph{limit maps} \cite{GGKW17}. Recall that for the hyperbolic group $\Gamma=\pi_1(S_g)$ with $g\geq2$, the Gromov boundary $\partial \Gamma$ is topologically a circle. 

\begin{proposition}[Limit maps]
\label{Prop:G2LimitMap}
Let $\rho:\Gamma\to \SL(n,\R)$ be a $P_k$-Anosov representation of a hyperbolic group $\Gamma$ with Gromov boundary $\partial \Gamma$. Then there exists a unique $\rho$-equivariant continuous map $\xi^k:\partial \Gamma \to \Gr_k(\R^{n})$ satisfying the following: for any $\gamma \in \Gamma$, if $\gamma^+\in \partial \Gamma$ is the attracting fixed point of $\gamma$, then $\xi^k(\gamma^+)$ is the unique attracting fixed point of $\rho(\gamma)$ on $\Gr_k(\R^{n})$.
\end{proposition}

For any $g\in \Gtwosplit$, the Jordan projection $\lambda(g) \in \g_2'$ lies in the Weyl chamber 
\[ \overline{\mathfrak{a}}^+ = \{ \diag(r+s, r, s, 0,-s,-r,-r-s) \in \g_2' \mid r\geq s \geq 0\}. \]
The two simple roots $\alpha, \beta$ of $\Gtwosplit$ satisfy the following, for any $g\in \Gtwosplit$:
\begin{align}
    \beta(\boldsymbol{\lambda}(g))&= \alpha_k(\boldsymbol{\lambda}(g)), \qquad k \in \{1,3,4,6\}, \label{eq:BetaEqualities}  \\
    \alpha(\boldsymbol{\lambda}(g))&= \alpha_k(\boldsymbol{\lambda}(g)), \qquad k \in \{2,5\}. \label{eq:AlphaEqualities}
\end{align} 
Note that $\lambda_4(g)=0$ for all $g\in \Gtwosplit$.

Now, let $\rho:\Gamma\to \Gtwosplit$ be a representation. 
As a consequence of \eqref{eq:BetaEqualities}, \eqref{eq:AlphaEqualities}, 
for the inclusion $\iota \circ \rho: \Gamma \rightarrow \SL(7,\R)$, the $P_k$-Anosov properties for $k\in\lbrace 1,3,4,6\rbrace$ are equivalent, and the $P_k$-Anosov properties for $k\in\lbrace 2,5\rbrace$ are  equivalent. Hence, by \cite[Proposition 4.1, Proposition 4.4]{GW12}, which relates the Anosov properties in $\Gtwosplit$ and those in $\SL(7,\R)$, we can take the following definition. 

\begin{definition}\label{Defn:AnosovG2'}
Let $\Gamma$ be a hyperbolic group and $\rho:\Gamma \to \Gtwosplit$ be a representation. 
\begin{itemize}
    \item If $\iota\circ \rho$ is $P_{k}$-Anosov for any $k \in \{1,3,4,6\}$, then it is called $\boldsymbol{P_\beta}$\textbf{-Anosov}.
    \item If $\iota\circ\rho$ is $P_{k}$-Anosov for any $k \in \{2,5\}$, then it is called $\boldsymbol{P_\alpha}$\textbf{-Anosov}.
\end{itemize}
\end{definition}

The Anosov limit maps of $\Gtwosplit$-representations land in the appropriate $\Gtwosplit$-flag manifolds. 

\begin{proposition}[{\cite[Proposition 4.4]{GW12}}]
\label{prop:LimitMapPropertiesG2}
Let $\Gamma$ be a hyperbolic group and $\rho:\Gamma\to \Gtwosplit$ be a representation.
\begin{itemize}
    \item If $\rho$ is $P_{\beta}$-Anosov, then $\xi^1(\partial \Gamma)\subset \Ein^{2,3}$. Also, for all $x\in \partial \Gamma$, we have 
    \[ \big[\xi^3(x)=\Ann(\xi^1(x))\big] \subset \big[\xi^4(x)=\xi^3(x)^{\perp_q} \big] \subset \big[ \xi^6(x)=\xi^1(x)^{\perp_q} \big] \] 
    In particular, $\xi^3:\partial \Gamma \rightarrow \Iso_3(\R^{3,4})$. 
    \item If $\rho$ is $P_\alpha$-Anosov, then $\xi^2(\partial \Gamma)\subset \Pho^\times$ and for all $x \in \partial\Gamma$, we have $\xi^5(x)=\xi^2(x)^{\perp_q}$.
\end{itemize}
\end{proposition}

\subsection{\texorpdfstring{$\alpha$}{alpha}-bundles and \texorpdfstring{$\beta$}{beta}-bundles}\label{Subsec:SigmaBundles}

In this subsection, we recall the notion of $\alpha$-bundles and $\beta$-bundles from \cite{DE26b} and some of their relevant features. These Higgs bundles are realized by imposing certain restrictions on a general class of cyclic $\Gtwosplit$-Higgs bundles as follows. 

\subsubsection{Definitions} We begin with a rank seven holomorphic vector bundle $\V = \bigoplus_{k=3}^{-3} \mathcal{L}_k$ that is a direct sum of holomorphic line subbundles with once-and-for-all fixed isomorphisms 
\begin{align}\label{FixedProductIsomorphism}
    \begin{cases}
        \mathcal{L}_0= \mathcal{O}\\
        \mathcal{L}_{-k}\cong \mathcal{L}_{k}^{-1}, \; k \in \{1,2,3\} \\
    \mathcal{L}_1 \mathcal{L}_2 \cong \mathcal{L}_{3}.
    \end{cases}
\end{align}
As explained in \cite[Section 3.2]{DE26b}, we can place a holomorphic cross-product $\times_{\V}:\V \otimes \V \rightarrow \V$ on the vector bundle $\V$. See \cite[Proposition 3.7]{DE26b} in particular. This procedure involves also defining a fixed non-degenerate holomorphic bilinear form $Q$ on $\V$ that is compatible with $\times_{\V}$. We shall uniformly apply the same conventions for $(Q, \times_{\V})$ from \cite{DE26b}. Thus, we notationally suppress these objects. We refer to \cite{DE26b} for more comprehensive details on $\Gtwosplit$-Higgs bundles. 

Now, here is the definition of interest. Recall that $\mathcal{K}=\left(\T^{(1,0)}\Sigma\right)^*$ is the canonical holomorphic line bundle on $\Sigma$.
\begin{definition}\label{Defn:G2Cyclic}
A \textbf{cyclic} $\Gtwosplit$-Higgs bundle over a Riemann surface $\Sigma$ is a pair $(\V,\Phi)$ such that $\V = \bigoplus_{k=3}^{-3} \mathcal{L}_k$ splits holomorphically as a sum of seven holomorphic line bundles satisfying \eqref{FixedProductIsomorphism}, and the Higgs field $\Phi\in H^0(\End(\V)\otimes \mathcal{K})$ is of the form 
\begin{equation}\label{CyclicG2'Higgs}
\begin{tikzcd}
 \mathcal{L}_3 \arrow[r, "\beta"] & \mathcal{L}_2\arrow[r, "\alpha"]& \mathcal{L}_1 \arrow[r, "-i\sqrt{2}\beta"]& \mathcal{L}_0 \arrow[r, "-i\sqrt{2}\beta"]&\mathcal{L}_{-1}\arrow[r, "\alpha"]&\mathcal{L}_{-2} \arrow[r, "\beta"]\arrow[lllll,bend left, "\delta"]&\mathcal{L}_{-3} \arrow[bend left, lllll, "\delta"]
\end{tikzcd}.
\end{equation}
Furthermore: 
\begin{itemize}[noitemsep]
    \item When $\mathcal{L}_{-1} = \K^{-1}$ and $\beta = 1$, we call the bundle $\boldsymbol{\beta}$-\textbf{cyclic}.
    \item When $\mathcal{L}_{3}\cong \mathcal{L}_2\K$ and $\alpha = 1$, we call the bundle $\boldsymbol{\alpha}$-\textbf{cyclic}.
\end{itemize} 
\end{definition}
In \cite{DE26b}, we explain why the Higgs field $\Phi$ in \eqref{CyclicG2'Higgs} respects the cross-product, i.e. $\Phi \in \Der(\times_{\V})$.  

We now consider an additional condition for these cyclic Higgs bundles. 
\begin{remark}
A $\beta$-cyclic bundle with either $\delta=0$ or $\alpha=0$ and similarly an $\alpha$-cyclic bundle with 
$\beta=0$ or $\delta=0$ is a \textbf{Hodge bundle}, namely a $\C^*$-fixed point in the moduli space of Higgs bundles. We refer to the former as a $\boldsymbol{\beta}$\textbf{-Hodge bundle} and the latter as an $\boldsymbol{\alpha}$\textbf{-Hodge bundle}. 
\end{remark}

\subsubsection{Properties of $\sigma$-bundles}
Now, to force conditions \eqref{FixedProductIsomorphism}, one easily observes $(\mathcal{L}_i)_{i=3}^{-3}$ must obtain the form: 
\[(\mathcal{L}_1\mathcal{L}_2, \;\mathcal{L}_2,\; \mathcal{L}_1, \; \mathcal{O}, \;\mathcal{L}_1^*, \;\mathcal{L}_2^*, \; \mathcal{L}_1^*\mathcal{L}_2^* ).\]

Recall that a \emph{harmonic metric} $h$ on a Higgs bundle $(\V,\Phi)$ is a hermitian metric such that the connection $\nabla =\nabla^h+\Phi+\Phi^{*h}$ is flat, where $\nabla^h$ is the Chern connection of $(\V,h)$ and $\Phi^{*h}$ is the adjoint of $\Phi$.
A harmonic metric $h$ on a cyclic $\Gtwosplit$-Higgs bundle has a restricted form.
\begin{proposition}[$\Gtwosplit$-harmonic metric {\cite{CT24}}] \label{Prop:HarmonicMetricDiagonal}
    Any harmonic metric $h$ on a polystable cyclic $\Gtwosplit$-Higgs bundle obtains the form:
    \[ h = \diag(h_i)_{i=3}^{-3}=\diag \left( g_1^{-1}g_2^{-1}, g_2^{-1}, g_1^{-1}, 1, \; g_1, \; g_2, \; g_1g_2 \right).\]

    The real locus $\V^{\R}$ of $\mathcal{E}$ preserved by the flat connection $\nabla=\nabla^h+\Phi+\Phi^{*h}$ is the real vector subbundle of elements $(z_3,z_2,\cdots, z_{-3})$ satisfying $z_i=(z_{-i})^{*h_i}$ for $-3\leq i\leq 3$. 
\end{proposition}
For a Higgs bundle from Definition \ref{Defn:G2Cyclic}, the connection $\nabla$ is a $\Gtwosplit$-connection. In other words, the flat connection $\nabla$ preserves the nondegenerate bilinear form $Q$ on $\V$, the bundle cross-product $\times_{\V}:\V\otimes \V \rightarrow \V$, and a $\C$-anti-linear involution $\lambda: \V\rightarrow \V$ with fixed point set $\V^{\R}$. The real locus $\V^{\R}$ is closed under $\times_{\V}$, and has the property that $(\V^{\R}, Q, \times_{\V})$ fiberwise is isomorphic to $(\R^{3,4},\times_{3,4})$.  

This real locus $\V^{\R}$ further decomposes into a bundle Frenet splitting, which fiberwise satisfies Definition \ref{defn:FrenetFrame}.
Indeed, $\lambda$ preserves the line bundle $\mathcal{L}_0$ and each of the rank two bundles $\mathcal{L}_i\oplus \mathcal{L}_{-i}$ for $i \in \{1,2,3\}$. We may then write $\mathscr{L},T,N,B$ for the following sub-bundles of $\V^\R$:
\begin{align}\label{LTNB}
    \begin{cases}
    \mathscr{L} &= \Fix(\lambda|_{\mathcal{L}_0})\\
    T &= \Fix(\lambda|_{\mathcal{L}_1\oplus \mathcal{L}_{-1}}) \\
    N&= \Fix(\lambda|_{\mathcal{L}_2\oplus \mathcal{L}_{-2}}) \\
    B &= \Fix(\lambda|_{\mathcal{L}_3\oplus \mathcal{L}_{-3}}).\end{cases}
\end{align} 
We will frequently appeal to this Frenet splitting $\V^{\R}=\mathscr{L}\oplus T \oplus N \oplus B$. 

As in \cite{DE26b}, we may use certain convenient local frames for $\Gtwosplit$-Higgs bundles. 
\begin{remark}\label{Remk:UnitaryCrossProd}
For any $p \in \Sigma$, or any open set $U \subset \Sigma$, we can find a basis for $\V|_p$ or a local frame for $\V|_U$, respectively, of the form $(x_k)_{k=3}^{-3}$ that is an $h$-unitary $\C$-cross product basis. See \cite[Proposition 3.7, Remark 3.12]{DE26b}. 
\end{remark}

We shall need the following general form of Hitchin's equations in terms of the holomorphic differentials $\alpha,\beta,\delta$. Here, we define the curvature function $\kappa_g: S \rightarrow \R$ as in \cite[Section 3.2]{DE26b}. 

\begin{lemma}[Hitchin's equations {\cite[Lemma 3.16]{DE26b}}]\label{Lem:HitchinsEquationsGeneral}
Let $g$ be a conformal metric on $\Sigma$ and $(\V,\Phi)$ a cyclic $\Gtwosplit$-Higgs bundle on $\Sigma$. Suppose that $h$ is a harmonic metric on $(\V,\Phi)$. Then $\alpha, \beta, \delta$ satisfy the following equations away from their zero respective sets: \begin{align}\label{HitchinEquationsDifferentials}
    \begin{cases}
   \Delta_g \log\lVert \alpha \rVert^2&= 2\lVert \alpha \rVert ^2-3\lVert \beta \rVert^2-\lVert \delta \rVert ^2 +\kappa_g,\\
 \Delta_g \log\lVert \beta \rVert^2&= 2\lVert \beta\rVert ^2-\lVert \alpha\rVert^2 + \kappa_g,\\
 \Delta_g \log\lVert \delta \rVert^2&=2\lVert \delta \rVert ^2-\lVert  \alpha\rVert ^2 + \kappa_g.
\end{cases}
\end{align}
\end{lemma}
\begin{remark}
Given a conformal metric $g$ on $\Sigma$, we may write $\lVert \eta\rVert$, for $\eta \in H^0(\End({\V})\otimes \K)$, to denote the function $\lVert \eta\rVert=\lVert \eta\rVert_{h,g}$. That is,
if $\eta = \hat{\eta} \,dz$ locally, then $\lVert\eta\rVert =\lVert \hat{\eta}\rVert_h\lVert dz\rVert _g$. 
\end{remark}
Let us recall the relevant maximum principles for $\beta$-Hodge bundles from \cite[Lemma 4.11]{DE26b}. 

\begin{lemma}[Maximum principle for $\beta$-Hodge bundles]\label{Lem:MaximumPrincipleBeta}
Let $(\V, \Phi)$ be a $\beta$-cyclic bundle and $h$ be a harmonic metric.
\begin{enumerate}[label=(\roman*), noitemsep]
    \item If $\delta=0$, then we have the following global inequality:
 \[ \lVert \alpha\rVert <  \frac{\sqrt{5}}{\sqrt{3}} \lVert \beta \rVert . \]

 \item If $\alpha=0$ and $\delta$ has at least one zero, then we have the following global inequality:
 \[ \lVert \delta\rVert <   \lVert \beta \rVert. \]
\end{enumerate}
\end{lemma}

For $\alpha$-bundles, we shall need an improved maximum principle from that of \cite{DE26b}, which we prove later in $\S$\ref{Subsec:AlphaNested}. \medskip 

Going forward, we now denote the non-vanishing one-form $\phi_f$ from the introduction by $\Phi_0$, as in \cite{DE26b}. 
From $\Phi_0$, we then build the object $\Psi_0 = \Phi_0+\Phi_0^{*h} \in \Omega^0(\T^*S\otimes \End(\V))$. We shall once again crucially use the parallelism of $\Phi_0$ and $\Psi_0$ from \cite[Section 4.2, Section 5.6]{DE26b}. 

\begin{proposition}[$\Psi_0$-parallelism]\label{Prop:Psi0Parallelism}
Let $(\V,\Phi)$ be a polystable $\sigma$-bundle, $h=\diag(h_i)_{i=3}^{-3}$ be a harmonic metric, and $\nabla^h$ be the Chern connection on $\text{End}(\mathcal{E})$ associated to $h$. 
\begin{itemize}[noitemsep]
    \item When $\sigma = \beta$, set $g := h_{-1}$, a Hermitian metric on $\T\Sigma\simeq \mathcal{K}^{-1}$. 
    \item When $\sigma = \alpha$, set $g :=h_{1}h_{2}^{-1}$, a Hermitian metric on $\T\Sigma\simeq \mathcal{K}^{-1}$.
\end{itemize}
In both cases, let $\nabla^g$ be the Levi-Civita connection of $g$.  
Then for any vector field $X \in \Omega^0(\Sigma, \T\Sigma)$,  
\begin{align}\label{Parallelism_BetaBundles}
    (\nabla^h \circ \Phi_0)(X)=& \,(\Phi_0\circ \nabla^g)(X),\\
    (\nabla^h \circ \Psi_0)(X)=& \,(\Psi_0\circ \nabla^g)(X).
\end{align}
In other words, $\Phi_0 \in \Omega^0(\K\otimes \End(\V))$ and $\Psi_0 \in \Omega^0(\T^*S\otimes \End(\V))$ are parallel. 
\end{proposition}

\section{The Anosov property}\label{Sec:AnosovProperty}

In this section, we prove that many representations associated to $\sigma$-bundles satisfy Anosov properties. 
The Anosov property is verified by a certain inclusion of cones in the (projective) flat bundle of $\rho$ under a geodesic flow that allows us to verify large eigenvalue gaps. We describe the framework for verifying the Anosov property in $\S$\ref{Sec:NestedCones} then implement it via Higgs bundles in $\S$\ref{sec:NestedConesHiggs}.

Here are the main theorems we will prove in this section. 
\begin{theorem}[Anosov property for $\beta$-Hodge bundles]\label{Thm:BetaAnosov}
Let $(\V,\Phi)$ be a polystable $\beta$-cyclic Hodge bundle. These are of one of the following form, corresponding to $\delta=0$ or $\alpha=0$, respectively:
\begin{equation}\label{eq:betaBundlesAnosovSection}
\begin{tikzcd}
\mathcal{B} \arrow[r, "1"] & \mathcal{B}\mathcal{K}^{-1}\arrow[r, "\alpha"]&\mathcal{K}\arrow[r, "-i\sqrt{2}"]&\mathcal{O} \arrow[r, "-i\sqrt{2}"]&\mathcal{K}^{-1} \arrow[r, "\alpha"]&\mathcal{B}^{-1}\mathcal{K} \arrow[r, "1"]&\mathcal{B}^{-1}
\end{tikzcd}.
\end{equation}

\begin{equation}\label{eq:betaBundlesAnosovSectionBis}
\begin{tikzcd}
\mathcal{B} \arrow[r, "1"] & \mathcal{B}\mathcal{K}^{-1}&\mathcal{K}\arrow[r, "-i\sqrt{2}"]&\mathcal{O} \arrow[r, "-i\sqrt{2}"]&\mathcal{K}^{-1} &\mathcal{B}^{-1}\mathcal{K} \arrow[r, "1"]\arrow[bend left, lllll, "\delta"]&\mathcal{B}^{-1}\arrow[bend left, lllll, "\delta"]
\end{tikzcd}.
\end{equation}
Suppose that $\delta$ has at least one zero. 
Then the associated representation $\rho:\pi_1S \rightarrow \Gtwosplit$ is $P_{\beta}$-Anosov. 
\end{theorem}

\begin{remark}\label{Remk:SL3}
    Let $(\V,\Phi)$ be a polystable $\beta$-bundle. 
    In fact, $\delta$ is non-vanishing exactly when $d=0$ (see \cite[Proposition 5.6, Theorem 5.14]{CT24}).
    In this case, the corresponding representations factor through a $\Z_2$-cover of a particular $\SL(3,\R)$-subgroup of $\Gtwosplit$ (see \cite[Remark 3.12]{CT24}). For a representation to this cover of $\SL(3,\R)$ to include into $\Gtwosplit$ as a $P_\beta$-Anosov-representation, its Jordan limit cone (in $\sllie_3\R$) must avoid the middle line of the model Weyl chamber, which is impossible for surface group representations.
\end{remark}

The companion result for $\alpha$-cyclic bundles is as follows. 

\begin{theorem}[Anosov property for $\alpha$-bundles]\label{thm:AlphaAnosov}
Let $(\V,\Phi)$ be a polystable $\alpha$-bundle. These are of the following form:
\begin{equation}\label{eq:alphaBundlesAnosovSection}
\begin{tikzcd}
\mathcal{T}^2\mathcal{K} \arrow[r, "\beta"] & \mathcal{T}\mathcal{K}\arrow[r, "1"]&\mathcal{T}\arrow[r, "-i\sqrt{2}\beta"]&\mathcal{O} \arrow[r, "-i\sqrt{2}\beta"]&\mathcal{T}^{-1} \arrow[r, "1"]&\mathcal{T}^{-1}\mathcal{K}^{-1} \arrow[r, "\beta"]\arrow[bend left, lllll, "\delta"]&\mathcal{T}^{-2}\mathcal{K}^{-1} \arrow[bend left, lllll, "\delta"]
\end{tikzcd}.
\end{equation}
Then the associated representation $\rho:\pi_1S \rightarrow \Gtwosplit$ is $P_{\alpha}$-Anosov. 
\end{theorem}

\subsection{Nested cones and eigenvalue gaps}\label{Sec:NestedCones}
We now present a sufficient geometric criterion to detect Anosov representations. First, given a transformation $g\in \GL(n,\R)$, we relate preservation by $g$ of \emph{cones} in projective space to the eigenvalue gaps of $g$. Next, given an equivariant map from the universal cover $\tilde{S}$ to a space of such cones in $\RP^{n-1}$, such that these cones are nested under a certain geodesic flow, we show the cones converge to the (projectivization of the) Anosov limit set of the corresponding representation. 
Similar ideas were used in \cite{Dav26}, but with a control on singular value gap. The techniques presented in this subsection are also developed in \cite{BDL26} to prove gaps for singular values, but there the control of the Anosov limit map is not exploited. 

Here, by `cone', we mean the following. 
\begin{definition}
We call $C \subset \mathbb{RP}^{n-1}$ a $\boldsymbol{d}$-\textbf{cone} when $C$ is the set of positive lines for some symmetric bilinear form $b$ of signature $(d,n-d)$ on $\R^n$, i.e. $C = \mathbb{P}Q_+(\R^n, b)$. 
\end{definition}
Observe that if $C$ is a $d$-cone associated to bilinear form $b$ of signature $(d, n-d)$, then the complement of the closure $\RP^{n-1}\backslash \overline{C} = \mathbb{P}Q_+(\R^n, -b)$ is an $(n-d)$-cone. 

Going forward, it will be important to quantitatively measure how much cones sit inside one another. The following definition makes this notion precise. 

\begin{definition}[{\cite[{Definition 5.1}]{Dav26}}]
    Let $C_1,C_2$ be $d$-cones such that $\overline{C_2}\subset C_1$. We define the \textbf{cross-ratio nestedness} $\Nest(C_1,C_2)$ as the infimum for all projective line $L \subset \mathbb{RP}^{n-1}$ and all quadruple of points $(x_1,x_2,y_2,y_1)$ on $L$ in this cyclic order such that $x_2\neq y_2\in \partial C_2$ and $x_1\neq y_1\in \partial C_1$ of the cross ratio $[x_1,x_2,y_2,y_1]$.
\end{definition}
The strict inclusion of $C_2$ into $C_1$ implies that this infimum is actually a minimum, see \cite{Dav26}, and hence this minimum is positive. 

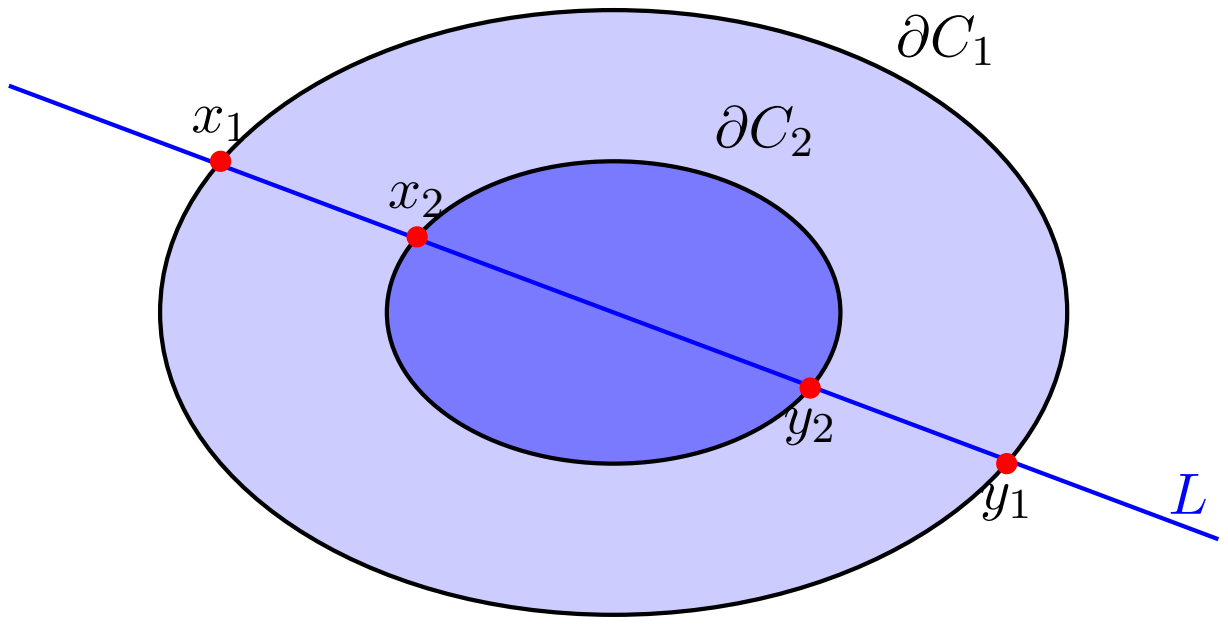
\begin{figure}[h]
    \centering
    \begin{tikzpicture}

    \draw[thick, fill=blue, fill opacity=.2,] (0,0) ellipse (3cm and 2cm);
    \node at (2.2,1.8) {$\partial C_1$};

    \draw[thick, fill=blue, fill opacity=.4] (0,0) ellipse (1.5cm and 1cm);
    \node at (1,1.2) {$\partial C_2$};

    \draw[blue, thick] (-4,1.5) -- (4,-1.5);
    \node[blue] at (3.8,-1.2) {$L$};

    \fill[red] (-2.6, 1.0) circle (2pt) node[above, black] {$x_1$};
    \fill[red] (2.6, -1.00) circle (2pt) node[below, black] {$y_1$};

    \fill[red] (-1.3, 0.5) circle (2pt) node[above, black] {$x_2$};
    \fill[red] (1.3, -0.5) circle (2pt) node[below, black] {$y_2$};

\end{tikzpicture}
    \caption{\emph{Illustration of the nestedness of $1$-cones in $\mathbb{RP}^2$.}}
    \label{fig:Nestedness1}
\end{figure}

Note that the nestedness is $\GL(n,\R)$-invariant since projective lines and cross-ratios are preserved by projective linear transformations. 
\medskip

We use the following convention regarding the cross-ratio: given any 4-tuple $(x_1,x_2,y_2,y_1)$ of points in $\RP^n$ on a projective line, for any affine chart where $x_1,x_2,y_2,y_1$ correspond to the real numbers $X_1, X_2, Y_2, Y_1 \in \R \cup \{\infty\}$, we define 
\[ [x_1,x_2,y_2,y_1]=\frac{Y_2-X_1}{Y_2-X_2}\frac{Y_1-X_2}{Y_1-X_1}. \]
This convention agrees with that of \cite{BGLPW24} and \cite{Dav26}.

\begin{figure}[h!]
    \centering
    \begin{tikzpicture}
        \draw[thick] (-4,0) -- (4,0) ;
        \draw[thick, blue] (-3,.2) -- ++(0,.2) -- ++(6,0) -- ++(0,-.2);
        \draw[thick, blue] (-1,.1) -- ++(0,.2) -- ++(2,0) -- ++(0,-.2);
        \fill[red] (-3,0) circle (2pt) node[below, black] {$X_1$};
        \fill[red] (3,0) circle (2pt) node[below, black] {$Y_1$};
        \fill[red] (-1,0) circle (2pt) node[below, black] {$X_2$};
        \fill[red] (1,0) circle (2pt) node[below, black] {$Y_2$};
    \end{tikzpicture}
    \caption{\emph{Illustration of the nestedness of segments in $\mathbb{RP}^1$.}}
    \label{fig:Nestedness2}
\end{figure}
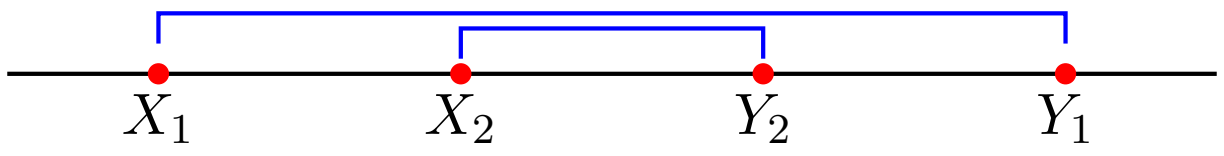

The following remark clarifies how the cross-ratio relates to nestedness.
\begin{remark}
Identify $\RP^1 = \R \cup \{\infty\}$ with $[1:x]$ corresponding to $x$ and $[0:1]$ to $\infty$. Under this convention, we have $[\infty,0,1,x] = x$. For $1<x<\infty$, the cross ratio $x=[\infty,0,1,x]$ measures how nested the segment $[0,1]$ is inside the projective segment $[\infty,x]$. The cross ratio converges to 1 as $x \searrow 1$, as these segments stop being nested, and the cross ratio converges to $+\infty$ as $x \nearrow \infty$ and the segment $[\infty, x]$ converges to $\RP^1$.
\end{remark}

We shall take for granted the following fact regarding 
nestedness of cones and positivity of linear combinations of the associated bilinear forms. 
\begin{remark}\label{Remk:Positivity}
To compare with the writing of \cite{Dav26}, note that if $B_i$ are bilinear forms on $\R^n$ with associated cones $C_i$, the condition $\overline{C_{2}} \subset C_{1}$ is equivalent to the existence of $c_1,c_2>0$ such that $c_1B_1>c_2B_2$  \cite{Fin36}. See also the survey \cite{PT07} for related and more modern results.
\end{remark}
   
The following proposition is (mostly) proven in \cite{Dav26}. 

\begin{proposition}[Nestedness properties]
\label{prop:Nestednessproperties} We have the following nestedness properties:
\begin{enumerate}[label=(\roman*), noitemsep]
    \item Let $C_1,C_2,C_3$ be $d$-cones with $\overline{C_3}\subset C_2$ and $\overline{C_2}\subset C_1$. Then:
    \[ \Nest(C_1,C_2)\Nest(C_2,C_3)\leq \Nest(C_1,C_3). \] 
  \item Suppose $(C_k)_{k\geq 0}$ is an infinite sequence of $d$-cones such that for all $k\geq 0$, $\overline{C_{k+1}}\subset C_k$ and $\lim_{n\rightarrow \infty} \Nest(C_1,C_n)= +\infty$. Then for some $d$-dimensional vector subspace $E\subset \R^n$:
    \[ \bigcap_{k\geq 0} C_k =\mathbb{P}(E). \] 
\end{enumerate}
\end{proposition}

\begin{proof}
    Using \cite[{Proposition 5.2 and 5.3}]{Dav26}, we obtain point (i) and regarding point (ii), we obtain that $\bigcap_{k\geq 0}C_k=\mathbb{P}(E)$ for a vector subspace $E\subset \R^n$. We just need to check that $E$ has dimension $d$. 

    Since $C_1$ is a $d$-cone, its closure is disjoint from $\mathbb{P}(F)$ for some subspace $F\subset \R^n$ of dimension $n-d$. In particular $F\cap E=\lbrace0\rbrace$, so $E$ has dimension at most $d$.

    Conversely, every cone $C_k$ admits a $d$-plane $E_k\subset \R^n$ such that $\mathbb{P}(E_k) \subset C_k$. Up to taking a subsequence, we can assume that the $d$-planes $E_k$ converge to some $d$-plane $E_\infty$. Hence, $\mathbb{P}(E_\infty)\subset \bigcap_{k\geq 0}C_k$, which means $E$ has dimension at least $d$.
\end{proof}

The following result is analogous to \cite[Lemma 3.11]{DB25}. Here, we use the Jordan projection $\boldsymbol{\lambda}$ from \eqref{eq:JordanProj}.

\begin{proposition}[Eigenvalue gaps via nestedness of cones]
\label{prop:Eigenvalue Gap}
    Let $C$ be a $d$-cone and let $g\in \GL(n,\R)$ satisfy $g\cdot \overline{C}\subset C$. Then: 
    \[e^{\lambda_d(g)-\lambda_{d+1}(g)}\geq \Nest(C,g\cdot C). \]
\end{proposition}

\begin{proof}
Let us set $C_0 =C$ and then $C_{i+1}= g \cdot C_i$. 
Using Proposition \ref{prop:Nestednessproperties}(i), we see that $\Nest(C,g^n\cdot C)$ diverges when $n$ goes to $+\infty$. 
Indeed, 
\[ \Nest(C, g^nC) \geq \prod_{i=0}^{n-1} \Nest(g^iC, g^{i+1}C)\geq \Nest(C,g\cdot C)^n,\]
by the invariance of cross-ratio nestedness.
Using Proposition \ref{prop:Nestednessproperties}(ii), there is a $d$-dimensional subspace $E$ of $\R^n$ such that:
$$\bigcap_{k\geq 0} C_k=\mathbb{P}(E).$$

This subspace $E$ is necessarily preserved by $g$. Similarly, considering the $(n-d)$ cones $\mathbb{RP}^{n-1}\setminus \overline{C_1}$ and $\mathbb{RP}^{n-1}\setminus \overline{C_2}$, there is an $(n-d)$-dimensional subspace $F\subset \R^n$ fixed by $g$, and moreover satisfying 
$\mathbb{P}(F) \subset \mathbb{RP}^{n-1}\setminus \overline{C_1}$. Thus, $\R^n = E \oplus F$ and $g$ preserves this splitting. \medskip

Suppose for now that $g$ admits real eigenvectors $e\in E$, $f \in F$ whose corresponding eigenvalues $\zeta_e, \zeta_f$ have minimum modulus $|\zeta_e|=\exp(a_e)$ and maximum modulus $|\zeta_f|=\exp(a_f)$, respectively, across all eigenvectors in $E$, $F$, respectively. Write $\boldsymbol{\lambda}(g) = (\lambda_1,\dots, \lambda_n)$. 
Then $\lambda_d=a_e$ and $\lambda_{d+1}=a_f$. 

In the projective line $L$ spanned by $e$ and $f$, we can find points $(x_1,x_2,e,y_2,y_1,f)$ in cyclic order with $x_2,y_2\in \partial \left( g\cdot C\right)$ and $x_1,y_1\in \partial C$. Note that $y_2 = g\cdot y_1$. 
The linear map $g|_L$ has attracting fixed point $[e]$ and repelling fixed point $[f]$ in $\mathbb{P}(L)$. 
Observe that $[f,e,g(z),z]=[e,f,z,g(z)] =e^{\lambda_d-\lambda_{d+1}}$ is actually independent of $z \in \mathbb{P}(L)$ as long as $z \notin \{\, [f],\,[e] \}$, as noted in \cite{BGLPW24}. Hence:
\[ e^{\lambda_d-\lambda_{d+1}}=[f,e,g\cdot y_1,y_1]=[f,e,y_2,y_1].\]

We now claim that $[f,e,y_2,y_1]\geq [x_1,x_2,y_2,y_1]$. 
To see this, we fix an affine chart of $L$ such that these points correspond to the following elements of $\R\cup \lbrace \infty\rbrace$: $Y_1,F,X_1,X_2,E,Y_2$, with $Y_2=\infty,Y_1=0$ and $F=1$. Then the inequality becomes: 
\[[x_1,x_2,y_2,y_1]=\frac{\infty-X_1}{\infty-X_2}\frac{0-X_2}{0-X_1}=\frac{X_2}{X_1}\leq E=\frac{\infty-1}{\infty-E}\frac{0-E}{0-1}=[f,e,y_2,y_1].\]
This inequality holds as $Y_1=0\leq F=1\leq X_1\leq X_2\leq E\leq \infty=Y_2$.

Hence, we obtain:
\[ e^{\lambda_d-\lambda_{d+1}}\geq [x_1,x_2,y_2,y_1]\geq \Nest(C, g\cdot C).\]
This concludes the proof of the proposition under our assumption about eigenvectors of $g$.
\medskip

If such real eigenvectors of $g$ in $E$ and $F$ do not exist, we can approximate. Take $g_o\in \GL(n,\R)$ arbitrarily close to $g$ such that $g_o$ still preserves $E\oplus F$ and $g_o$ has complex eigenvalues whose arguments are of the form $\frac{2k\pi}{n}$ for some positive integers $n,k$. The proposition applies to $g_o^n$ by the previous argument. Now, by Proposition \ref{prop:Nestednessproperties} we have $\text{Nest}(C, g_o^n\cdot C)\geq \text{Nest}(C, g_o\cdot C)^n$ and hence:
\[e^{\lambda_d(g_o)-\lambda_{d+1}(g_o)}= \left(e^{\lambda_d(g_o^n)-\lambda_{d+1}(g_o^n)}\right)^{\frac{1}{n}}\geq \text{Nest}(C, g_o^n\cdot C)^\frac{1}{n}\geq \text{Nest}(C, g_o\cdot C) . \]

Finally, when $g_o$ approaches $g$ the inequality remains valid, so the proposition applies to $g$.
\end{proof}

Let $\mathscr{B}_{k,n}$ denote the space of symmetric bilinear forms on $\R^n$ of signature $(k,n-k)$.  
We shall use the following sufficient criterion to detect the $P_{k}$-Anosov property of a representation. 
\begin{lemma}[Nested Cones \& Anosov Representations]\label{Lem:ConesAnosov}
Let $\rho: \pi_1S \rightarrow \SL(n,\R)$ be a representation.
Suppose that $B: \T^1\tilde{S} \rightarrow \mathscr{B}_{k,n}$ is a $\rho$-equivariant continuous map. Define $C(v) = \mathbb{P}Q_+(B_v) \subset \mathbb{P}(\R^n)$ as the cone associated to $B_{v}$. 
\begin{enumerate}[label=(\roman*)]
    \item If there is a negatively curved metric $g$ on $S$ such that the induced geodesic flow $\mathcal{G}_t$ on $\tilde{S}$ satisfies 
\[  \overline{C(\mathcal{G}_t(v))} \subset C(v), \; \;  (\forall t > 0) ,\]
then $\rho$ is $P_{k}$-Anosov.  
    \item Moreover, if $v$ points towards $x \in \partial \pi_1S$, then the limit map $\xi^k:\partial\pi_1S\rightarrow \Gr_k(\R^n)$ satisfies 
\[ \xi^k(x) = \bigcap_{t \geq 0} C(\mathcal{G}_t(v)).\]
\end{enumerate} 
\end{lemma}
This Lemma generalizes \cite[{Theorem 5.4}]{Dav26}, and is analogous to \cite[{Theorem 3.6}]{DB25}.

\begin{proof}
(i) First, we observe uniform contraction. 
Now, since $\T^1S$ is compact and $B$ is equivariant and continuous, there exists $\delta>0$ such that $\Nest\big(C(v),C({\mathcal{G}_1(v)})\big)\geq \delta$ for all $v\in \T^1\widetilde{S}$. We imagine the Gromov boundary $\partial \pi_1S$ of $\pi_1S$ to be realized as the visual boundary of $(\widetilde{S},g)$. 

Fix any element $\gamma \in \pi_1S$. Now, let $\gamma_{s} \subset (S,g)$ be the shortest geodesic representative of any element of the free homotopy class of $\gamma$ and $\ell$ be the length of $\gamma_s$. That is, 
$\ell = \ell([\gamma ])$. Define $\tilde{\gamma}_s \subset \tilde{S}$ to be a lift of $\gamma_s$ and set $v$ to be a tangent vector to $\tilde{\gamma}_s$. 
By equivariance of $C$, 
\[ \rho(\gamma)\cdot C(v)=C(\gamma \cdot v)= C(\mathcal{G}_{\ell}(v)). \]
Using Propositions \ref{prop:Nestednessproperties} and Proposition \ref{prop:Eigenvalue Gap}, we obtain
\[  e^{\lambda_k(\rho(\gamma))-\lambda_{k+1}(\rho(\gamma))}\geq \Nest(C(v), \rho(\gamma)\cdot C(v)) = \Nest(C(v),C(\mathcal{G}_\ell(v)))\geq \delta^{\ell-1}= c \cdot \delta^{\ell([\gamma])}, \]
for $c= \frac{1}{\delta}$. 

This proves that $\rho$ is $P_{k}$-Anosov by Definition \ref{defn:Anosov} and Proposition \ref{Prop:Anosov}.\medskip 

(ii) Let us first prove the statement when $x = \gamma^+ \in \partial\pi_1S$ is an attracting fixed point of $\gamma \in \pi_1S$. In this case, a tangent vector pointing towards $\gamma^+$ is of the form $v = \tilde{\gamma}'(t_0)$ for $\tilde{\gamma}$ a lift of $\gamma$ to $\tilde{S}$.  However, $\rho(\gamma)$ sends each cone $C(\mathcal{G}_t(v))$ inside itself. Thus, $\xi^k(x)$, the unique attracting $k$-plane fixed by $\rho(\gamma)$ must lie in $C(\mathcal{G}_t(v))$. The result then follows in this case by Proposition \ref{prop:Nestednessproperties}(ii). 

Now, let $x \in \partial \pi_1S$ be arbitrary and choose $v \in \T^1\tilde{S}$ pointing towards $x$. Since $g$ is negatively curved, we can approach $v$ with vectors $v_n \in \T^1\tilde{S}$ pointing towards $\gamma_n^+$, the attracting fixed point for some $\gamma_n\in \Gamma$. The sequence $\gamma_n^+$ converges to $x$, hence $\xi^k(\gamma_n^+)$ converges to $\xi^k(x)$ by the continuity of $\xi^k$. Thus, $\bigcap_{t \geq 0} C(\mathcal{G}_t(v_n))$ also converges to $\bigcap_{t \geq 0} C(\mathcal{G}_t(v))$, so $\xi^k(x)=\bigcap_{t \geq 0} C(\mathcal{G}_t(v))$.
\end{proof}

\subsection{Nested cones via Higgs bundles}\label{sec:NestedConesHiggs}

Let $(\mathcal{E},\Phi)$ be a cyclic $\sigma$-bundle as in \eqref{eq:betaBundlesAnosovSection},\eqref{eq:betaBundlesAnosovSectionBis} or \eqref{eq:alphaBundlesAnosovSection}. 
Recall the metric $g=g_{\sigma}$ on $S$ from Proposition \ref{Prop:Psi0Parallelism} for $\sigma$-bundles. For $\beta$-bundles, $g=g_{\beta}$ is known to be strictly negatively curved in the Hitchin case by \cite{DaiLi24}. We verify this occurs generally for $\beta$-bundles, and we show this is also the case for $g_{\alpha}$. 
We then construct a $\pi_1S$-equivariant map that associates to $v\in \T^1\widetilde{S}$ a $1$-cone, respectively a $2$-cone, $C_v \subset \mathbb{RP}^6$, for $\beta$-bundles and $\alpha$-bundles. We check that these cones are nested with respect to the $g_{\sigma}$-geodesic flow, and use the results of the previous subsection to verify the corresponding Anosov properties. 

\subsubsection{Nestedness cones for \texorpdfstring{$\beta$}{b}-Hodge bundles}  

Recall the sub-endomorphism $\Phi_0$ of the Higgs field $\Phi$ for the bundle $\mathcal{E}$ given by the sum of the non-vanishing components in \eqref{BetaCyclic}. For every $x\in \Sigma$ and $v\in \T^1_{x}\Sigma$, the endomorphism $\Psi_0(v)=\Phi_0(v)+\Phi_0^{*h}(v)$ is a Hermitian endomorphism of $\mathcal{E}$. Since $\beta (\boldsymbol{\mu}(\Psi_0(v))) \neq 0$, we find that $\Psi_0(v)$ has a unique attracting line $E$ and a unique repelling hyperplane $F$. The line and hyperplane $E$ and $F$ are \emph{transverse}: $E \oplus F = \mathcal{E}$. 
Thus, every vector $w\in \mathcal{E}$ decomposes as $w=w^++w^-$ with $w^+\in E$ and $w^-\in F$. Using the harmonic metric $h$ on $\V$, along with this splitting, we define the following cones in the flat bundle $\V^{\R}$.

\begin{definition}[Cones for $\beta$-bundles]\label{Defn:ConeBeta}
Let $v \in \T^1_xS$. We define $C(v)\subset \mathbb{P}(
\mathcal{E}_{\mid x}^\mathbb{R})$ to be the set of positive lines for the following bilinear form $B_v$ of signature $(1,6)$ on the real locus $\mathcal{E}^\mathbb{R}$:
\[ B_v(w,w) \coloneqq h(w^+,w^+)-h(w^-,w^-) . \] 
In other words, $B_v = h|_{E_{\Psi_0(v)}}-h|_{F_{\Psi_0(v)}} $ and $C(v) = \{ \;[X ]\in \mathbb{P}(
\mathcal{E}_{\mid x}^\mathbb{R}) \mid B_v(X) > 0\}$. 

Similarly, let us introduce the following alternate cones $
{C}'(v)\subset \mathbb{P}(\mathcal{E}^\R_{\mid x})$ of positive lines for the following bilinear form $B'_{v}$ of signature $(1,6)$, where $w_N,w_B$ are the components of $w$ in $N$ and $B$, respectively:
\begin{align}\label{B'Defn}
    B'_{v}(w,w) \coloneqq h(w^+,w^+)-h(w^-,w^-) -h(w_N,w_N)-h(w_B,w_B). 
\end{align}
\end{definition}

We quickly note another description of $B_v'$. As we shall see in the proof of Lemma \ref{Lem:NestedBeta}, $B'_v$ is equivalently described via the harmonic metric and the $\Psi_0(v)$-eigenspaces:
\begin{align}\label{B'DefnAlternate}
     B'_v= h|_{\mathbb{E}_{2}(v)}-h|_{\mathbb{E}_{0}(v)}-h|_{\mathbb{E}_{-2}(v)}-2h|_{\mathbb{E}_{1}(v)}-2h|_{\mathbb{E}_{-1}(v)}.
\end{align}
Here, $\mathbb{E}_k(v)$ denotes the $k$-eigenspace of $\Psi_0(v)$. 

\begin{remark}
    The reason why we introduce two different cones is that the cones $C(v)$ are more natural and yield simpler computations, but the cones $C'(v)$ give a slightly better control on the limit curve that we exploit in the proof of Theorem \ref{thm:BetaAvoidsLimitSet}.
\end{remark}

\begin{remark}
\label{rem:ParallelTransportCones}
    We may lift the cones $C$ in $ \mathbb{P}(\mathcal{E}_{}^\mathbb{R})$ 
    to corresponding cones $\hat{C}$ in $\mathbb{P}( \tilde{\mathcal{E}}^\R) \cong \tilde{S}\times \RP^6$, where the latter identification is naturally given by the flat connection $\nabla$. 
    For $v_0 \in \T^1\tilde{S}$, the resulting map $v_0\mapsto \hat{C}(v_0)\subset \RP^6$ is equivariant with respect to the representation $\rho:\pi_1S\to \Gtwosplit$ associated to the Higgs bundle. Similarly the cones $C'$ in $\mathbb{P}(\mathcal{E}^{\mathbb{R}})$ lift to cones $\hat{C}'$ in $ \mathbb{P}(\tilde{\mathcal{E}}^\R)$. 
\end{remark}

As we shall see, the following lemma verifies the infinitesimal nestedness of the cones $C(v)$ and is the technical backbone of the proof of the Anosov condition.

\begin{lemma}[Nestedness of cones for $\beta$-bundles]
\label{Lem:NestedBeta}
Let $(\V,\Phi)$ be a polystable $\beta$-Hodge bundle on $\Sigma$ such that $\delta$ has at least one zero. 
Fix any point $x \in \Sigma$, any nonzero tangent vector $v \in \T_x\Sigma$, and define $\psi:=\Psi(v) \in \End(\V^{\R}_{\mid x})$. Then the following hold: 
\begin{enumerate}[label=(\roman*)]
    \item If $w\in \V_{\mid x}$ satisfies $B_v(w,w)=0$, then the derivative at $t=0$ of $B_v(e^{t\psi}\cdot w,e^{t\psi}\cdot w)$ is positive. 
    \item If $\delta =0$, and $B_v'(w,w)=0$, then derivative at $t=0$ of $B_v'(e^{t\psi}\cdot w,e^{t\psi}\cdot w)$ is also positive. 
\end{enumerate}
\end{lemma}

\begin{proof}
Take any $x \in \Sigma$ and $v \in \T_x\Sigma$ as in the statement of the lemma. 

Let us fix a local $h$-unitary complex cross-product basis $(x_i)_{i=3}^3$ for the Higgs bundle by Remark \ref{Remk:UnitaryCrossProd}. In such a basis, one can write the endomorphism $\psi:=\frac{\Psi(v)}{||\beta||}$ as follows, for some $\alpha_0 \in \C$ and $\delta_0\in \C$. For a Hodge bundle at least one of $\alpha_0$, $\delta_0$ vanishes, but we shall not use this fact until later.
$$ \psi=\begin{pmatrix}
     0& 1 &0&0&0&\color{blue}{\delta_0}&0\\
    1& 0 &\color{blue}{\overline{\alpha_0}}&0&0&0&\color{blue}{\delta_0}\\
    0& \color{blue}{\alpha_0} &0&i\sqrt{2}&0&0&0\\
    0& 0 &-i\sqrt{2}&0&i\sqrt{2}&0&0\\
    0& 0 &0&-i\sqrt{2}&0&\color{blue}{\overline{\alpha_0}}&0\\
    \color{blue}{\overline{\delta_0}}& 0 &0&0&\color{blue}{\alpha_0} &0&1\\
    0& \color{blue}{\overline{\delta_0}} &0&0&0&1&0
       \end{pmatrix}$$
In this notation, we have $|\alpha_0| = \frac{\lVert \alpha \rVert}{\lVert 1\rVert}$ and $|\delta_0| = \frac{\lVert \delta \rVert}{\lVert 1\rVert}$. 
The matrix $\psi_0$ representing $\frac{\Psi_0(v)}{\lVert \beta\rVert}$ is given by considering all constant coefficient terms, i.e., by removing the blue terms. 

Relative to the fixed basis $(x_i)_{i=3}^{-3}$, we consider the following unitary eigenbasis $U$ of $\psi_0$:
\[ U=\frac{1}{\sqrt{2}}\left(\begin{array}{rrrrrrr}
0 & 0 &  \frac{i}{\sqrt{2}} & - \,  1 & - \frac{i}{\sqrt{2}} & 0 & 0 \\
0 & 0 & \frac{i}{\sqrt{2}} &   \,  1 & -  \frac{i}{\sqrt{2}} & 0 & 0 \\
0 & 0 &   \,  1 & 0 &   \,  1 & 0 & 0 \\
  \,  1 &   \,  1 & 0 & 0 & 0 & 0 & 0 \\
  \,  1 & -  \,  1 & 0 & 0 & 0 & 0 & 0 \\
0 & 0 & 0 & 0 & 0 &   \,  1 &   \,  1 \\
0 & 0 & 0 & 0 & 0 &   \,  1 & -  \,  1
\end{array}\right)\ . \] 
Conjugating $\psi$ by this new basis, we obtain:
\[U\psi U^{-1}=\left(\begin{array}{rrrrrrr}
2 & 0 & 0 & \color{blue}{\frac{1}{4} i \, \sqrt{2} \alpha_0} & \color{blue}{-\frac{1}{4} i \, \sqrt{2} \alpha_0} & \color{blue}{-\frac{1}{4} i \, \sqrt{2} \overline{\alpha_0}} & \color{blue}{-\frac{1}{4} i \, \sqrt{2} \overline{\alpha_0}} \\
0 & -2 & 0 & \color{blue}{\frac{1}{4} i \, \sqrt{2} \alpha_0} & \color{blue}{-\frac{1}{4} i \, \sqrt{2} \alpha_0} & \color{blue}{-\frac{1}{4} i \, \sqrt{2} \overline{\alpha_0}} & \color{blue}{-\frac{1}{4} i \, \sqrt{2} \overline{\alpha_0}} \\
0 & 0 & 0 & \color{blue}{\frac{1}{2} \, \alpha_0} & \color{blue}{-\frac{1}{2} \, \alpha_0} & \color{blue}{\frac{1}{2} \, \overline{\alpha_0}} & \color{blue}{\frac{1}{2} \, \overline{\alpha_0}} \\
\color{blue}{-\frac{1}{4} i \, \sqrt{2} \overline{\alpha_0}} & \color{blue}{-\frac{1}{4} i \, \sqrt{2} \overline{\alpha_0}} & \color{blue}{\frac{1}{2} \, \overline{\alpha_0}} & 1 & 0 & \color{blue}{\delta_0} & 0 \\
\color{blue}{\frac{1}{4} i \, \sqrt{2} \overline{\alpha_0}} & \color{blue}{\frac{1}{4} i \, \sqrt{2} \overline{\alpha_0}} & \color{blue}{-\frac{1}{2} \, \overline{\alpha_0}} & 0 & -1 & 0 & \color{blue}{\delta_0} \\
\color{blue}{\frac{1}{4} i \, \sqrt{2} \alpha_0} & \color{blue}{\frac{1}{4} i \, \sqrt{2} \alpha_0} & \color{blue}{\frac{1}{2} \, \alpha_0} & \color{blue}{\overline{\delta_0}}& 0 & 1 & 0 \\
\color{blue}{\frac{1}{4} i \, \sqrt{2} \alpha_0} & \color{blue}{\frac{1}{4} i \, \sqrt{2} \alpha_0} & \color{blue}{\frac{1}{2} \, \alpha_0} & 0 & \color{blue}{\overline{\delta_0}} & 0 & -1
\end{array}\right).\]
Let us now verify the two claims of the lemma. \medskip  

\textbf{(i)} In particular, note that
$U\psi_0 U^{-1} = \diag(2,-2,0,-1,1,-1,-1)$. 

Thus, in the basis expressed by the column vectors of $U^{-1}$, 
the bilinear form $B_v$ is represented by the following matrix:
\[ Q := \diag(1,-1,-1,-1,-1,-1,-1).\]
That is, $B_v(U^{-1}y, U^{-1}y)=Q(y,y)$. 
By adding a multiple of the identity matrix and a matrix $C$ that belongs to the Lie algebra $\mathfrak{so}(Q) \cong \mathfrak{so}(1,6)$, we obtain a new matrix $A$: 

\[A := U\psi U^{-1} - 2I + C=\left(\begin{array}{rrrrrrr}
0 & 0 & 0 & 0 & 0 & 0 & 0 \\
0 & -4 & 0 & \color{blue}{\frac{1}{4} i \, \sqrt{2} \alpha_0} & \color{blue}{-\frac{1}{4} i \, \sqrt{2} \alpha_0} & \color{blue}{-\frac{1}{4} i \, \sqrt{2} \overline{\alpha_0}} & \color{blue}{-\frac{1}{4} i \, \sqrt{2} \overline{\alpha_0}} \\
0 & 0 & -2 & \color{blue}{\frac{1}{2} \, \alpha_0} & \color{blue}{-\frac{1}{2} \, \alpha_0} & \color{blue}{\frac{1}{2} \, \overline{\alpha_0}} & \color{blue}{\frac{1}{2} \, \overline{\alpha_0}} \\
0 & \color{blue}{-\frac{1}{4} i \, \sqrt{2} \overline{\alpha_0}} & \color{blue}{\frac{1}{2} \, \overline{\alpha_0}} & -1 & 0 & \color{blue}{\delta_0} & 0 \\
0 & \color{blue}{\frac{1}{4} i \, \sqrt{2} \overline{\alpha_0}} & \color{blue}{-\frac{1}{2} \, \overline{\alpha_0}} & 0 & -3 & 0 & \color{blue}{\delta_0} \\
0 & \color{blue}{\frac{1}{4} i \, \sqrt{2} \alpha_0} & \color{blue}{\frac{1}{2} \, \alpha_0} & \color{blue}{\overline{\delta_0}} & 0 & -1 & 0 \\
0 & \color{blue}{\frac{1}{4} i \, \sqrt{2} \alpha_0} & \color{blue}{\frac{1}{2} \, \alpha_0} & 0 & \color{blue}{\overline{\delta_0}} & 0 & -3
\end{array}\right). \]

Let $A_-$ be the $6 \times 6$ bottom block of $A$. We check that $A_-$ is negative definite, which shall yield the lemma. We have two cases, corresponding to the two type of Hodge bundles.  \medskip 

\textbf{Case 1}: $\alpha_0=0$. We have $|\delta_0| <1$ by Lemma \ref{Lem:MaximumPrincipleBeta}. We 
compute $D =\det(A_-) = 72-80|\delta_0|^2+8|\delta_0|^4$. 
For $\abs{\delta_0}<1$, we have $D > 0$. Now, $A_-$ is evidently negatively-definite when $\delta_0=0$. As $A_-$ cannot change signature for $|\delta_0|<1$, this means $A_-$ remains negative-definite for all such $\delta_0$. \medskip 

\textbf{Case 2}: $\delta_0=0$. We reason just as in Case 1. The determinant $D$ of $A_-$ is $D= 2 |\alpha_0|^{4}- 30|\alpha_0|^2 + 72$. By Lemma \ref{Lem:MaximumPrincipleBeta}, we have $|\alpha_0| \leq \sqrt{2}$. 
Note that $D>0$ for all $|\alpha_0|< \sqrt{3}$, so $A_-$ has constant signature. 
Again for $\alpha_0=0$, the matrix $A_-$ is negative-definite, hence it is negative-definite for all $\alpha_0$ such that $|\alpha_0|< \sqrt{3}$. 
\medskip

Thus, $A_-$ is negative definite in all cases. \medskip 

Now, denote $w'=Uw$. By hypothesis, $Q(w',w') = B_v(w,w)=0$. 
Let us write $w'=w'_{+}+w'_{-}$ relative to $\R^{1,6} = \R^{1,0} \oplus \R^{0,6}$, so that $w'_- \in \R^{0,6}$ is nonzero.
Here, $\R^{1,0} = E_{\psi}$ and $\R^{0,6} =F_{\psi}$ are the attracting and repelling line-hyperplane of $\psi$. 
Finally, we conclude: 
\begin{align*}
    \frac{d}{dt} \bigg|_{t=0}B_v(e^{t\psi}w, e^{t\psi }w )&= 2B_v(\psi w, w)=2 Q(U\psi w, Uw)= 2Q(U\psi U^{-1}w',w')
=_{(\star)}2Q(Aw',w')\\
& =2Q(A_-w'_-, w'_-)>0.
\end{align*}
In the penultimate step $(\star)$, we use that $Q\big( (U\psi U^{-1}-A)w',w'\big)=Q(2w',w')+Q(Cw',w')=0$. Thus, the proof of (i) is complete for $B_v$. \medskip

\textbf{(ii)} When $\delta=0$, we consider the bilinear form $B'_v$ represented in the basis $(u_i)_{i=1}^7$ by the matrix:
\[ Q':=\diag(1,-1,-1,-2,-2,-2,-2).\]

As above, we find 
\[ \delta':=\frac{d}{dt}\bigg|_{t=0} B'_v(e^{t\psi}\cdot w,e^{t\psi}\cdot w) = 2Q'(U\psi U^{-1} w',w'). \]  
Thus, to prove $\delta' >0$, we need to prove the associated bilinear form is positive on all vectors $w \in \partial C'(v)$. In symmetrizing, we consider the bilinear form $\mathcal{B}$ given by 
\[ \mathcal{B}(\cdot,\cdot) = Q'\big(U\psi U^{-1}(\cdot), (\cdot )\big) + Q'\big(\cdot , U\psi U^{-1}( \cdot )\big).\]
The bilinear form $\mathcal{B}$ is represented by the 
Hermitian matrix $A':=U\psi U^{-1} Q'+Q'U\psi U^{-1}$. In other words, $\mathcal{B}(\cdot,\cdot ) = h(A'(\cdot), \cdot)$.

We wish to prove that $\mathcal{B}$ is positive on $\partial C'(v)$. 
We use the same trick as above. Note that
$h(Q'(w), w)=0$ by definition for $w \in \partial C'(v)$. Hence, it suffices to prove instead that the bilinear form $h(M(\cdot), (\cdot))$ is positive on $\partial C'(v)$, 
for the following matrix $M$: 
\[ M=U\psi U^{-1} Q'+Q'U\psi U^{-1}-\frac{7}{2}Q'.\]
In coordinates:
\[ M=\left(\begin{array}{rrrrrrr}
\frac{1}{2} & 0 & 0 & \color{blue}{-\frac{1}{4} i \, \sqrt{2} \alpha_0} & \color{blue}{\frac{1}{4} i \, \sqrt{2} \alpha_0} & \color{blue}{\frac{1}{4} i \, \sqrt{2} \overline{\alpha_0}} & \color{blue}{\frac{1}{4} i \, \sqrt{2} \overline{\alpha_0}} \\
0 & \frac{15}{2} & 0 & \color{blue}{-\frac{3}{4} i \, \sqrt{2} \alpha_0} & \color{blue}{\frac{3}{4} i \, \sqrt{2} \alpha_0} & \color{blue}{\frac{3}{4} i \, \sqrt{2} \overline{\alpha_0}} & \color{blue}{\frac{3}{4} i \, \sqrt{2} \overline{\alpha_0}} \\
0 & 0 & \frac{7}{2} & \color{blue}{-\frac{3}{2} \, \alpha_0} & \color{blue}{\frac{3}{2} \, \alpha_0} & -\color{blue}{\frac{3}{2} \, \overline{\alpha_0}} & -\color{blue}{\frac{3}{2} \, \overline{\alpha_0}} \\
\color{blue}{\frac{1}{4} i \, \sqrt{2} \overline{\alpha_0}} & \color{blue}{\frac{3}{4} i \, \sqrt{2} \overline{\alpha_0}} & \color{blue}{-\frac{3}{2} \, \overline{\alpha_0}} & 3 & 0 & 0 & 0 \\
\color{blue}{-\frac{1}{4} i \, \sqrt{2} \overline{\alpha_0}} & \color{blue}{-\frac{3}{4} i \, \sqrt{2} \overline{\alpha_0}} & \color{blue}{\frac{3}{2} \, \overline{\alpha_0}} & 0 & 11 & 0 & 0 \\
\color{blue}{-\frac{1}{4} i \, \sqrt{2} \alpha_0} & \color{blue}{-\frac{3}{4} i \, \sqrt{2} \alpha_0} & \color{blue}{-\frac{3}{2} \, \alpha_0} & 0 & 0 & 3 & 0 \\
\color{blue}{-\frac{1}{4} i \, \sqrt{2} \alpha_0} & -\color{blue}{\frac{3}{4} i \, \sqrt{2} \alpha_0} &\color{blue}{-\frac{3}{2} \, \alpha_0} & 0 & 0 & 0 & 11
\end{array}\right). \]

The determinant of $M$ is :
\[ \det(M) =  2646 \, \abs{\alpha_0}^{4} - \frac{50589}{4} \, \abs{\alpha_0}^{2} + \frac{114345}{8}.\]
The polynomial $\det(M)$ is positive in particular for $\abs{\alpha_0}\leq \sqrt{5/3}$, and hence $M$ is positive-definite for all such $\alpha_0$. By Lemma \ref{Lem:MaximumPrincipleBeta}, we conclude that $\mathcal{B}$ is positive on $C'(v)$, proving (ii).
\end{proof}

Before we verify the Anosov property, we show that $g=g_{\beta}$ is always negatively curved. We show this generally, not just for $\beta$-Hodge bundles. 

\begin{lemma}[Negative curvature]\label{Lem:NegativeCurvatureBeta}
Let $(\V,\Phi)$ be a polystable $\beta$-bundle on a closed Riemann surface $\Sigma=(S,J)$ of genus at least two and $h=\diag(h_i)_{i=3}^{-3}$ be a harmonic metric. Then the metric $g_{\beta}=h_{-1}$ on $S$ is negatively curved. 
\end{lemma}

\begin{proof}
We may write $\beta =1$ for simplicity. Let us apply Lemma \ref{Lem:HitchinsEquationsGeneral} for the conformal metric $g= \lVert 1 \rVert_h$. Under this choice, $\lVert 1\rVert_{h,g}=1$. 
Now, \cite[Lemma 4.11]{DE26b} provides the bound $\lVert \alpha \rVert_{h,g} \leq 2\lVert \beta \rVert_{h,g} = 2$, with equality at one point if and only if global equality holds. The second equation of 
Lemma \ref{Lem:HitchinsEquationsGeneral}, combined with this inequality, shows that $\kappa_g \leq 0$.

We now verify the inequality is strict. If $\kappa_g(p)=0$ at one point $p \in S$, then by the above, $\lVert \alpha\rVert =2\lVert\beta \rVert$ everywhere, which means $\kappa_g \equiv 0$. However, this is impossible as $S$ is closed of genus at least two. Thus, $\kappa_g < 0$.
\end{proof}

Consequently, the associated $J$-holomorphic curves are strictly negatively curved, partially generalizing \cite[Theorem 8.4]{DaiLi24}, which treated the $\Gtwosplit$-Hitchin case. 

\begin{corollary}
Let $\nu:\tilde{S} \rightarrow \mathbb{S}^{2,4}$ be an equivariant alternating $J$-holomorphic curve that is immersed. Then the metric 
$g= -\nu^*g_{\mathbb{S}^{2,4}}$ is strictly negatively curved. 
\end{corollary}

\begin{proof}
This follows immediately from Lemma \ref{Lem:NegativeCurvatureBeta} and \cite{CT24}. Indeed, \cite{CT24} showed that every equivariant (immersed) alternating $J$-holomorphic curve arises from a polystable $\beta$-bundle. 
\end{proof}

Lemma \ref{Lem:NestedBeta} now allows us to verify the $P_{\beta}$-Anosov property, and also get information about the Anosov limit map.

\begin{theorem}[Anosov property for $\beta$-Hodge bundles]
 Let $(\V, \Phi)$ be a polystable $\beta$-Hodge bundle such that $\delta$ has a nonzero divisor. Then 
 \begin{enumerate}[label=(\roman*)]
     \item The cones $\hat{C}(v)$ associated to $B_v$ are nested under the $g=g_{\beta}$-geodesic flow. When $\delta =0$, the same is true for the cones $\hat{C}'(v)$ associated to $B_v'$. 
     \item The associated representation $\rho:\pi_1S\rightarrow \Gtwosplit$ is $P_{\beta}$-Anosov.
 \end{enumerate}
\end{theorem}

\begin{remark}
   The representations corresponding to $\delta= 0$, $\alpha\neq 0$ are irreducible in $\Gtwosplit$, and Theorem \ref{Thm:BetaAnosov} is new in this case. On the other hand, the case $\alpha=0$, $\delta\neq 0$ is known, as the corresponding representations factor through a copy of $\SO(2,2) < \Gtwosplit$. However, in this second case our proof also gives a control on the limit curve, which we use later.
\end{remark}

\begin{proof}
Write $h = (h_i)_{i=3}^{-3}$ for the harmonic metric on $\V$ in the given line bundle splitting $\V = \bigoplus_{i=3}^{-3} \mathcal{L}_i$. 
We then lift the vector bundle $\V$ to $\widetilde{\Sigma}$, along with the corresponding objects $B,\Psi_0,h$. As a slight abuse of notation, in this proof for $x_0\in \widetilde{\Sigma}$ and $v_0\in \widetilde{\Sigma}$ we write $\mathcal{E}_{x_0}, B_{v_0}, \Psi_0(v_0)$ these lifts. Let us also denote by $\mathcal{G}_t$ the geodesic flow on $\T^1\widetilde{\Sigma}$ for the metric $g = h_{-1}$ on $\mathcal{K}^{-1}$. \medskip 

We first prove the nestedness of cones $\hat{C}(v_0)\subset \RP^6$ defined via $B_v$ in \eqref{B'Defn} and Remark \ref{rem:ParallelTransportCones}, along the flow $\mathcal{G}_t$. The same argument will be valid for the cones $\hat{C}'(v)$ when $\delta =0$.

Fix any $v_0 \in \T^1_{x_0} \widetilde{\Sigma}$. Let $\gamma_t$ be the unique $g$-geodesic with $\dot{\gamma}(0)=v$. We write $\mathcal{P}_{s,t}^{\mathcal{C}}: {\V}_{\gamma(s)}\rightarrow {\V}_{\gamma(t)}$ to denote parallel translation along $\gamma$ with respect to a connection $\mathcal{C}$. \medskip 

Fix any point $w \in \partial \hat{C}(v_0)$. Our goal is to compute 
\[ A := \frac{d}{dt}\bigg|_{t=0} B_{\mathcal{G}_t(v_0)}(\mathcal{P}_{0,t}^{\nabla}(w)).\]
That is, we fix $w\in \R^7$ for the identification $\R^7\simeq \mathcal{E}^\R_{x_0}$, which equivalently means that $w$ moves by $\nabla$-parallel transport in $\mathcal{E}^\R_{\gamma(t)}$, and we ask how the norm of $w$ is changing with respect to the quadratic forms $B_{\mathcal{G}_t(v_0)}$. 

Now, let us recall how the quadratic form $B_{v_0}$ is built from the endomorphism $\Psi_0(v_0) \in \End(\V^{\R}_{\mid v})$: 
\begin{align}\label{B'Alternate}
     B_v= h|_{\mathbb{E}_{2}(v_0)}-h|_{\mathbb{E}_{0}(v_0)}-h|_{\mathbb{E}_{-2}(v_0)}-h|_{\mathbb{E}_{1}(v_0)}-h|_{\mathbb{E}_{-1}(v_0)},
\end{align}
where $\mathbb{E}_k(u)$ denotes the $k$-eigenspace of $\Psi_0(u)$ in $\V^{\R}_{\mid u}$. 
By the $\Psi_0$-parallelism in Proposition \ref{Prop:Psi0Parallelism}, 
\begin{align}\label{Psi0GeodesicFlow}
    \mathcal{P}_{0,t}^{\nabla^h} \Psi_0(v_0) = \Psi_0(\mathcal{G}_t(v_0)). 
\end{align}
Due to \eqref{B'Alternate} and \eqref{Psi0GeodesicFlow}, we obtain the key property $B_{\mathcal{G}_t(v_0)}=\mathcal{P}_{0,t}^{\nabla^h}\cdot B_{v_0}$, meaning 
\begin{align}
    B_{\mathcal{G}_t(v_0)}(\cdot ,\cdot)= B_{v_0}(\mathcal{P}_{t,0}^{\nabla^h}(\cdot), \mathcal{P}_{t,0}^{\nabla^h}(\cdot)),
\end{align}
This allows us to treat the moving quadratic forms as living on the common vector space $\V^{\R}_{ \gamma(0)}$. 

Hence, the quantity of interest $A$ is given by 
\[ A = \frac{d}{dt}\bigg|_{t=0} B_{v_0}\big( \mathcal{P}_{t,0}^{\nabla^h}\circ \mathcal{P}_{0,t}^{\nabla}(w)\big) = 2 B_{v_0}\bigg(w, \frac{d}{dt}\bigg|_{t=0}\mathcal{P}_{t,0}^{\nabla^h}\circ \mathcal{P}_{0,t}^{\nabla}(w)\bigg).\]
Recall that parallel transport gives rise to covariant differentiation by $\frac{d}{dt}\big|_{t=0}\mathcal{P}_{t,0}^{\mathcal{C}} =D^{\mathcal{C}}_{t}$, for any connection $\mathcal{C}$.  
Additionally, we have $D^{\nabla}_t= D^{\nabla^h}_t+ \Psi(\mathcal{G}_t(v_0))$. Let us write $w_t = \mathcal{P}_{0,t}^{\nabla}(w)$. By Lemma \ref{Lem:NestedBeta},  $B_{v_0}(w, \Psi(v_0)w)> 0$.  As a consequence:
\begin{align}
    A =B_{v_0}\big(w, D_{t}^{\nabla^h}(w_t)\big) = B_{v_0}\big(w, -\Psi(v_0)(w)\big) < 0.
\end{align}
Hence, for any vector $w \in \partial \hat{C}(v_0)$, we have that $[w] \notin \hat{C}'(\mathcal{G}_t(v))$ for $t$ sufficiently small, meaning $\overline{\hat{C}(\mathcal{G}_t(v_0))} \subset \hat{C}(v_0)$ holds for all $t>0$.

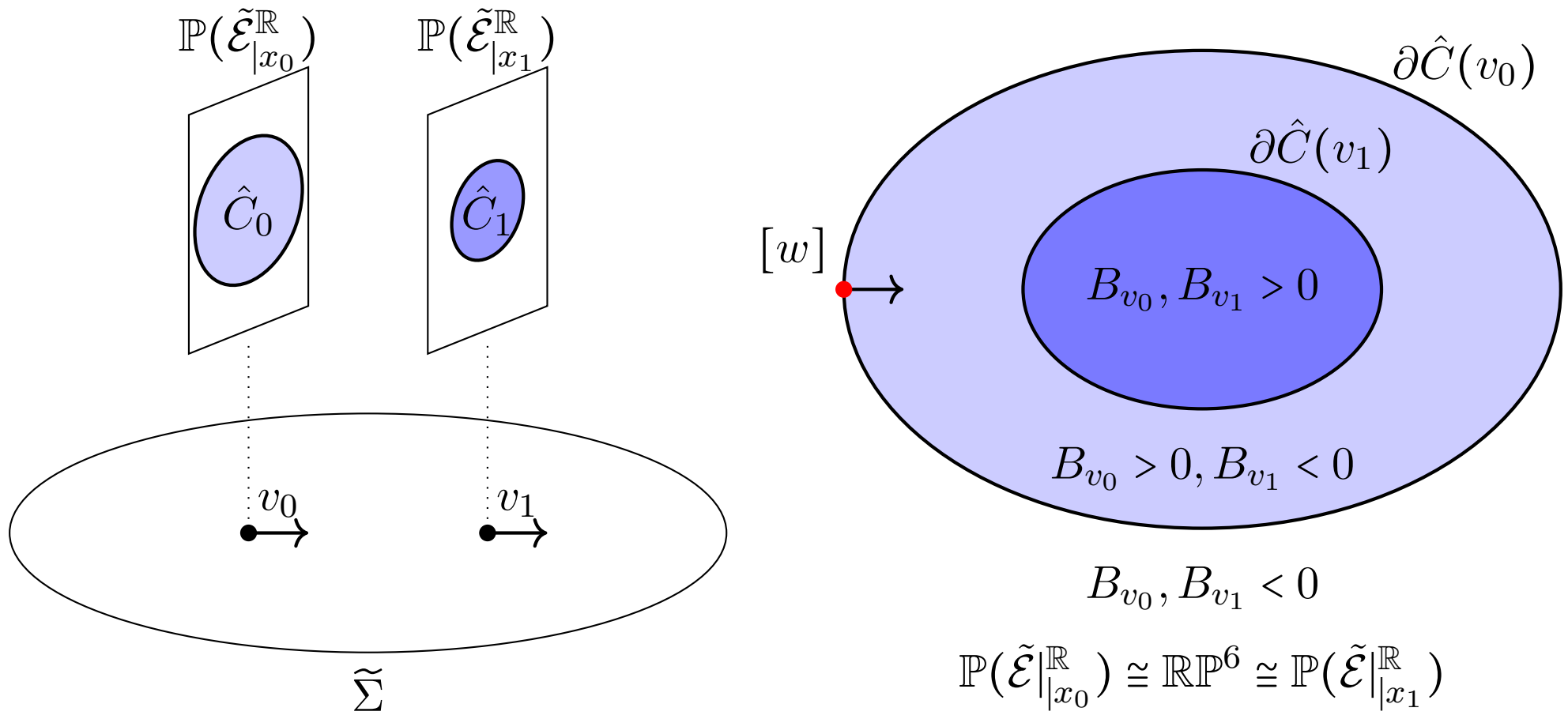
\begin{figure}[h]
    \begin{center}

    \begin{tikzpicture}

    \draw (0,0) ellipse (3cm and 1cm);
    \node[below] at (0,-1) {$\widetilde{\Sigma}$};
    \fill (-1, 0) circle (2pt);
     \node[above] at (-.75,0) {$v_0$};
      \draw[thick,->] (-1,0)  -- (-.5,0);
    \fill (1, 0) circle (2pt);
     \node[above] at (1.25,0) {$v_1$};
     \draw[thick,->] (1,0)  -- (1.5,0);

\begin{scope}[
    shift={(.5,1.5)},
    x={(.5,.2)},
    y={(0,1)}
]
\node[above] at (1,2) {$\mathbb{P}(\tilde{\mathcal{E}}^\R_{|x_1})$};
    \draw
        (0,0) -- (2,0) -- (2,2) -- (0,2) -- cycle;
        \draw[thick,fill=blue, fill opacity=.4] (1,1) ellipse (.6 and .4);
    \node (A) at (1,0) {};
    \node at (1,1) {$\hat{C}_1$};
\end{scope}

\begin{scope}[
    shift={(-1.5,1.5)},
    x={(.5,.2)},
    y={(0,1)}
]
\node[above] at (1,2) {$\mathbb{P}(\tilde{\mathcal{E}}^\R_{|x_0})$};
    \draw
        (0,0) -- (2,0) -- (2,2) -- (0,2) -- cycle;
        \draw[thick,fill=blue, fill opacity=.2] (1,1) ellipse (.9  and .6);

        \node (B) at (1,0) {};
        \node at (1,1) {$\hat{C}_0$};
\end{scope}
    \draw[dotted] (B) -- (-1,0) ;
    \draw[dotted] (A) -- (1,0) ;

\end{tikzpicture}
\begin{tikzpicture}
    \draw[thick, fill=blue, fill opacity=.2] (0,0) ellipse (3cm and 2cm);
    \node at (2.2,1.9) {$\partial \hat{C}(v_0)$};

    \draw[thick,fill=blue, fill opacity=.4] (0,0) ellipse (1.5cm and 1cm);
    \node at (1,1.2) {$\partial \hat{C}(v_1)$};

    \node at (0,0) {$B_{v_0},B_{v_1}>0$};
    \node at (0,-1.5) {$B_{v_0}>0, B_{v_1}<0$};
        \draw[thick,->] (-3,0)  -- (-2.5,0);

            \node at (0,-2.5) {$B_{v_0},B_{v_1}<0$};
            \node at (0,-3.25) {$\mathbb{P}(\tilde{\V}|_{\mid x_0}^{\R})\cong \mathbb{RP}^6\cong \mathbb{P}(\tilde{\V}|_{\mid x_1}^{\R})$};
        \fill[red] (-3, 0) circle (2pt) node[above left, black] {$[w]$};

\end{tikzpicture}
\end{center}
    \caption{\emph{Illustration of the nestedness of the cones, where $v_1=\mathcal{G}_t(v_0)$. For $t >0$, the boundary element $[w] \in \partial \hat{C}({v_0})$ lies on the outside of the smaller cone $\hat{C}({v_1})$. }}
    \label{fig:NestednessParallel}
\end{figure}

The projected metric $g=h_{-1}$ is strictly negatively curved by Lemma \ref{Lem:NegativeCurvatureBeta}. 
Thus, we can now apply Lemma \ref{Lem:ConesAnosov} for $k=1$ to conclude that the holonomy representation $\rho$ of the Higgs bundle is 1-Anosov in $\SL(7,\R)$, and the proof is complete by Definition \ref{Defn:AnosovG2'}. 
\end{proof}

Consequently, we find the limit set is trapped inside of these nested cones as follows.

\begin{corollary}[Trapping the $\beta$-limit set]\label{Cor:BetaLimitSetTrapped}
Let $(\V,\Phi)$ and $\rho:\pi_1S\rightarrow \Gtwosplit$ be as in Theorem \ref{Thm:BetaAnosov}, and $\xi^{\beta}: \partial \pi_1S \rightarrow \Ein^{2,3}$ be the $P_{\beta}$-Anosov limit map of $\rho$. 
Fix $p \in \widetilde{\Sigma}$ and a tangent vector $v_0 \in \T^1_p\widetilde{\Sigma}$ that points towards $x \in \partial \pi_1S$. Then 
\begin{enumerate}[label=(\roman*)]
    \item $\xi^\beta(x) \in \hat{C}(v_0)$. 
    \item If $[u_{\mathscr{L}}+u_N+v_T+v_B]=\xi^\beta(x)$, where $u_{\mathscr{L}}\in \mathscr{L}_p$, $u_N\in N_p$, $v_T\in T_p$, $v_{B}\in B_p$, then:
$$\lVert u_{\mathscr{L}}\rVert^2+\lVert v_T\rVert^2>2\lVert u_N\rVert^2+2\lVert v_B\rVert^2.$$
    \item Moreover, $\xi^\beta(x) \in \hat{C}'(v_0)$.
\end{enumerate}
\end{corollary}

Note that here for a fixed point $p\in \widetilde{\Sigma}$, we identify the fiber $\mathcal{E}^{\R}_{\mid p}=\mathscr{L}_p\oplus T_p\oplus N_p\oplus B_p$ with $\R^{3,4}$.

\begin{proof}
    The fact that the limit curve belongs to $\hat{C}(v_0)$ is an additional consequence of Lemma \ref{Lem:ConesAnosov}, applied to the proof of Lemma \ref{Lem:NestedBeta}. The same argument can be applied to $\hat{C}'(v_0)$ in the case when $\delta=0$. Hence, we have shown (i) in all cases, and (iii) in the case $\delta=0$. 

    (ii) Let $v_0$ and $x$ be given as in the statement. We proceed in two cases.\medskip 
    
    \textbf{Case 1}: $\delta=0$. Let us use the notations from the proof of Lemma \ref{Lem:NestedBeta} we have the general inequality:
    \[B'_{v_0}(u_{\mathscr{L}}+u_N+v_T+v_B)\leq \lVert u_{\mathscr{L}}\rVert^2+\lVert v_T\rVert^2-2\lVert u_N\rVert^2-2\lVert v_B\rVert^2.\]
    This inequality comes from the fact that $B'_{v_0}=-2h$ on $N_p\oplus B_p$, $B'_{v_0}=-h$ on $E_{-1}$, and $B'_v=h$ on $E_1$ for some decomposition $E_1\oplus E_{-1}=\mathscr{L}_p\oplus T_p$. This proves the desired inequality, which then implies $\xi^\beta(x)\in \hat{C}'(v_0)$. Since the inequality is independent of $v_0$, it holds for all $x\in \partial \Gamma$.\medskip 
    
    \textbf{Case 2}: $\alpha =0$. In this case, note that the associated representation $\rho$ preserves the decomposition $\R^{3,4}=\R^{1,2}\oplus \R^{2,2}$ where $\R^{1,2}=\mathscr{L}_p\oplus T_p$ and $\R^{2,2}=N_p\oplus B_p$ holds for every $p \in \tilde{S}$. 
    
    For any element $\gamma\in \pi_1S$, the eigenline $E^1(\rho(\gamma))$ of $\rho(\gamma)$ associated to the highest eigenvalue must belong to $\mathbb{P}(\R^{2,2})$ or $\mathbb{P}(\R^{1,2})$. 
    Hence, for every attracting fixed point $\gamma^+ \in \partial \pi_1S$, we have $\xi^{\beta}(\gamma^+)\in \Ein(\R^{1,2})$ or $\xi^{\beta}(\gamma^+)\in \Ein(\R^{2,2})$. In fact, each of them must make the same choice. That is, $\xi^{\beta}(\gamma^+) \in \Ein(\R^{2,2})$ for every $\gamma$ or $\xi^{\beta}(\gamma^+) \in \Ein(\R^{1,2})$ for every $\gamma$. 
    Hence, by density, $\im(\xi^{\beta})\subset \Ein(\R^{2,2})$ or $\im(\xi^{\beta})\subset \Ein(\R^{1,2})$. 
    The former possibility is ruled out by the fact that $\xi^\beta(x) \in \hat{C}(v_0)$ for some $v_0\in \T_{p}^1\widetilde{\Sigma}$. In the latter case, we clearly have the desired inequality.
    Now, note by \eqref{B'Defn} that $B'_v(X)=B_v(X)-\lVert X_{N}\rVert^2-\lVert X_{B}\rVert^2$. Since $\xi^{\beta}(x) \subset \R^{1,2} = \mathscr{L}_{p}\oplus T_{p}$, we see that 
    $\xi^{\beta} \in \hat{C}({v_0})$ implies $\xi^{\beta}(x) \in \hat{C}'(v_0)$. Thus, (ii) and (iii) hold when $\alpha=0$. 
\end{proof}

We note a further consequence of nestedness that will useful later. Namely, the $J$-holomorphic curve $\nu:\tilde{\Sigma}\rightarrow \Gtwosplit$ associated to a $\beta$-Hodge bundle is transverse to the nested cones $\hat{C}$ and thus transverse to the Anosov limit set.

\begin{corollary}[$\beta$-limit set in an affine chart]\label{Cor:BetaLimitAffineChart}
Let $(\V,\Phi)$ and $\rho$ be as in Theorem \ref{Thm:BetaAnosov} and $\nu:\tilde{\Sigma}\rightarrow \Gtwosplit$ be the associated $J$-holomorphic curve. Then for any point $p \in \tilde{\Sigma}$ and tangent vector $v_0 \in \T^1_p\tilde{\Sigma}$, the cone $\hat{C}(v_0)$ does not intersect $\nu(p)^\bot$. Thus, $\Lambda = \im(\xi^{\beta}) $ is never orthogonal to $\nu$. 
\end{corollary} 

\begin{proof}
Fix any point $p \in \tilde{\Sigma}$. Then for any $x \in \partial \pi_1S$, there is a corresponding tangent vector $v \in \T^1_{p}\widetilde{\Sigma}$ such that $v$ points towards $x$. Hence, $\xi^{\beta}(x) \in \hat{C}(v)$ by Lemma \ref{Lem:ConesAnosov}(ii). Now, similar to the previous proof, we note that  
\[B_{v}(u_{\mathscr{L}}+u_N+v_T+v_B)\leq \lVert u_{\mathscr{L}}\rVert^2+\lVert v_T\rVert^2-\lVert u_N\rVert^2-\lVert v_B\rVert^2.\]
If we then write $\xi^1(x) = [u_{\mathscr{L}} + u_N + v_T +v_B]$, where $\lVert u_{\mathscr{L}}\rVert +\lVert u_{N}\rVert =1=\lVert v_{T}\rVert +\lVert v_{B}\rVert $, we immediately conclude $\lVert u_{\mathscr{L}} \rVert \neq 0$ by Corollary \ref{Cor:BetaLimitSetTrapped}(i). However, $\lVert u_{\mathscr{L}} \rVert = | \nu(p) \cdot (u_{\mathscr{L}} + u_N + v_T +v_B) |$. 
\end{proof}

Finally, we have one further corollary. For $\beta$-Hodge bundles, the associated almost-complex curves $\nu: \tilde{\Sigma} \rightarrow \X$ and their associated minimal surfaces $f: \tilde{\Sigma} \rightarrow \X$ are embedded. 
When $\alpha=0$, these curves are totally geodesic (see \cite{CT24}), so we record only the non-trivial case. 
\begin{corollary}
Let $\nu: \tilde{\Sigma} \rightarrow \quadric$ be an equivariant immersed $J$-holomorphic curve that is superminimal $(\delta =0)$. Then both $\nu$ and the associated minimal immersion $f: \tilde{\Sigma} \rightarrow \X_{\Gtwosplit}$ are embeddings. 
\end{corollary}

\begin{proof}
The key is the relationship between the Frenet frame of $\nu$ and the 1-cones just constructed. If we fix any point $p \in \tilde{\Sigma}$ and any unit tangent vector $v \in \T^1_p\tilde{\Sigma}$, then the proof of Lemma \ref{Lem:NestedBeta} shows that, in fact, each component subspace $\mathscr{L},T,N,B$ of the Frenet frame lies in the null cone of the signature (1,6) bilinear form $B_v$ associated to $v$.  

Now, fix any pair of points $p \neq  p' \in \tilde{\Sigma}$. Then there is a unique $g$-geodesic $\gamma: \R \rightarrow \tilde{\Sigma}$ containing both points. Without loss of generality, write $\gamma(t) = p, \gamma(t') = p'$ for $t < t'$. Then since $\overline{\hat{C}(v')} \subset \hat{C}(v)$, for $v =\gamma'(t), v'=\gamma'(t')$, it follows that $\nu(p) \neq \nu(p')$. Similar reasoning shows that, for $f(p) = \mathscr{L}_p\oplus N_p$, we find that $f(p) \cap f(p') = \{0\}$. 
\end{proof}

\subsubsection{Nestedness cones for \texorpdfstring{$\alpha$}{a}-bundles}\label{Subsec:AlphaNested}
We proceed completely analogously for $\alpha$-bundles. Let $\Phi_0$ be the $\alpha$-component of the Higgs field. For every $v\in \T_x^1S$, then $\Psi_0(v)=\Phi_0(v)+\Phi_0^*(v)$ is a hermitian endomorphism of $\mathcal{E}$. We denote by $E(v)$ its fixed attracting 2-plane and $F(v)$ its fixed repelling codimension $2$-subspace. Since the planes $E(v)$, $F(v)$ are transverse, every vector $w\in \mathcal{E}$ decomposes as $w=w^+ + w^-$, with $w^+\in E(v)$ and $w^-\in F(v)$.
We can then build the following cones. 

\begin{definition}[Cones for $\alpha$-bundles]\label{Defn:ConeAlpha}
Given any unit tangent vector $v \in \T^1_xS$, let $C(v)\subset \mathbb{P}(\V^\R_{\mid x})$ be the cone of positive lines for the following bilinear form of signature $(2,5)$:
\begin{align}\label{Alpha_2cone}
    B_v(w,w)=h(w^+,w^+)-h(w^-,w^-) .
\end{align}
\end{definition}

By Remark \ref{Remk:UnitaryCrossProd}, we can fix a complex cross-product basis for $\V$. In this basis, one can represent $\Psi(v)$ by a matrix $\psi$ as follows: 

\begin{align}\label{Psialpha}
    \psi=\left(\begin{array}{rrrrrrr}
0 & \color{blue}{\overline{\beta_0}} & 0 & 0 & 0 & \color{blue}{\delta_0} & 0 \\
\color{blue}{\beta_0} & 0 & 1 & 0 & 0 & 0 & \color{blue}{\delta_0} \\
0 & 1 & 0 & \color{blue}{i \, \sqrt{2} \,\overline{\beta_0}} & 0 & 0 & 0 \\
0 & 0 & \color{blue}{-i \, \sqrt{2}\beta_0} & 0 & \color{blue}{i \, \sqrt{2} \,\overline{\beta_0}} & 0 & 0 \\
0 & 0 & 0 & \color{blue}{-i \, \sqrt{2} \beta_0} & 0 & 1 & 0 \\
\color{blue}{\overline{\delta_0}} & 0 & 0 & 0 & 1 & 0 & \color{blue}{\overline{\beta_0}} \\
0 & \color{blue}{\overline{\delta_0}} & 0 & 0 & 0 & \color{blue}{\beta_0} & 0
\end{array}\right)\ .
\end{align}
Once again, removing the coefficients different from $1$, i.e., the blue terms, we obtain the matrix associated with $\Psi_0(v)$.

Next, we prove an estimate for the norms $\lVert\delta\rVert$ and $\lVert \beta\rVert$ for $\alpha$-bundles, which will have two key consequences. First, it will imply the projected metric $g=g_{\alpha}$ from Proposition \ref{Prop:Psi0Parallelism} is negatively curved. Second, it will provide the key estimate for verifying the nestedness of cones.

\begin{lemma}[Maximum principle for $\alpha$-bundles]\label{Lem:AlphaBundles_KeyEstimate}
Let $(\V, \Phi)$ be a polystable $\alpha$-cyclic Higgs bundle on a closed Riemann surface $\Sigma$ of genus $g\geq 2$ and $h$ be a harmonic metric. Then we have the following global inequality:
\[ 6\lVert\beta\rVert^2+2\lVert\delta\rVert^2<4\lVert 1\rVert^2. \]
\end{lemma}

\begin{proof}
Let $g$ be a conformal metric on $\Sigma$. As before, we write $\Delta_g$ for the global Laplacian, locally written $\Delta_g= \frac{1}{g_0} \partial_{z}\partial_{\zbar}$ for $g= g_0|dz|^2$. Let us write $P=6\lVert\beta\rVert^2+2\lVert\delta\rVert^2$. To prove the desired inequality, we will compute the Laplacian of $\log(P)$. To this end, suppose $f_1,f_2\geq 0$ are non-negative functions on $S$, each with isolated zeros. To compute $\Delta_g \log P$, we will use the following inequality, which holds away from the union of the finite zero sets of $f_1, f_2$:
\begin{align}\label{Elliptic_Inequality}
    \Delta_g\log(f_1+f_2)\geq \frac{f_1}{f_1+f_2}\Delta_g\log(f_1)+\frac{f_2}{f_1+f_2}\Delta_g\log(f_2).
\end{align}
For the proof of \eqref{Elliptic_Inequality}, we write simply $\Delta$ rather than $\Delta_g$, as the equations are conformally invariant. 

First, a simple calculation finds
\[ \Delta\log(f)=\frac{\Delta f}{f}-\frac{1}{2}\frac{\lVert\nabla f \rVert^2}{f^2}. \] 
Hence, we compute that 
\[ (f_1+f_2)\Delta \log(f_1+f_2)-f_1\Delta\log(f_1)-f_2\Delta\log(f_2)= \frac{1}{2} \left ( \frac{\lVert\nabla f_1 \rVert^2}{f_1}+\frac{\lVert\nabla f_2 \rVert^2}{f_2}-\frac{\lVert\nabla f_1 +\nabla f_2\rVert^2}{f_1+f_2} \right).\]
Multiplying the RHS by $f_1+f_2$ and using the Cauchy-Schwarz inequality, we obtain: 
\begin{align*}
    \left(1+\frac{f_2}{f_1}\right)\lVert\nabla f_1 \rVert^2+\left(1+\frac{f_1}{f_2}\right)\lVert\nabla f_2 \rVert^2-\lVert\nabla f_1 +\nabla f_2\rVert^2&\geq \frac{f_2}{f_1}\lVert\nabla f_1 \rVert^2+\frac{f_1}{f_2}\lVert\nabla f_2 \rVert^2-2\lVert\nabla f_1 \rVert\lVert\nabla f_2 \rVert \\
    &= \frac{(f_1||\nabla f_2|| - f_2||\nabla f_1||)^2}{f_1f_2} \geq 0 .
\end{align*} 
We conclude that \eqref{Elliptic_Inequality} holds on $S \setminus (\mathcal{Z}(f_1) \cup \mathcal{Z}(f_2))$, where $\mathcal{Z}(f_i)$ is the zero set of $f_i$. 

Hitchin's equation from Lemma \ref{Lem:HitchinsEquationsGeneral} applied to $\alpha$-bundles yields the following:
\begin{align}
\label{eq:alphabundle1}
    \Delta_g\log\left(\frac{\lVert\beta\rVert ^2}{\lVert1\rVert^2}\right)&=5\lVert\beta\rVert^2+\lVert\delta\rVert^2-3\lVert1\rVert^2, \\
    \label{eq:alphabundle2}
    \Delta_g\log\left(\frac{\lVert\delta\rVert ^2}{\lVert1\rVert^2}\right)&=3\lVert\beta\rVert^2+3\lVert\delta\rVert^2-3\lVert1\rVert^2.
\end{align}
Now, to prove the desired result, we compute the Laplacian of $\log(\frac{P}{\lVert 1 \rVert})$. 
Using \eqref{Elliptic_Inequality}, we obtain:
\[ \Delta_g\log\left(\frac{P}{\lVert 1\rVert^2}\right)\geq \frac{1}{P}\left(6\lVert \beta\rVert^2\Delta_g\log\left(\frac{\lVert\beta\rVert ^2}{\lVert1\rVert^2}\right)+2\lVert \delta\rVert^2\Delta_g\log\left(\frac{\lVert\delta\rVert ^2}{\lVert1\rVert^2}\right)\right)=:\frac{X}{P}. \]

As already emphasized, this inequality only makes sense away from the isolated zeros of $\beta$ and $\delta$. Now we can compute the term $X$:
\begin{align*}
    X&=6\lVert \beta\rVert^2\left(5\lVert\beta\rVert^2+\lVert\delta\rVert^2-3\lVert 1\rVert^2\right)+2\lVert \delta\rVert^2\left(3\lVert\beta\rVert^2+3\lVert\delta\rVert^2-3\lVert1\rVert^2\right)\\
    &= Q-3 P\lVert1\rVert^2,
\end{align*}
where $Q$ is given by 
\[Q=30\lVert\beta\rVert^4+12\lVert\beta\rVert^2\lVert\delta\rVert^2+6\lVert\delta\rVert^4.\] Next, we note the following remarkable factorization: 
\begin{equation}
\label{eq:Factorization}Q-\frac{3}{4}P^2=3(\lVert \beta\rVert^2-\lVert\delta\rVert^2)^2\geq 0.\end{equation}
Combining \eqref{eq:Factorization} with the inequality for $\Delta_g \left(\frac{P}{\lVert 1\rVert^2} \right)$, we obtain following inequality, that is valid 
away from the common zeros of $\delta$ and $\beta$:
\[ \Delta_g\log\left(\frac{P}{\lVert 1\rVert^2}\right)\geq \frac{3}{4}\left(P-4\lVert 1\rVert^2\right). \]

In particular, this lower bound is non-positive at a global maximum of $\frac{P}{\lVert 1 \rVert^2}$. Thus, at a maximum point $p_0$ on $S \setminus (\mathcal{Z}(\beta) \cup \mathcal{Z}(\delta))$, we have $P(p_0)\leq 4\lVert 1\rVert^2(p_0)$. Since $\log \left( \frac{P}{\lVert 1\rVert ^2} \right) \rightarrow -\infty$ as $P$ approaches a zero of $\beta$ or $\delta$, we conclude $P \leq 4\lVert 1 \rVert^2$ holds globally on $S$. \medskip

Now, let us consider the elliptic PDE on $S \setminus (\mathcal{Z}(\beta) \cup \mathcal{Z}(\delta))$ given by 
\[ F[u]:=\Delta_g u - \Delta_g \log(\lVert 1\rVert^2)-\frac{3}{4}(e^u-4\lVert 1 \rVert^2) = 0. \]
The operator $F$ is elliptic in the sense of \cite[{Definition 3.3.1}]{jost2007partial}. Define $\eta = \log P$ and $\eta_0 = \log(4\lVert 1 \rVert^2)$. By the previous step, we have $F[\eta] \geq 0$. Also, one sees $F[\eta_0 ] =0$ directly. 
Now, the strong maximum principle \cite[Theorem 3.3.1]{jost2007partial} says we either have $P< 4\lVert 1\rVert^2$ globally or $P = 4\lVert 1 \rVert^2$ globally on $S \setminus (\mathcal{Z}(\beta) \cup \mathcal{Z}(\delta))$. 
We show the latter is impossible. 

Now, suppose $P = 4\lVert 1 \rVert^2$ holds on $S \setminus (\mathcal{Z}(\beta) \cup \mathcal{Z}(\delta))$. Then $P = 4\lVert 1 \rVert^2$ holds on all of $S$ by continuity. However, retracing the inequalities for $\Delta_g \log\frac{P}{\lVert 1\rVert ^2}$, one finds we must have equality in \eqref{eq:Factorization}. Hence, $\lVert \delta\rVert=\lVert \beta\rVert$. As a consequence, $\lVert \delta\rVert^2=\lVert \beta\rVert^2= \frac{1}{2}\lVert 1\rVert^2$, so both of $\beta$ and $\delta$ must be non-vanishing. However, for an $\alpha$-bundle $(\V,\Phi)$ on a closed Riemann surface $\Sigma$ of genus $g\geq2$, it is impossible that both $\beta \in H^0(\mathcal{T}^{-1}\K)$ and $\delta \in H^0(\mathcal{T}^3\K^3)$ are non-vanishing. 

Therefore, we must have $P< 4\lVert 1\rVert^2$ globally on $S\setminus (\mathcal{Z}(\beta) \cup \mathcal{Z}(\delta))$. One can easily check $P< 4\lVert 1\rVert^2$ holds also on  $\mathcal{Z}(\beta) \cup \mathcal{Z}(\delta)$, so the desired inequality holds globally on $S$.
\end{proof}

\begin{corollary}[Curvature of $\alpha$-projected metric]\label{Cor:AlphaMetricNegativelyCurved}
Let $(\V,\Phi)$ and $h =\diag(h_i)_{i=3}^{-3}$ be as in Lemma \ref{Lem:AlphaBundles_KeyEstimate}. Then the Riemannian metric $g_{\alpha}= h_{2}^{-1}h_1$ is negatively curved. 
\end{corollary}

\begin{proof}
We consider Hitchin's equations \eqref{HitchinEquationsDifferentials} for $g=g_{\alpha}$. 
We may write $\alpha=1$. 
Hence, $\lVert 1 \rVert=\lVert 1 \rVert_{h,g}= 1$ holds by definition. Thus, 
\[ \kappa_g = 3\lVert \beta \rVert^2+\lVert\delta \rVert^2-2\lVert 1 \rVert^2 < 0. \qedhere \]
\end{proof}

Next, we prove the analogue of Lemma \ref{Lem:NestedBeta}, which shall yield the $P_{\alpha}$-Anosov property. 

\begin{lemma}[Nestedness of cones for $\alpha$-bundles]
\label{Lem:NestedAlpha}
Let $(\V,\Phi)$ be a polystable $\alpha$-bundle on $\Sigma$ and $h$ be a harmonic metric.
Choose any point $p \in \Sigma$ and tangent vector $v\in \T^1\Sigma$ and set $\psi = \Psi(v)$. If $w\in \V_{\mid p}$ satisfies $B_v(w,w)=0$, then the derivative at $t=0$ of $B_v(e^{t\psi}\cdot w,e^{t\psi}\cdot w)$ is positive. 
\end{lemma}

\begin{proof}
We once again work in a local $h$-unitary complex cross product basis $(x_i)_{i=3}^{-3}$ for $\V$ such that $\psi=\Psi(v)$ obtains the form \eqref{Psialpha}. Observe that $|\delta_0|=\frac{\lVert \delta \rVert}{\lVert1\rVert}$ and $|\beta_0|= \frac{\lVert \beta \rVert}{\lVert 1 \rVert }$ in this notation. 

Consider the following unitary change of base:
\begin{align}\label{Psi0AlphaEigenvectors}
    U=\frac{1}{\sqrt{2}}\left(\begin{array}{rrrrrrr}
0 & 1 & 1 & 0 & 0 & 0 & 0 \\
0 & 0 & 0 & 0 & 1 & 1 & 0 \\
0 & 1 & -1 & 0 & 0 & 0 & 0 \\
0 & 0 & 0 & 0 & 1 & -1 & 0 \\
\sqrt{2} & 0 & 0 & 0 & 0 & 0 & 0 \\
0 & 0 & 0 & \sqrt{2} & 0 & 0 & 0 \\
0 & 0 & 0 & 0 & 0 & 0 & \sqrt{2}
\end{array}\right)
\end{align}
In the new basis $\psi$ obtains the form:
\[U\psi U^{-1}=\left(\begin{array}{rrrrrrr}
1 & 0 & 0 & 0 & \color{blue}{\frac{1}{\sqrt{2}}  \beta_0} & \color{blue}{i \, \overline{\beta_0}} & \color{blue}{\frac{1}{\sqrt{2}}  \delta_0} \\
0 & 1 & 0 & 0 & \color{blue}{\frac{1}{\sqrt{2}} \overline{\delta_0}} & -\color{blue}{i \, \beta_0} & \color{blue}{\frac{1}{\sqrt{2}}  \overline{\beta_0}} \\
0 & 0 & -1 & 0 & \color{blue}{\frac{1}{\sqrt{2}} \beta_0} & \color{blue}{-i \, \overline{\beta_0}} & \color{blue}{\frac{1}{\sqrt{2}}  \delta_0} \\
0 & 0 & 0 & -1 & \color{blue}{-\frac{1}{\sqrt{2}} \overline{\delta_0}} & \color{blue}{-i \, \beta_0} & \color{blue}{-\frac{1}{\sqrt{2}}  \overline{\beta_0}} \\
\color{blue}{\frac{1}{\sqrt{2}}  \overline{\beta_0}} & \color{blue}{\frac{1}{\sqrt{2}}  \delta_0} & \color{blue}{\frac{1}{\sqrt{2}}  \overline{\beta_0}} & \color{blue}{-\frac{1}{\sqrt{2}}  \delta_0} & 0 & 0 & 0 \\
\color{blue}{-i \, \beta_0} & \color{blue}{i \, \overline{\beta_0}} & \color{blue}{i \, \beta_0} & \color{blue}{i \, \overline{\beta_0}} & 0 & 0 & 0 \\
\color{blue}{\frac{1}{\sqrt{2}}  \overline{\delta_0}} & \color{blue}{\frac{1}{\sqrt{2}}  \beta_0} & \color{blue}{\frac{1}{\sqrt{2}}  \overline{\delta_0}} & \color{blue}{-\frac{1}{\sqrt{2}}  \beta_0} & 0 & 0 & 0
\end{array}\right).\]
The bilinear form $B_v$ associated to $v$ is represented by the matrix $Q:=\diag(1,1,-1,-1,-1,-1,-1)$ in the basis given by the rows of $U$, or the columns of $U^{-1}$.  

As in the proof of Lemma \ref{Lem:NestedBeta}, we consider a convenient alteration to $U\psi U^{-1}$ as follows, where $C$ is some element of $\mathfrak{so}(Q)\simeq\mathfrak{so}(2,5)$: 
\[A'=U\psi U^{-1}-I+C=\left(\begin{array}{rrrrrrr}
0 & 0 & 0 & 0 & 0 & 0 & 0 \\
0 & 0 & 0 & 0 & 0 & 0 & 0 \\
0 & 0 & -2 & 0 & \color{blue}{\frac{1}{\sqrt{2}}  \beta_0} & \color{blue}{-i \, \overline{\beta_0}} & \color{blue}{\frac{1}{\sqrt{2}} \delta_0} \\
0 & 0 & 0 & -2 & \color{blue}{-\frac{1}{\sqrt{2}} \,\overline{\delta_0}} & \color{blue}{-i \, \beta_0} & \color{blue}{-\frac{1}{\sqrt{2}} \;\overline{\beta_0}} \\
0 & 0 & \color{blue}{\frac{1}{\sqrt{2}}\, \overline{\beta_0}} & \color{blue}{-\frac{1}{\sqrt{2}} \delta_0} & -1 & 0 & 0 \\
0 & 0 & \color{blue}{i \, \beta_0} & \color{blue}{i \; \overline{\beta_0}} & 0 & -1 & 0 \\
0 & 0 & \color{blue}{\frac{1}{\sqrt{2}}  \;\overline{\delta_0}} & \color{blue}{-\frac{1}{\sqrt{2}} \beta_0} & 0 & 0 & -1
\end{array}\right ) \]
The determinant $D$ of the $5$ by $5$ block $A''$ remaining satisfies: 
\begin{align}
    -D &=4+\frac{1}{4} |\delta_0|^4+\frac{1}{2}|\delta_0|^2|\beta_0|^2+\frac{5}{4} \, |\beta_0|^4 - 6 \, |\beta_0|^2 -2|\delta_0|^2 +\beta_0^3\overline{\delta_0}+\overline{\beta_0^3}\delta_0 \\
    -D &\geq 4+\frac{1}{4} |\delta_0|^4+\frac{1}{2}|\delta_0|^2|\beta_0|^2+\frac{5}{4} \, |\beta_0|^4 - 6 \, |\beta_0|^2 -2|\delta_0|^2-2\abs{\beta_0^3\delta_0}  .
\end{align}
We note here that the lower bound can be factorized. 
In particular, we obtain 
\[-D\geq (5\abs{\beta_0}^2+\abs{\delta_0}^2+2\abs{\beta_0\delta_0}-4)\bigg(\frac{1}{4}(\abs{\beta_0}-\abs{\delta_0})^2-1\bigg). \]

The condition $6|\beta_0|^2+2|\delta_0|^2 < 4$ from Lemma \ref{Lem:AlphaBundles_KeyEstimate} implies both 
$5\abs{\beta_0}^2+\abs{\delta_0}^2+2\abs{\beta_0\delta_0}< 4$ and $(\abs{\beta_0}-\abs{\delta_0})^2< 4$. Hence, $-D>0$. We conclude that for all $|\alpha_0|,|\delta_0|$ such that $5\abs{\beta_0}^2+\abs{\delta_0}^2+2\abs{\beta_0\delta_0}< 4$, the signature of $A''$ is constant. As $A''$ is negative-definite for $\delta_0=\beta_0=0$, we conclude it is negative-definite for all such $\delta_0,\beta_0$. Hence, $A'$ is negative-definite, which implies the desired result by the same reasoning as in Lemma \ref{Lem:NestedBeta}.
\end{proof}

We can now prove the Anosov property for $\alpha$-bundles.

\begin{proof}[{Proof of Theorem \ref{thm:AlphaAnosov}}]
The proof is completely analogous to the proof of Theorem \ref{Thm:BetaAnosov}. Here, we use the same reasoning as in Theorem \ref{Thm:BetaAnosov} to conclude the nestedness of cones from Lemma \ref{Lem:NestedAlpha}. Due to Corollary \ref{Cor:AlphaMetricNegativelyCurved}, we can apply Lemma \ref{Lem:ConesAnosov}(i) to conclude that $\rho$ is $P_2$-Anosov in $\SL(7,\R)$, or equivalently, $P_{\alpha}$-Anosov in $\Gtwosplit$. 
\end{proof}

Via the nestedness of cones, the limit set is again trapped. 

\begin{corollary}[Trapping the $\alpha$-limit set]\label{Cor:AlphaLimitSetTrapped}
Let $(\V, \Phi)$ and $\rho:\pi_1S\rightarrow \Gtwosplit$ be as in Theorem \ref{thm:AlphaAnosov} and $\xi^{\alpha}: \partial \pi_1S \rightarrow \Pho^\times$ be the $P_{\alpha}$-Anosov limit map. 
Fix a point $p \in \widetilde{\Sigma}$ and a tangent vector $v \in \T^1\widetilde{\Sigma}$ such that $v$ points towards $x \in \partial \pi_1S$. Then $\mathbb{P}(\xi^\alpha(x)) \subset C(v)$.

Moreover, if $(u_{\mathscr{L}}+u_N+v_T+v_B)\in \xi^{\alpha}(x)$, with $u_{\mathscr{L}}\in \mathscr{L}_p$, $u_N\in N_p$, $v_T\in T_p$, $v_B\in B_p$, then:
\[ \lVert u_N\rVert^2+\lVert v_T\rVert^2>\lVert u_{\mathscr{L}}\rVert^2+\lVert v_B\rVert^2. \] 
\end{corollary}

\begin{proof}
Let $p,v,x$ be given as in the statement and suppose $w:=(u_{\mathscr{L}}+u_N+v_{T}+v_B) \in \xi^{\alpha}(x)$.

By Lemma \ref{Lem:ConesAnosov}(ii), we have $w \in C(v)$, meaning $B_v(w,w) >0$. 
Denote $E(v)$ and $F(v)$ as the attracting and repelling $2$ and $5$-planes of $\Psi_0(v)$, respectively. 
Let $w_+$ be the orthogonal projection of $w$ to $E(v)$ and $w_-$ the projection to $F(v)$. The positivity of $B_v(w,w)$ means that $\lVert w_+\rVert^2 >\lVert w_-\rVert^2$. Note that $E(v)\subset N_p\oplus T_p$ and $\mathscr{L}_p\oplus B_p \subset F(v)$, by considering
the $\Psi_0(v)$-eigenvector matrix \eqref{Psi0AlphaEigenvectors}. Hence:
\[\lVert u_N \rVert^2+\lVert v_T \rVert^2\geq \lVert w_+ \rVert^2> \lVert w_- \rVert^2 \geq \lVert u_\mathscr{L} \rVert^2+\lVert v_B \rVert^2. \qedhere\]
\end{proof}

\section{Limits of Frenet frames}

In order to understand limiting properties of Frenet frames, we will use the following result.

\begin{proposition}[Extending limit maps]
\label{prop:Limit maps}
    Let $\rho:\pi_1S\to \SL(n,\R)$ be a $d$-Anosov representation, for $1\leq d\leq n-1$, and let $\xi^{i}:\partial \pi_1S\to \Gr_{i}(\R^n)$ be the corresponding limit maps for $i \in \{d, n-d\}$. 
    
    Suppose $f^d:\widetilde{S}\to \Gr_{d}(\R^n)$ is a $\rho$-equivariant continuous map such that for all $p\in \widetilde{S}$ and $x\in \partial \pi_1S$, the $d$-plane $f^d(p)$ is transverse to $\xi^{n-d}(x)$. Then $f^d$ extends continuously to 
    \[ f^d \sqcup \xi^d: \widetilde{S} \sqcup \partial\pi_1S \rightarrow \Gr_d(\R^n) .\] 
\end{proposition}

This proposition is an application of the North-South dynamics of Anosov representations from \cite{KLP18}. We include an elementary proof for completeness.
\begin{proof}
Let $(p_i)$ be any sequence of points on $\widetilde{S}$ converging to $x\in \partial\pi_1S$. We show that for any such sequence, up to taking a further subsequence, $f^d(p_i)$ converges to $\xi^d(x)$. This will conclude the proof by the sub-subsequence lemma. 

Let $K$ be a  compact fundamental domain for the action of $\pi_1S$ on $\widetilde{S}$ and fix for all $i\geq 0$ an element $\gamma_i\in \pi_1S$ such that $p_i\in \gamma_i\cdot K$. We consider the Cartan decomposition \[\rho(\gamma_i) =k_i\exp(a_i)k'_i \]
with $k_i, k'_i\in \SO(n,\R)$ and $a_i=\diag(\sigma_1^{(i)},\sigma_2^{(i)},\cdots, \sigma_n^{(i)})$ with $\sigma_1^{(i)}\geq\sigma_2^{(i)}\geq \cdots\geq \sigma_n^{(i)}$. In other words, $\boldsymbol{\mu}(\rho(\gamma_i))= a_i$ is the Cartan decomposition of $\rho(\gamma_i)$ in the group $\SL(n,\R)$. 

The sequence $(\gamma_i)$ converges to $x$ and $(\gamma_i^{-1})$ converges up to taking a subsequence to some $y\in \partial \pi_1S$. Let $E$ be the subspace generated by the first $d$ vectors of the standard basis of $\R^n$ and $F$ the subspace generated by the complementary $n-d$ standard basis vectors. The $d$-Anosov condition of $\rho$ implies that $\sigma_d^{(i)}-\sigma_{d+1}^{(i)}\to +\infty$, along with 
$k_i\cdot E\to \xi^d(x)$ and $(k_i')^{-1}\cdot F\to \xi^{n-d}(y)$ (see \cite{Tso26} and the references therein). Up to taking a subsequence we assume that $k_i,k'_i$ converge to $k,k'\in \SO(n,\R)$. Let $K_0\subset \Gr_d(\R^n)$ be a compact neighborhood of $k'\cdot f^d(K)$ containing only subspaces transverse to $F$, which one can find since $F=k'\cdot \xi^{n-d}(y)$. 

The gap in eigenvalues of $a_i$ implies that since $K_0$ is a compact set of subspaces transverse to $F$, then the sequence of compact subsets $\exp(a_i)\cdot K_0$ converge uniformly to the single point $E\in\Gr_d(\R^n)$. Since $k_i\cdot E$ converges to $\xi^d(x)$ then the compact sets $\rho(\gamma_i)\cdot f^d(K)=f^d(\gamma_i\cdot K)$ converge uniformly to the point $\xi^d(x)$. In particular $f^d(p_i)$ converges to $\xi^d(x)$. 
\end{proof}

As an application, we can use the nested cones to verify transversality between certain subspaces of the ambient vector space with the Anosov limit sets. For the proof, recall that $(\xi^i)^\perp=\xi^{7-i}$ for $1\leq i\leq 6$, by Proposition \ref{prop:LimitMapPropertiesG2}. 

Now, writing $\Gamma =\pi_1S$, we show how Frenet frames degenerate to the limit sets of interest. 

\begin{theorem}[Frenet frames to $\beta$-limit set]\label{thm:BetaLimitExplicitly}
Let $(\mathscr{L}_p, T_p,N_p,B_p)_{p\in \widetilde{S}}$ be the Frenet frame associated to a polystable $\beta$-Hodge bundle such that $\delta$ has at least one zero. Let $x\in \partial \Gamma$ and $\xi^d:\partial \Gamma\to \Gr_d(\R^{3,4})$, for $d\in \{1,3,4,6\}$, be the Anosov limit maps of the associated representation $\rho:\pi_1S \rightarrow \Gtwosplit$. Then when $p\in \widetilde{S}$ converges to $x$, we have:
    \begin{enumerate}[label=(\roman*)]
\item $\mathscr{L}_p\to \xi^1(x)$
\item $\mathscr{L}_p\oplus B_p\to \xi^3(x)$
\item $T_p\oplus N_p\to \xi^4(x)$
\item $T_p\oplus N_p \oplus B_p\to \xi^6(x)$.
    \end{enumerate}
\end{theorem}
 
\begin{proof}
Using Proposition \ref{prop:Limit maps}, it suffices to check certain transversalities for all $p\in \widetilde{S}$ and $x\in \partial\pi_1S$. 

(iv) The relevant transversality property holds by Corollary \ref{Cor:BetaLimitAffineChart}. \medskip 

(iii) By Remark \ref{rem:CircleEinstein}, $\Ein(T_p\oplus N_p) \subset \dev(\overline{{\Fiber}_\sigma}_{|p})$. Hence, Theorem \ref{thm:BetaAvoidsLimitSet} gives that for any $p \in \tilde{S}$, the torus $\Ein(N_p\oplus T_p)$ avoids the $\beta$-thickening of $\Lambda = \im(\xi^1)$. However, $K_{\xi^1(x)} =\mathbb{P}(\xi^3(x))$. Thus, we have 
$\Ein(N_p\oplus T_p) \cap \xi^3(x)=\{0\}$
for any $x \in \partial \pi_1S$. Since the 3-plane $\xi^3(x)$ is isotropic, this means $\xi^3(x) \cap (N_p\oplus T_p) =\{0\}$.  

Now, the relevant transversality needed to verify (i) follows from that of (iv) and the transversality needed for (ii) follows from that of (iii).
\end{proof}

We have the following parallel result for $\sigma =\alpha$. 
When $\beta \neq 0$ for an $\alpha$-bundle, the following result
says that the tangent spaces $(T_p)_{p \in \tilde{S}}$ of the (branched) $J$-holomorphic curve $\nu: \tilde{S} \rightarrow \mathbb{S}^{2,4}$ (recall \S\ref{Subsec:SigmaBundles}) associated to an $\alpha$-bundle converge to the limit set $\Lambda =\im(\xi^{\alpha})\subset \Pho^\times \cong\Gtwosplit/P_{\alpha}$. 

\begin{theorem}[Frenet frames to $\alpha$-limit set]\label{thm:AlphaLimitExplicit}
Let $(\mathscr{L}_p, T_p,N_p,B_p)_{p\in \widetilde{S}}$ be the Frenet frame associated to a polystable $\alpha$-bundle. Let $x \in \partial \Gamma$ and $\xi^{d}:\partial \Gamma \rightarrow \Gr_i(\R^7)$ be the Anosov limit maps, for $d \in \{2,5\}$, of the associated representation $\rho:\pi_1S\rightarrow \Gtwosplit$. 
Then when $p \in \tilde{S}$ converges to $x \in \partial\Gamma$, we have  

    \begin{enumerate}[label=(\roman*)]
\item $T_p\to \xi^2(x)$
\item $N_p\to \xi^2(x)$
\item $\mathscr{L}_p \oplus T_p\oplus B_p\to \xi^5(x)$
\item $\mathscr{L}_p\oplus N_p\oplus B_p\to \xi^5(x)$
    \end{enumerate}
\end{theorem}

\begin{proof}Using Proposition \ref{prop:Limit maps}, it suffices once again to check transversality conditions. \medskip 

(iii) If a vector $u_\mathscr{L}+v_T+v_B$, with $u_\mathscr{L}\in \mathscr{L}_p$, $v_T\in T_p$, and $v_B\in B_p$, lies in $\xi^2(y)$ for some $y\in \partial\pi_1S$, then it is isotropic. Hence, $\lVert u_\mathscr{L}\rVert^2=\lVert v_T\rVert^2+\lVert v_B\rVert^2$, which is incompatible with the condition $\lVert v_T\rVert^2>\lVert u_\mathscr{L}\rVert^2+\lVert v_B\rVert^2$ from Corollary \ref{Cor:AlphaLimitSetTrapped}. \medskip 

(iv) Apply the same reasoning as in (iii), but exchanging $N$ and $T$. \medskip 

Now, (i) follows from (iv) and (ii) follows from (iii). 
\end{proof}
\section{Comparison of geometric structures}\label{Sec:Comparison}

In this section, we show that the image of the developing maps for the geometric structures constructed in \cite{DE26b}, for certain representations $\rho:\pi_1S\rightarrow \Gtwosplit$ associated to $\sigma$-bundles, is exactly the \cite{KLP18}-domain $\Omega_\sigma^{\rho}$. Moreover, the $(\Gtwosplit,\mathcal{F}_{\sigma})$-manifold ${\Fiber}_{\sigma}={\Fiber}^{\rho}_{\sigma}$ constructed in \cite{DE26b}\footnote{In \cite{DE26b} the manifold on which the geometric structure is constructed is $B_{\Psi_0}$, but we write ${\Fiber}$ in this paper.} is a (possibly trivial) finite cover of the quotient $M^{\rho}_{\sigma}=\rho(\pi_1S)\backslash \Omega^{\rho}_{\sigma}$. The proofs for $\alpha$ and $\beta$-bundles are analogous, with the key idea being to relate the image of a fiber by the developing map and the cones trapping the Anosov limit set.

\subsection{For \texorpdfstring{$\beta$}{b}-bundles}

Let $\rho:\pi_1S\to \Gtwo'$ be a $P_{\beta}$-Anosov representation. Recall from Proposition \ref{Prop:G2LimitMap} that $\rho$ is $P_{1,3}$-Anosov in $\SL(7,\R)$ and has limit maps $\xi^\beta=\xi^1:\partial \pi_1S\to \Ein^{2,3}$ and $\xi^3: \partial \pi_1S \rightarrow \Iso_3(\R^{3,4})$ that are related by $\xi^3(x) = \Ann(\xi^{\beta}(x))$. 

We introduce the following thickening in our $\Gtwosplit$-setting.
\begin{definition}[Thickening]
Given $\ell \in \Ein^{2,3}$, we denote by $K_{\ell}$ the following thickening: 
\[ K_{\ell } = \{ \ell' \in \Ein^{2,3} \mid \ell \times \ell'=0\}.\]
Similarly, we define $K_{\Lambda} =\bigcup_{x \in \Lambda} K_x$, where $\Lambda \subset \Ein^{2,3}$. 
\end{definition}

By our above discussion, we have $K_{\xi^{\beta}(x)}=\mathbb{P}(\xi^3(x)) \cong \RP^2$ by Proposition \ref{Prop:Annihilators}. Thus, this notion of thickening is the same as in \cite{GW12,KLP18}. We then conclude the following. 

\begin{theorem}[{\cite{GW12, KLP18}}]
    The representation $\rho$ acts properly and cocompactly on the following open set in $\Ein^{2,3}$:
    \[\Omega_\beta^{\rho} = \Ein^{2,3} \setminus K_{\mathrm{im}(\xi^{\beta})}=\Ein^{2,3} \setminus \bigcup_{x\in \partial \pi_1S}K_{\xi^{\beta}(x)}.\]
\end{theorem}

Let us consider a representation $\rho:\pi_1S\rightarrow \Gtwosplit$ associated to a $\beta$-Hodge bundle. Recall that such a representation comes with associated vector bundles $\mathscr{L},T,N,B$ on $S$. In fact, the flat vector bundle $\V^\R = \tilde{S} \times_{\rho} \R^{3,4}$ on $S$ decomposes as a direct sum as $\V^{\R} = \mathscr{L} \oplus T \oplus N \oplus B$. In \cite{DE26b} we constructed a developing map for a  $(\Gtwosplit,\Ein^{2,3})$-structure on the fiber bundle $\Ein(\mathscr{L} \oplus T\oplus B)\to S$ with descended holonomy $\rho$. We 
write ${\Fiber}_{\beta} = \Ein(\mathscr{L} \oplus T\oplus B)$ for brevity. Let $\overline{{\Fiber}_{\beta}}$ be the corresponding fiber bundle over $\widetilde{S}$ achieved via pullback. The developing map is defined as follows:
\[ \dev:\overline{{\Fiber}_{\beta}}\to \Ein^{2,3}, \]
\[ [u_\mathscr{L}+u_N+v_B ]\longmapsto\left[ u_\mathscr{L}+u_N+\frac{2u_\mathscr{L}\times v_B-u_N\times v_B}{\sqrt{4\lVert u_\mathscr{L} \rVert^2+ \lVert u_N \rVert^2}}\right]. \]

To make sense of this formula, we lift the bundles $\mathcal{E}^{\R}, \mathscr{L}, T,N$ to $\widetilde{S}$ and we identify for $p\in \widetilde{S}$ the bundle $\mathcal{E}^\R_{\mid p}$ with $\R^{3,4}$ using the flat connection. 

Setting $u_N=0$ and $u_{\mathscr{L}}=0$, and recalling $T=N\times B$, one then notes the following.
\begin{remark}
\label{rem:CircleEinstein}
    In particular, $\Ein(\mathscr{L}_p\oplus B_p)$ and $\Ein(T_p\oplus N_p)$ lie in the developed image of the fiber over $p\in \widetilde{S}$.
\end{remark}

As defined, there is no a priori relation between the map $\dev$ above and the domain $\Omega^{\rho}_{\beta}$. 
The following technical lemma allows us to prove our developing map avoids the thickening of the limit set for $\beta$-bundles. 

\begin{theorem}[Avoiding the $\beta$-thickening]\label{thm:BetaAvoidsLimitSet}
Let $(\mathcal{E},\Phi)$ be a polystable $\beta$-Hodge bundle such that $\delta$ has at least one zero and $\rho:\pi_1S\rightarrow \Gtwosplit$ be the associated representation. Then:
\begin{enumerate}[label=(\roman*)]
    \item The developing map $\dev: \overline{{\Fiber}_{\beta}}\rightarrow \Ein^{2,3}$ has image contained in $\Omega_\beta^{\rho} = \Ein^{2,3} \setminus K_{\mathrm{im}(\xi^{\beta})}$. 
    \item For all $p\in \widetilde{S}$ and $v\in \T_p\widetilde{S}$, the developed fiber $\dev(\overline{\Fiber_{\beta}}|_{p})$ is disjoint from the thickening $K_\ell$ for any $\ell \in \Ein^{2,3}$ such that $\ell\in \hat{C}'(v)$.
\end{enumerate}
\end{theorem}

\begin{proof}
Let us first show that (ii) implies (i). Fix a point $p \in \tilde{S}$ and any point $[Z_0] \in \dev(\overline{{\Fiber}_\beta}|_{p})$. We want to show that $[Z_0]$ does not lie in the thickening of $\im(\xi^{\beta})$. That is, for any isotropic line $[W] \in \im(\xi^{\beta})$, we need to show that $Z_0 \times W \neq 0$. 

Choose $x \in \partial \pi_1S$. Then we can find a tangent vector $v \in \T^1_{p}\widetilde{S}$ pointing towards $x$. By Corollary \ref{Cor:BetaLimitSetTrapped}, the isotropic line $[W]:=\xi^{\beta}(x)$ lies in the cone $\hat{C}'(v)$. 

Hence, assuming (ii), $[Z_0] \notin K_{[W]}$. Since $x$ was arbitrary, $[Z_0]\notin K_{\Lambda}$ for $\Lambda = \im(\xi^{\beta})$, which concludes (i). \medskip 

The remainder of the proof verifies part (ii).
\medskip 

Take any $p \in \widetilde{S}$, $v\in \T^1_{p}\widetilde{S}$, and select a unit generator of $\mathscr{L}_p$, denoted $\mathbf{e}$. For the proof, we switch perspectives: to show the developed fiber $\mathfrak{F} = \dev(\overline{\Fiber}_{\beta}|_{p})$ avoids the thickening of $\hat{C}'(v) \cap \Ein^{2,3}$, we can equivalently show that the restricted cone $\hat{C}'(v) \cap \Ein^{2,3}$ avoids the thickening $K_{\mathfrak{F}}$. 

Choose any element $[Z]\in \Ein^{2,3}$ in the thickening of an element $[Z_0]\in \mathfrak{F}$. We may decompose $Z$ as $Z=Z_{\mathscr{L}}+Z_N+Z_T+Z_B$ with respect to the Frenet splitting of $\V^\R$.
Let $\lVert \cdot \rVert$ denote the $h$-norm of an element in $\V^\R$. By the coordinate description of the bilinear form $B'_v$ from Lemma \ref{Lem:NestedBeta}, we note that for any $X \in \V^\R$, 
\[ \lVert{X_\mathscr{L}}\rVert^2 + \lVert X_{T}\rVert^2 - 2 \lVert X_N \rVert^2 -2 \lVert X_B\rVert^2 \geq B'_v(X). \] 
Thus, to guarantee that $[Z] \notin \hat{C}'(v)$, it suffices to prove that: 
\begin{align}\label{betaSufficient}
 2\lVert Z_N \rVert^2+2\lVert Z_B \rVert^2\geq \lVert Z_L \rVert^2+\lVert Z_T \rVert^2. 
\end{align}
We will now verify this inequality.

Write $U = \mathscr{L}\oplus N$, so that $\V^{\R}=U\oplus U^\bot$. 
Let us suppose the representative $Z_0$ has $\lVert (Z_0)_{U}\rVert =1$. Then since $[Z_0] \in \mathfrak{F}$, by the proof of \cite[Lemma 4.15]{DE26b}, for some $x,y\in \R$ such that $x^2+y^2=1$ and some unit elements $\mathbf{n}\in N_p$ and $\mathbf{b}\in B_p$ and $\lambda^{-1}=\sqrt{1+3x^2}$, we can express $Z_0$ as 
\[ Z_0=x\mathbf{e}+y\mathbf{n}+\lambda\big(2 x\mathbf{e}\times \mathbf{b}- y\mathbf{n}\times \mathbf{b}\big).\]

Now, recall that $Z \in \V^\R_p$ satisfies $Z_0\times Z=0$, since $[Z] \in K_{[Z_0]}$. We can rescale the representative $Z$ by some nonzero real number, if needed, so that $\lVert Z_U \rVert=1=\lVert Z_{U^\bot} \rVert$, where $Z_{U}=Z_\mathscr{L}+Z_N$ and $Z_{U^{\bot}}=Z_T+Z_B$. After doing so, we can write
\[ Z_U = c\mathbf{e}+ d \cos(\theta) \mathbf{n}+d\sin(\theta) \mathbf{e}\times \mathbf{n},\]
for some real numbers $c,d$ such that $c^2+d^2=1$. 

So far, we have expressed $Z_U = Z_{\mathscr{L}}+ Z_N$ in coordinates. 
The condition $Z \in \Ann(Z_0)$ implies that $Z_{U^\bot}$ is determined by $Z_U$, see \cite[Proposition 2.8]{DE26b}. More precisely, we have the relation 
\[ Z_{U^\bot} = -(Z_0)_{U^\bot}[ (Z_0)_{U} Z_U].\] 
Be advised: juxtaposition means split octonion multiplication in $\Oct'$ and not the cross-product, a meaningful distinction for this calculation.  

Thus, we have: 
\[ Z_T+Z_B=Z_{U^\bot} =-(2\lambda x\mathbf{e}\times \mathbf{b}-\lambda y\mathbf{n}\times \mathbf{b}) \big( (x\mathbf{e}+y\mathbf{n})(c\mathbf{e}+d\cos(\theta) \mathbf{n}+d\sin(\theta)\mathbf{e}\times\mathbf{n})\big).\]

We first compute that 
\[ (Z_0)_UZ_U= \big[-xc-yd\cos(\theta) \big]1_{\Oct'}+yd\sin(\theta)\mathbf{e}-xd\sin(\theta)\mathbf{n}+(xd\cos(\theta)-cy)\mathbf{e}\times \mathbf{n}.\]
Now, let us decompose the product $-(Z_0)_{U^{\bot}}((Z_0)_UZ_U)$ into $8$ terms:
\begin{align*}
-(Z_0)_{U^\perp}((Z_0)_UZ_U) &=- \lambda y(xc+yd\cos(\theta)) \mathbf{n}\times \mathbf{b} \\ 
 &+\lambda y^2d\sin(\theta)(\mathbf{n}\times \mathbf{b})\times \mathbf{e} \\
 &-\lambda xyd\sin(\theta)(\mathbf{n}\times \mathbf{b})\times \mathbf{n} \\
 &+\lambda y (xyd\cos(\theta)-cy)(\mathbf{n}\times \mathbf{b})\times (\mathbf{e}\times \mathbf{n}) \\
 &+2\lambda x(cx+xyd\cos(\theta)) \mathbf{e}\times \mathbf{b}\\
 &-2\lambda xyd\sin(\theta)(\mathbf{e}\times \mathbf{b})\times \mathbf{e} \\
 &+2\lambda x^2d\sin(\theta)(\mathbf{e}\times \mathbf{b})\times \mathbf{n} \\
 &-2\lambda x(xd\cos(\theta)-cy)(\mathbf{e}\times \mathbf{b})\times (\mathbf{e}\times \mathbf{n}). 
\end{align*}

Define $\mathbf{t}:=\mathbf{n}\times \mathbf{b}$. Using the fact that $(\mathbf{e},\mathbf{n},\mathbf{b})$ are the same as the model Stiefel triple $(\mathbf{i},\mathbf{j},\mathbf{l})$ up to the $\Gtwosplit$-action by \cite[Proposition 2.3]{DE26b}, we obtain the following relations:
\[ \begin{cases}
    (\mathbf{n}\times \mathbf{b})\times \mathbf{e}= -\mathbf{e}\times \mathbf{t}, \\
    (\mathbf{n}\times \mathbf{b})\times \mathbf{n}=\mathbf{b}, \\
    (\mathbf{n}\times \mathbf{b})\times (\mathbf{e}\times \mathbf{n})= -\mathbf{e}\times \mathbf{b}, \\
    (\mathbf{e}\times \mathbf{b})\times \mathbf{e}= \mathbf{b},\\
    (\mathbf{e}\times \mathbf{b})\times \mathbf{n}=\mathbf{e}\times \mathbf{t},\\
    (\mathbf{e}\times \mathbf{b})\times (\mathbf{e}\times \mathbf{n})= \mathbf{t}.
\end{cases} \] 

After expanding the previous product for $Z_{U^\bot}$ and decomposing into $T$ and $B$ parts, we obtain 
\begin{align}
    Z_T&=\lambda\big(cxy-d\cos(\theta)(2x^2+y^2)\big)\mathbf{t}+\lambda d \sin(\theta)\big(2x^2-y^2\big)\mathbf{e}\times\mathbf{t}.\\
    Z_B&=-3\lambda xyd \sin(\theta)\,\mathbf{b}+\lambda \big(c(1+x^2)+xyd\cos(\theta) \big)\mathbf{e}\times\mathbf{b}, 
\end{align}

Now let us consider the following quantity $A$:
\[A:=\lambda^{-2}\big(2\lVert Z_N \rVert^2+2\lVert Z_B \rVert^2- \lVert Z_L \rVert^2-\lVert{Z_T} \rVert^2 \big)=\lambda^{-2}\big(2d^2+2\lVert Z_B \rVert^2-c^2-\lVert Z_T \rVert^2 \big).\]
As noted above, to show $Z \notin C'(v)$ and hence complete the proof, it suffices to see $A \geq 0$. 
\medskip 

First note that the transformation $(d,\theta)\mapsto (-d,\theta+\pi)$ leaves $Z_T$, $Z_B$ and $A$ invariant. Hence one can assume that $cdxy\geq 0$. A calculation shows that, while only $\theta$ varies, the minimum value of $\lVert Z_B \rVert$ and the maximum value of $\lVert Z_T\rVert$ is attained at $\theta=\pi$. 

When $\theta=\pi$ we obtain:
\begin{align}
Z_T&=\lambda\big(cxy+d(1+x^2)\big)\mathbf{t}.\\
Z_B&=\lambda \big(c(1+x^2)-xyd\big)\mathbf{e}\times\mathbf{b}.
\end{align}
Hence setting $\theta=\pi$ when $cdxy\geq 0$ gives a lower bound on $A$: 
\begin{align*}
    A&\geq 2d^2(1+3x^2)+2c^2(1+x^2)^2+2d^2(xy)^2+4cd(1+x^2)(xy) \\ 
& -c^2(1+3x^2)-c^2(xy)^2-d^2(1+x^2)^2+2cd(xy)(1+x^2).
\end{align*}
We now regroup the terms in the lower bound, writing $A\geq c^2A_1+2cdA_2+d^2A_3$, where 
\begin{align*}
    A_1&=2(1+x^2)^2-(xy)^2-(1+3x^2)=1+3x^4, \\
 A_2&=2(1+x^2)(xy)+(1+x^2)(xy), \\
 A_3&=2(1+3x^2)+2(xy)^2-(1+x^2)^2=1+6x^2-3x^4.
\end{align*}

Observe that $A_1,A_3\geq 0$ for $0\leq x\leq 1$. Thus, it suffices to check that $\Delta=A_1A_3-A_2^2\geq 0$ to prove $A \geq 0$. Now, we conclude:
\begin{align*}
    \Delta&=(1+3x^4)(1+6x^2-3x^4)-9x^2y^2(1+x^2)^2=1-3x^2-9x^4+27x^6 \\
 &=(3x^2-1)^2(3x^2+1)\ge 0. \qedhere
\end{align*}
\end{proof}

\begin{corollary}\label{Cor:FiniteCoverBeta}
The map $\dev:\overline{{\Fiber}_{\beta}}\rightarrow \Omega_\beta^{\rho}$ is surjective and moreover is a finite covering. 
\end{corollary}

\begin{proof}
The claim follows from Theorem \ref{thm:BetaAvoidsLimitSet}, as $\dev:\overline{{\Fiber}_{\beta}}\rightarrow \Omega_\beta^{\rho}$ descends to a proper local diffeomorphism $\overline{\dev}:\Fiber_{\beta} \rightarrow M_{\beta}^{\rho}= \rho(\pi_1S)\backslash \Omega_{\beta}^{\rho}$, which must be a covering map. Hence, $\dev$ upstairs is also a covering map.  
For further details, see the similar argument in \cite[Theorem 4.25]{DE26b}.
\end{proof}

\subsection{For \texorpdfstring{$\alpha$}{a}-bundles}

We first recall the relevant domain of discontinuity. 
Given $\omega\in \Pho^{\times}$, we denote by $K_{\omega}\subset\Pho^\times$ the set of annihilator photons $\omega'$ such that $\omega \, \bot \, \omega'$. 
We proved in \cite[Proposition 5.10]{DE26b} that $K_{\omega}$ is topologically a projective plane. 

\begin{theorem}[{\cite{KLP18}}]
Let $\rho:\pi_1S\to \Gtwosplit$ be a $P_{\alpha}$-Anosov representation with limit map $\xi^\alpha:\partial \pi_1S\to \Pho^{\times}$. The representation $\rho$ acts properly and cocompactly on the following open set:
\begin{align}\label{DOD-alpha}
    \Omega_\alpha^{\rho} = \Pho^{\times} \setminus K_{\mathrm{im}(\xi^{\alpha})}=\Pho^{\times} \setminus \bigcup_{x\in \partial \pi_1S}K_{\xi^{\alpha}(x)}.
\end{align}

\end{theorem}
From the \cite{KLP18}-perspective, $K_{\omega}$ is a $\frac{\pi}{2}$-neighborhood of $\omega$ in $\Pho^{\times}$, seen in the visual boundary ${\vis}\X_{\Gtwosplit}$, with respect to the Tits angle metric.
For further details, see \cite[Section 5.2]{DE26b}. \medskip 

Let us consider a representation $\rho:\pi_1S\rightarrow \Gtwosplit$ associated to an $\alpha$-bundle $(\V,\Phi)$. Recall from \S\ref{Subsec:SigmaBundles} that the flat bundle $\V^{\R}\cong \tilde{S}\times_{\rho}\R^{3,4}$ of such a representation comes with associated vector bundles $\mathscr{L},T,N,B$ on $S$ giving a Frenet splitting of $\V^{\R}$. In \cite{DE26b} we constructed a developing map for a $(\Gtwosplit,\Pho^{\times})$-structure  with descended holonomy $\rho$ on the following $(\RP^2\times \sphere^1)$-fiber bundle:
\begin{align}\label{Malpha}
     {\Fiber}_{\alpha}:=\Gr_2(\mathscr{L}\oplus N)\oplus \Hom_\C(N,B)\to S. 
\end{align}
Let $\overline{{\Fiber}_{\alpha}}$ be the corresponding bundle over $\widetilde{S}$, namely the pullback of ${\Fiber}_{\alpha}$ under the universal cover $\pi:\tilde{S} \rightarrow S$. We constructed a developing map\footnote{With this model our developing map is only a local homeomorphism. Indeed, the preferred identification of $\mathscr{M}_{\alpha}$ with $B_{\Psi_0}$ in \cite{DE26b} is only a homeomorphism. The map:
$\dev_{honest}: \overline{B}_{\Psi_0} \rightarrow \Pho^{\times}$ is smooth and a local diffeomorphism.}
\[ \dev:\overline{{\Fiber}_{\alpha}}\to \Pho^{\times}. \]
The description of this map is explicit but lengthy. The only facts we will use about this developing map are as follows:

\begin{lemma}\label{Lem:OrthogonalPhoton}
Let $\omega \in \dev\left(\overline{{\Fiber}_{\alpha}}_{\mid p}\right)$ for some $p\in \widetilde{S}$. Then there is an isotropic line $\ell\in \Ein(N_p\oplus B_p)$ such that $\ell \subset\omega$. Moreover, $\Pho^{\times}(N_p\oplus B_p)$ lies in the image of the developed fiber over $p$.
\end{lemma}

\begin{proof}
    This is a consequence of \cite[Section 5.2]{DE26b}, see especially \cite[Lemma 5.16]{DE26b}.
\end{proof}
Now, as in the $\Ein^{2,3}$-case, the map $\dev:\overline{{\Fiber}_\sigma}\rightarrow \Pho^\times$ has no a priori relation to $\Omega^{\rho}_{\alpha}$. 
The following technical lemma allows us to prove that $\dev$ avoids the thickening of the limit set for $\alpha$-bundles. 

\begin{theorem}[Avoiding the $\alpha$-thickening]\label{thm:AlphaAvoidsLimitSet}
Let $(\mathcal{E},\Phi)$ be a polystable $\alpha$-Hodge bundle with associated representation $\rho:\pi_1S\rightarrow \Gtwosplit$. Then the developing map $\dev: \overline{{\Fiber}_{\alpha}} \rightarrow \Pho^{\times}$ has image contained in $\Omega_\alpha^{\rho}$. 
\end{theorem}

\begin{proof}
Fix any point $p \in \tilde{S}$ and any tangent vector $v \in \T^1_p\tilde{S}$. Let $(\mathscr{L},T,N,B)$ be the Frenet splitting of $\V^{\R}$ at $p$. Suppose that $\omega \in \Pho^\times$ is any annihilator photon such that $B_{v}|_{\omega} >0$, meaning $\mathbb{P}(\omega)\subset \hat{C}(v)$.

We will prove the following:\medskip 

\textbf{Claim}: For any isotropic line $\ell \in \Ein(N\oplus B)$, the 2-plane $\omega$ is not orthogonal to $\ell$. \medskip 

Most of our work is to verify this claim, from which we will easily deduce the lemma. \medskip

\emph{Proof of claim:} Let $\ell \in \Ein(N\oplus B)$ be arbitrary. To begin, we introduce some notation. 
First, recall that the $\pm 1$-eigenspaces $\mathbb{E}_{\pm 1}$ of $\psi = \Psi_0(v)$ form a pair of transverse annihilator photons $\omega_{\pm} \in \Pho^\times(N\oplus B)$. The kernel of $\psi$, denoted $\mathbb{E}_0$, is given by $\mathbb{E}_0=\mathscr{L}\oplus B$. \medskip 

Next, we observe a structural fact about $\mathbb{E}_1$ and $\mathbb{E}_{-1}$. Since $\psi$ exchanges $T$ and $N$, 
\begin{align}\label{Psi0Eigenvectors}
    (n+t) \in \mathbb{E}_{1} \iff (n-t) \in \mathbb{E}_{-1}.
\end{align}
We shall invoke this observation later.

Recall from Definition \ref{Defn:ConeAlpha} that if we introduce $h = q|_{U} \oplus (-q|_{U^\bot})$, the Euclidean metric associated to $U = \mathscr{L}\oplus N$, then $B_v$ is described by 
\[ B_{v} = h|_{\mathbb{E}_1} \oplus (-h|_{\mathbb{E}_{0}} ) \oplus (-h|_{\mathbb{E}_{-1}}).\]
Since $\mathbb{E}_1$, $\mathbb{E}_0$, $\mathbb{E}_{-1}$ are $h$-orthogonal, any 2-plane $\omega$ such that $B_v|_{\omega} > 0$ must then be a graph of a linear map $\chi: \mathbb{E}_1 \rightarrow \mathbb{E}_0\oplus \mathbb{E}_{-1}$ that is strictly $1$-Lipschitz for $h$. That is, $h(v,v)> h(\chi(v),\chi(v))$.

By hypothesis on $\ell$, we may write 
\[ \ell = [n+b], \;\text{for some} \; n \in Q_+(N), \;b\in Q_-(B).\]
Our goal is to show $\ell$ cannot be orthogonal to $\omega$. That is, we want to show there is some element $w \in \omega $ such that $\langle n+b, w \rangle_q \neq 0$. We introduce a number of notations to aid in this endeavor. 

Without loss of generality, we suppose that the Frenet splitting of $\V^{\R}_{\mid p}$ is as follows: 
\[ \mathscr{L} = \R\{ \mathbf{i}\}, \;T = \spann \langle \mathbf{l},\mathbf{li}\rangle, \;N = \spann \langle \mathbf{j},\mathbf{k}\rangle, \;B = \spann \langle \mathbf{lj},\mathbf{lk}\rangle.\]
We now choose bases of $\mathbb{E}_{+1}$, $\mathbb{E}_0$,  $\mathbb{E}_{-1}$ adapted to the situation. Since $\mathbb{E}_1$, $\mathbb{E}_{-1} \in  \Pho^\times(N \oplus T)$, each of these photons is the graph of an anti-isometry from $N$ to $T$. Let us select some element 
$(n+t) \in \mathbb{E}_1$, meaning $t \in Q_-(T)$. 
Then since $\mathbb{E}_1 \in \Pho^\times(N\oplus T)$, we have 
\[ \mathbb{E}_1 = \spann_{\R} \{ n+t, \mathbf{i}n+\mathbf{i}t\}.\]
Using \eqref{Psi0Eigenvectors}, we conclude that $\mathbb{E}_{-1}$ is completely determined:
\[ \mathbb{E}_{-1} = \spann_{\R} \{ n-t, \mathbf{i}n-\mathbf{i}t\}.\]
Next, we use the basis $(\mathbf{i},b,\mathbf{i}b)$ for $\mathbb{E}_0=\mathscr{L} \oplus B$. 

Consider the element $w \in \omega$ given by $w= (n+t)+\chi(n+t)$. 
Writing all of the constraints on $w$ out in the fixed basis will lead us to an elementary calculation from which we can conclude. Now, first express $w$ as follows:
\[ w= w_{+1} +w_{0}+w_{-1} \in \mathbb{E}_{1} \oplus \mathbb{E}_{0} \oplus \mathbb{E}_{-1}. \]
In other words, $w_{0}+w_{-1}= \chi(n+t)$ and $w_{+1}=n+t$.  
We now expand $w$ in coordinates. Namely, 
\[ w_0 + w_{-1} = c_1 (n-t) + c_2(\mathbf{i}n-\mathbf{i}t)+c_3\mathbf{i}+c_4b+c_5\mathbf{i}b.\]
The quantity of interest, the $q_{3,4}$-pairing of $w$ and $n+b$, is given by 
\[ \langle n+b, w\rangle_{q}= \langle n+b, w_{+1}+w_0+w_{-1} \rangle_{q}= 1+\langle n,w_0\rangle_q+\langle b, w_{-1}\rangle_q = 1+c_1-c_4.\]

We now consider the constraints on the coordinates $(c_i)_{i=1}^5$. 
First, we have $q(w)=0$. 
Using $q(n-t)=q(n+t)=0$ and $q(n-t,n+t)=2$, we obtain 
\begin{align}\label{Constraint1}
     0 = q(w,w) = 2q(w_1, w_{-1})+q(w_0,w_0) = 4c_1+c_3^2-c_4^2-c_5^2.
\end{align}
Next, using that $\chi$ is strictly 1-Lipschitz, we find 
\begin{align}\label{Constraint2}
 2c_1^2+2c_2^2+c_3^2+c_4^2+c_5^2=h(w_0+w_{-1}, w_0+w_{-1})<2.
\end{align}

Finally, we proceed to show $\langle n+b, w\rangle_q > 0$. To this end, we consider the minimum of $1+c_1-c_4$ relative to the constraints \eqref{Constraint1}, \eqref{Constraint2}. 
Clearly, the worst case possible occurs when $c_5=0=c_2$. In this case, $q(w,w)=0$ gives the 
relation $c_3^2=c_4^2-4c_1$. Substituting this back into $h(w_0+w_{-1}) <2$,  
\[ 2c_1^2-4c_1+2c_4^2< 2,\]
\[ 2(c_1-1)^2+2c_4^2<4\,.\]
From this, we can conclude:
\[1+c_1-c_4=2-(1-c_1)-c_4\geq 2-\sqrt{2(1-c_1)^2+2c_4^2}> 0.\]
Thus, the claim has finally been proven. $\square_{\text{claim}}$ \medskip 

We now deduce the lemma. 
Fix any point $p \in \tilde{\Sigma}$. Select any annihilator photon $\omega_{\mathrm{fib}} \in \dev(\overline{{\Fiber}_\sigma}_{|_p})$. By Lemma \ref{Lem:OrthogonalPhoton}, there is some line $\ell \subset \omega_{\mathrm{fib}}$ of the form $\ell \in \Ein(N_p\oplus B_p)$. 

Theorem \ref{thm:AlphaAnosov} asserts that $\rho$ is $P_{\alpha}$-Anosov. 
Choose any $x \in \partial \pi_1S$. We can select a tangent vector $v \in \T^1_p\tilde{\Sigma}$ such that $v$ points towards $x$. Now, let us write $\omega_{\lim} = \xi^{\alpha}(x) $. 
By Lemma \ref{Lem:NestedAlpha} and Proposition \ref{prop:Nestednessproperties}, we conclude that the limiting annihilator photon $\omega_{\lim}$ lies in the cone $C(v)$ of positive lines for the bilinear form $B_{v}$ from \eqref{Alpha_2cone}. 

Finally, the claim asserts that $\ell$ is not orthogonal to $\omega_{\lim}$. Hence, we conclude that $\omega_{\mathrm{fib}}$ is not orthogonal to $\omega_{\lim}$. 
As $\omega_{\lim} \in \image(\xi^{\alpha})$ was arbitrary, by definition \eqref{DOD-alpha} of $\Omega^{\rho}_{\alpha}$, this means $\omega_{\mathrm{fib}} \in \Omega^{\rho}_{\alpha}$. Hence, $\dev$ has image in $\Omega^{\rho}_{\alpha}$. 
\end{proof}

\begin{corollary}
The map $\dev: \overline{\mathscr{M}_{\alpha}}\rightarrow \Omega^{\rho}_{\alpha}$ is surjective and moreover is a finite covering.
\end{corollary}

\begin{proof}
The claim follows from Theorem \ref{thm:AlphaAvoidsLimitSet} by the same reasoning as in Corollary \ref{Cor:FiniteCoverBeta}. 
\end{proof}

\subsection{Linking number in \texorpdfstring{$\mathbb{S}^2\times \mathbb{S}^3$}{S2xS3}}

We recall from \cite{DE26b} that both $\Ein^{2,3}$ and $\Pho^\times $ are diffeomorphic to $\RP^2 \times \sphere^3$. Thus, the universal cover of each flag manifold is diffeomorphic to $\mathbb{S}^2\times \mathbb{S}^3$. Now, for $\alpha$-Hodge bundles and $\beta$-Hodge bundles, the Anosov limit maps $\xi$ takes values in the same flag manifold as the developing map $\dev$. This allows us to compare the relative position of a developed fiber and of the limit set. For the remainder of this section, we work exclusively in the double cover $\widehat{\mathcal{F}}_\sigma$ of $\mathcal{F}_{\sigma} = \Gtwosplit/P_{\sigma}$. 

Now, a bit more precisely, using the $\Gtwosplit$ geometry, we define a notion of linking number for a pair of disjoint loops in $\widehat{\mathcal{F}}_\sigma$. We show that the limit map $\xi^{\sigma}$ lifts to a loop $\hat{\xi}^{\sigma}$ in $\widehat{\mathcal{F}}_{\sigma}$ that is simply linked with a loop $\gamma$ generating the fundamental group of a developed fiber, i.e., $\hat{\xi}^{\sigma}$ and $\gamma$ have linking number $\pm1$. A consequence that we obtain in  Theorem \ref{thm:DomainPi1} is that the fundamental group of the domain of discontinuity is a semi-direct product of the fundamental group $\Z$ of the fiber with a finite group $N$. This allows us to get enough control on the fundamental group of the quotient $\hat{M}_{\sigma}$ of the domain of discontinuity $\hat{\Omega}_{\sigma}$ in Corollary \ref{Cor:FundamentalGroupSES}. 
\medskip

We begin by introducing the following notion:

\begin{definition}\label{Defn:LinkingNumber}
    Let $\gamma:\mathbb{S}^1\to \mathbb{S}^2\times \mathbb{S}^3$ and $\tau:\sphere^1\times \sphere^2\to \mathbb{S}^2\times \mathbb{S}^3$ be continuous maps such that 
    \begin{itemize}
        \item $\gamma \pitchfork \tau$ (i.e., their images are disjoint),
        \item $\tau$ defines a trivial homology class in $H_3(\mathbb{S}^2\times \mathbb{S}^3)$.
    \end{itemize}
     We define the \textbf{linking number} $\ln(\gamma, \tau)\in \mathbb{Z}$ of $\gamma$ and $\tau$ as the intersection number between $\tau$ and any map $\tilde{\gamma}:\mathbb{D}^2\to \mathbb{S}^2\times \mathbb{S}^3$ such that $\tilde{\gamma}_{|\partial \mathbb{D}^2}=\gamma$. 
\end{definition}

We now remark why Definition \ref{Defn:LinkingNumber} is well-defined. Since $\sphere^2\times \sphere^3$ is simply connected, an extension $\tilde{\gamma}:\mathbb{D}^2\to \mathbb{S}^2\times \mathbb{S}^3$ exists. To make sense of the intersection number between $\tau$ and any such disk $\tilde{\gamma}$, one can homotope both $\tau$ and $\tilde{\gamma}$, the latter rel boundary, until the two  maps are smooth (away from $\partial \mathbb{D}^2$ for $\tilde{\gamma}$), and intersect transversely \cite[$\S$2.3]{GP10}. Then one can count the sum of signed intersections between $\tau$ and $\tilde{\gamma}$. This number is independent of the smooth disk chosen. Indeed, if we pick a different smooth disk $\tilde{\gamma}'$ bounding $\gamma$, then gluing $\tilde{\gamma}$ and $\tilde{\gamma}'$, after changing the orientation of one disk, defines a map $f(\tilde{\gamma}, \tilde{\gamma}'):\mathbb{S}^2\to \mathbb{S}^2\times \mathbb{S}^3$, and its homological intersection number with $\tau$ must vanish, so the sum of signed transverse intersections between $f$ and $\tau$ is equal to zero. 

In particular, if $\gamma$ and $\tau$ are smooth embeddings and one can extend $\gamma$ to a smooth embedding $\tilde{\gamma}$ that transversely intersects $\tau$ exactly once, then $\ln (\gamma, \tau)=\pm 1$. 
\medskip

We now verify that the Anosov limit maps lift to the double covers. 

\begin{lemma}[Lifting the limit map]
    Let $\rho:\pi_1S\to \Gtwosplit$ be the representation associated to a stable $\beta$-Hodge bundle or a polystable $\alpha$-bundle. Denote $\sigma \in \{\alpha,\beta\}$. Then the Anosov limit map $\xi^{\sigma}:\partial \pi_1S\to \mathcal{F}_{\sigma}$ lifts to a map $\hat{\xi}^{\sigma}: \mathbb{S}^1\simeq \partial \pi_1S\to \widehat{\mathcal{F}}_{\sigma} \simeq  \mathbb{S}^2\times \mathbb{S}^3$.
    \end{lemma}

\begin{proof}
Fix any point $p \in \tilde{S}$ and denote $U:= (\mathscr{L}_p\oplus N_p) \in \X_{\Gtwosplit}$. 

\textbf{Case 1}: ($\sigma= \beta$). 
Recall that associated to $U$, we obtain an explicit diffeomorphism $F_{U}: \RP^2\times \sphere^3 \cong \Gr_1(U)\times Q_-(U^\bot) \rightarrow \Ein^{2,3}$ by $F_U([x], y) = [x+ x\times y]$ \cite[Proposition 4.2]{DE26b}. Now, $\pi_1(\Ein^{2,3}) \cong \pi_1\RP^2$. Using  $F_U$ to identify $\Ein^{2,3} \cong \RP^2 \times \sphere^3$, it then suffices to prove that for $\Lambda: = \im(\xi^{\beta})$, the projection $\pi_{\RP^2}(F_U^{-1}(\Lambda))$ is nullhomotopic. However, this follows immediately by Corollary \ref{Cor:BetaLimitAffineChart}, which says that the projected limit set $\pi_{\RP^2}(F_U^{-1}(\Lambda))$ is trapped in an affine chart of $\Gr_1(U)$. Hence, $\xi^{\beta}$ lifts to the double cover $\widehat{\Ein}^{\, 2,3} = \widehat{\mathcal{F}}_{\beta}$. 

\textbf{Case 2}: ($\sigma= \alpha$). We proceed similarly. From \cite[Section 5.1]{DE26b}, there is a smoothly trivial $\sphere^3$-fiber bundle $\pi_{U}: \Pho^\times \rightarrow \Gr_2(U)$ given by the orthogonal projection of $\omega$ to $U$. Hence, we may take a diffeomorphism $\Pho^\times \cong \RP^2 \times \sphere^3$ such that the projection to $\RP^2$ is given by $\pi_{U}$.
As in Case 1, we must prove that the projected limit set $\pi_{U}(\im(\xi^{\alpha}))$ is nullhomotopic. Let us take $\omega \in \im(\xi^{\alpha})$ and write 
$\omega = \spann\langle u+z, v+(uv)z\rangle $ for a triple $(u,v,z) \in V_{(+,+,-)}(\imoct)$ (see \cite[Corollary 2.9(2)]{DE26b}). Since $\spann \langle u,v \rangle \cap N_p \neq \emptyset$, we may assume without loss of generality that $u \in N$. Now, Corollary \ref{Cor:AlphaLimitSetTrapped} gives the inequality
\[  \lVert v_N\rVert+\lVert((uv)z)_T\rVert > \lVert v_{\mathscr{L}}\rVert+\lVert((uv)z)_B\rVert, \] 
where each individual term is at most one. In particular, $v_N=0$ implies $||v_{\mathscr{L}} ||= 1$, which is impossible. Hence, $v_{N} \neq 0$, which means $\pi_{U}(\omega) = \spann \langle u,v\rangle$ is transverse to $\mathscr{L}_p$. We conclude the projection $\pi_U(\im(\xi^{\alpha}))$ is nullhomotopic and thus $\xi^{\alpha}$ lifts to $\widehat{\mathcal{F}}_{\alpha}$. \end{proof}

We now discuss thickenings of lifts. Let $\hat{\xi}^{\sigma}: \sphere^1\rightarrow \widehat{\mathcal{F}}_{\sigma}$ be the lifted limit map.\medskip 

First, we consider $\sigma =\beta$. Thickenings upstairs in the double cover $\widehat{\mathcal{F}}_{\beta}$ are similar to those downstairs. Identify $\widehat{\mathcal{F}}_{\beta}$ with the space of isotropic rays in $\R^{3,4}$. 
Given a point $\ell \in \widehat{\mathcal{F}}_{\beta}$, define the thickening $K_{\ell}$ via annihilators:
\begin{align}\label{eq:BetaThickening}
     K_{\ell} = \{ \ell' \in \widehat{\mathcal{F}}_{\beta} \;\mid \;\ell'\times \ell = 0\}.
\end{align}

It will be useful to view the construction of the `thickening' of a point $x \in \widehat{\mathcal{F}}_{\beta}$ as defining a fiber bundle over the flag manifold, as the following result describes. 

\begin{proposition}[$\beta$-thickenings, upstairs]\label{Prop:BetaThickenings}
For any $\ell \in \widehat{\mathcal{F}}_{\beta}$, we have $K_\ell \cong \sphere^2$. 

Additionally, the space
\[ \mathscr{K} = \{ (\ell,\ell') \in \widehat{\mathcal{F}}_\beta \times \widehat{\mathcal{F}}_{\beta} \mid  \ell\times \ell' = 0\}.\]
is a trivial $\sphere^2$-fiber bundle over 
$ \widehat{\mathcal{F}}_{\beta} $ via the map $(\ell,\ell') \mapsto \ell$. 
\end{proposition}

\begin{proof}
By definition, we have $K_{\ell} \cong \sphere(\Ann(\ell))$. 
Following similar reasoning to \cite[Proposition 4.24]{DE26b}, the vector bundle $\Ann \rightarrow \widehat{\mathcal{F}}_{\beta}$ by $\Ann|_{\ell} =\Ann(\ell) $ is trivial and of rank three. 
Thus, the total space of $\sphere(\Ann) \cong \widehat{\mathcal{F}}_{\beta} \times \sphere^2$ is diffeomorphic to the total space of $\mathscr{K}$. 
\end{proof}

For $\sigma = \alpha$, a similar story applies. The double cover $\widehat{\mathcal{F}}_{\alpha}$ is the space of \emph{oriented} annihilator photons, which can be identified with the total space of the pullback of the fiber bundle $\Pho^\times \rightarrow \Gr_{2}(U)$ under the double cover $\Gr_2^{or}(U) \rightarrow \Gr_2(U)$. For $\hat{\omega} \in\widehat{\mathcal{F}}_{\alpha}$, the thickening is again given by
 \[ K_{\hat{\omega}} =  \{ \omega' \in \widehat{\mathcal{F}}_{\alpha} \mid \omega \perp \omega'\}.\]
The thickening $K_{\omega}$ is also topologically a 2-sphere. For clarity, we summarize the analogous results in the $\widehat{\mathcal{F}}_{\alpha}$ case, which are a bit more involved. 

To begin, we must recall a topological model for $\widehat{\mathcal{F}}_{\alpha}$ from \cite[Lemma 5.6]{DE26b}. 

\begin{lemma}[A geometric model for $\widehat{\mathcal{F}}_{\alpha}$]
\label{lem:ModelPho}
Fix any point $U \in \X_{\Gtwosplit}$ and a splitting $U = \mathscr{L} \oplus N$. Then there is a corresponding homeomorphism 
\[ F: \Gr_2^{or}(U) \times \sphere(\Hom_{\C}(N, U^\bot))\rightarrow \widehat{\Pho}^\times, \qquad (W, \delta) \mapsto \spann\langle u+ \delta(u), v+(uv)\delta(u)\rangle , \]
where $(u,v)$ is an orthonormal basis for $W$ and $u \in Q_+(N\cap W)$ is arbitrary. 
\end{lemma}
Here, we endow the image photon $F(W,\delta)$ with the orientation induced by that of $W$. 

In fact, the map $F$ in Lemma \ref{lem:ModelPho} is not smooth. Nevertheless, it suffices for the present topological considerations. By using this model, we can validate that the universal bundle of thickenings $\mathscr{K} \rightarrow \widehat{\mathcal{F}}_{\alpha}$ is also topologically trivial. 

\begin{proposition}[$\alpha$-thickenings, upstairs]\label{Prop:AlphaThickenings}
For any $\omega \in \widehat{\mathcal{F}}_{\alpha}$, we have $K_{\omega} \cong \sphere^2$. 

Additionally, the space $\mathscr{K}$ given by 
\[ \mathscr{K} = \{ (\omega,\omega') \in \widehat{\mathcal{F}}_\alpha \times \widehat{\mathcal{F}}_{\alpha} \mid  \omega \,\bot\, \omega' \}\]
is topologically a trivial $\sphere^2$-fiber bundle over 
$ \widehat{\mathcal{F}}_{\alpha} $ via the map $(\omega,\omega') \mapsto \omega$. 
\end{proposition}

\begin{proof}
Fix a Frenet splitting $\R^{3,4}=(\mathscr{L}\oplus N)\oplus (T\oplus B)$ as in Definition \ref{defn:FrenetFrame} and form the associated model $F$ 
from Lemma \ref{lem:ModelPho}. Write $U = \mathscr{L}\oplus N$. Now, choose any oriented annihilator photon $\omega = F(W, \delta)$. The condition that $\omega\, \bot\, \omega'=0$ is equivalent to demanding $\dim(\omega + \omega' )\leq 3$ by \cite[Proposition 5.8]{DE26b}. Let us then fix any line $[u+\delta(u)] \subset \omega$, where $u \in N$. 

The model $F$ is very convenient to describe thickenings. In fact, for $F(W',\delta') =\omega' \in \widehat{\mathcal{F}}_{\alpha}$, we have $\omega\,\bot \, \omega' =0$ if and only if $\delta = \delta'$, see the proof of Proposition \cite[{Proposition 5.11}]{DE26}.

From this fact, we conclude that $F^{-1}(K_{\omega}) = \Gr_2^{or}(U)\times \{ \delta\} \cong \sphere^2$. In particular, this implies that, via $F$, there is a canonical identification between $K_{\omega}$ and $\Gr_2^{or}(U)\cong \sphere^2$. In particular, there is a fiber bundle homeomorphism between $\mathscr{K}$ and $\widehat{\mathcal{F}}_{\alpha} \times \sphere^2$. 
\end{proof}

We conclude the topology of the thickening of any \emph{transverse} loop $\gamma$. Recall that this means if $x \neq x'$ in $\sphere^1$, then the flags $\gamma(x)$ and $\gamma(x')$ are transverse in the target flag manifold. 

\begin{corollary}[Bundle of thickenings]\label{Cor:BundleThickenings}
Let $\sigma \in \{\alpha, \beta\}$ and $\gamma: \sphere^1 \rightarrow \widehat{\mathcal{F}}_{\sigma}$ be continuous and transverse. Then $K_{\im(\gamma)}$ is homeomorphic to $\sphere^1\times \sphere^2$. 
\end{corollary}

\begin{proof}
We have $\gamma^*\mathscr{K} \cong \sphere^1\times \sphere^2$ by Proposition \ref{Prop:BetaThickenings} and \ref{Prop:AlphaThickenings}, respectively. The map $\gamma^*\mathscr{K} \rightarrow \widehat{\mathcal{F}}_{\sigma}$ by $(x,y) \mapsto y$ is injective since $\gamma$ is transverse. The result follows.
\end{proof}

For the following proof, given a curve $\gamma: \sphere^1 \rightarrow\widehat{\mathcal{F}}_{\sigma}$, we can pull back the fiber bundle $\mathscr{K}\rightarrow\widehat{\mathcal{F}}_{\sigma}$ to obtain an $\sphere^2$-fiber bundle over the circle, regardless of whether $\gamma$ is transverse or not. 
We can then understand homological properties of any bundle of thickenings by studying a model case. 

\begin{lemma}[Homologically trivial bundle of thickenings]
\label{Lem:BundleThickeningsH3}
Let $\sigma \in \{\alpha,\beta\}$ and $\gamma: \sphere^1\rightarrow \widehat{\mathcal{F}}_{\sigma}$ be a continuous curve. 
The inclusion map $\tau:\gamma^*\mathscr{K}\rightarrow \widehat{\mathcal{F}}_{\sigma}$ induces a trivial map of third homology groups.
\end{lemma}

\begin{proof}
One can find a homotopy that contracts $\gamma$ to a constant map $\gamma_f:\mathbb{S}^1\to \lbrace f\rbrace \subset \widehat{\mathcal{F}}_\sigma$. The homotopy class of the corresponding inclusion map $\tau_f:\gamma_f^*\mathscr{K}\rightarrow \widehat{\mathcal{F}}_{\sigma}$ is equal to that of $\tau$.
However, $\tau_f$ factors through $\mathbb{S}^2$ and hence induces a trivial map on the level of $H_3$.
\end{proof}

We can now realize the main purpose of Definition \ref{Defn:LinkingNumber}: to give a notion of `linking' for two appropriate loops in the double covers $\widehat{\mathcal{F}}_{\sigma}$. 
\begin{definition}[Linking in oriented $\Gtwosplit$-flag manifolds]\label{Defn:LinkingG2Flags}
Let $\gamma_1, \gamma_2:\mathbb{S}^1 \rightarrow \widehat{\mathcal{F}}_{\sigma}$ be continuous maps such that 
$\gamma_1$ avoids the thickening of $\gamma_2$. Then we define the \textbf{linking number} of $\gamma_1$ and $\gamma_2$ by 
$\ln(\gamma_1, \gamma_2):=\ln (\gamma_1, \tau_2)$, where
$\tau_2$ is the inclusion of the bundle of thickenings $\gamma_2^*\mathscr{K}$. 
\end{definition}

\begin{remark}
    We do not prove that this linking number is symmetric in general, as we will not need such considerations going forward. 
\end{remark}

Note that by Corollary \ref{Cor:BundleThickenings} and Lemma \ref{Lem:BundleThickeningsH3}, the above notion is well-defined. \medskip 

For later, we compute the linking number in a simple model case. 
Here, we use the oriented Frenet splittings introduced in Definition \ref{defn:FrenetFrame}.

\begin{lemma}[Linking number, model case]\label{Lem:BetaFibers+LimitSet}
   Fix any oriented Frenet splitting $(e,T,N,B)$ for $\R^{3,4}$, denote
   $\hat{\mathscr{L}} = \R_+\{e\}$, and consider the inclusion maps:
    \[ \gamma_T: \widehat{\Ein}(\hat{\mathscr{L}}\oplus T) \rightarrow \widehat{\mathcal{F}}_{\beta}, \qquad \gamma_B: \widehat{\Ein}(\hat{\mathscr{L}}\oplus B) \rightarrow \widehat{\mathcal{F}}_{\beta}.\]
    Then $\ln(\gamma_{T}, \gamma_{B}) = \ln(\gamma_{B}, \gamma_{T}) =\pm1 $.
\end{lemma}

Note that $\widehat{\Ein}(\R^{1,2}) \cong \sphere^1\times \sphere^0$ is the disjoint union of two circles. 
This is why we consider $\widehat{\Ein}(\hat{\mathscr{L}}\oplus T)$ and $\widehat{\Ein}(\hat{\mathscr{L}}\oplus B) $, which are each just one circle. 

\begin{proof}
    Again, identify $\widehat{\Ein}^{\,2,3} \cong\widehat{\mathcal{F}}_{\beta}$. We use the notation $\widehat{\Ein}(W)=\{x \in \widehat{\Ein}^{\,2,3} \mid x \subseteq W\}$, where $W \subset \R^{3,4}$ is any subset.  
    First note that $K_{\widehat{\Ein}(\hat{\mathscr{L}}\oplus T)}$ and $\widehat{\Ein}(\hat{\mathscr{L}}\oplus B)$ are disjoint since $T\times B= N$. The entire argument is symmetric in $T$ and $B$, so we need only prove $\ln(\gamma_{T},\gamma_{B})=\pm 1$. 

    Write $U = \mathscr{L}\oplus N$ and $U^\bot = T \oplus B$. We shall use the diffeomorphic model $\widehat{\Ein}^{\,2,3}\cong Q_+(U)\times Q_-(U^\bot)$, writing each ray $x \in \widehat{\Ein}^{2,3}$ uniquely as $ x=\R_+\{u+v\}$, with $u\in Q_+(U), v\in Q_-(U^\bot)$. 
    Consequently, we may canonically identify 
    \begin{align}\label{TangSpaceEinHat}
         \T_{(u,v)}\widehat{\Ein}^{\,2,3} \simeq (u^\bot \cap P) \oplus (v^\bot \cap U^\bot).
    \end{align}
   
    Select a unit vector $b_0\in B$. We can fill the loop $\gamma_T$ by a smooth disk $\tilde{\gamma}_T: \mathbb{D}^2 \rightarrow \widehat{\Ein}^{\,2,3}$ as follows: 
    \begin{align*}
        \tilde{\gamma}_{T}:\lbrace t\in T:  \lVert t\rVert\leq 1\rbrace &\longrightarrow \widehat{\Ein}(\hat{\mathscr{L}}\oplus T\oplus \langle b_0\rangle) \\
     t&\longmapsto  \big(e,\,t+\sqrt{1-\lVert t \rVert^2}\,b_0 \big). 
    \end{align*}
    Let $\tau_B: \sphere^1\times \sphere^2 \rightarrow K_{\gamma_B}$ parametrize the thickening $K_{\gamma_B}$. 
    As noted earlier, $\tau_B$ and $\gamma_T$ do not intersect. 
    One can check that $\tilde{\gamma}_{T}$ and $\tau_B$ intersect exactly once, namely at $(e,b_0) \in \widehat{\Ein}^{2,3}$. 
    
    Moreover, the maps $\tilde{\gamma}_{T}$ and $\tau_B$ intersect transversely at this point, which is seen as follows. First, note that 
    $d\tilde{\gamma}_T|_{\mathbf{0}}= 0 \oplus T$ in the model 
    \eqref{TangSpaceEinHat}. Next, denote 
    $(\theta,y) \in \sphere^{1}\times \mathbb{S}^2$
    as the unique preimage of $(e,b_0)$ under $\tau_{B}$. 
    Then observe that $d\tau_{B}|_{(\theta,y)}$ has image that orthogonally projects surjectively onto $N \oplus \langle b_0\rangle^\bot \subset \T_{(e,b_0)}\widehat{\Ein}^{2,3}$. 
    Thus, $\image(d\tau_{B}|_{(\theta,y)}) \oplus \image (d\gamma_T|_{\mathbf{0}}) = \T_{(e,b_0)}\widehat{\Ein}^{2,3}$.  
\end{proof}

The analogous result for $\sigma =\alpha$ is as follows. 

\begin{lemma}\label{Lem:AlphaFibers+LimitSet}
   Fix any oriented Frenet splitting $(e,T,N,B)$ for $\R^{3,4}$ and  $\gamma_T, \gamma_B: \sphere^1 \rightarrow \widehat{\mathcal{F}}_{\alpha}$ be lifts of the loops $\Pho^{\times}(N \oplus T)$ and $\Pho^{\times}(N\oplus B)$, respectively. 
   Then $\ln(\gamma_T, \gamma_B) =\ln(\gamma_B, \gamma_T)=\pm1 $.
\end{lemma}

\begin{proof}
    First, we note that $K_{\Pho^{\times}(N\oplus T)}$ and  $\Pho^{\times}(N\oplus B)$ do not intersect. 
    Thus, if we consider $\tau_{B}: \sphere^1\times \sphere^2 \rightarrow \widehat{\mathcal{F}}_{\alpha}$ to parametrize the thickening $K_{\gamma_B}$, then $\tau_{B}$ and $\gamma_{T}$ are disjoint.  
    Again, the argument is symmetric in $T$ and $B$, and we shall choose to extend $\gamma_T$. 
    
    Choose a norm one conformal map $\phi_0: N\to B$ whose graph is an oriented annihilator photon, with orientation induced by $\mathcal{C}_{e}$, such that $\graph(\phi_0) \in \image(\gamma_B)$. 
    We can fill in the loop $\gamma_T$ by a disk $\tilde{\gamma}_{T}:\mathbb{D}^2 \rightarrow \widehat{\mathcal{F}}_{\alpha}$ as follows. Let $\Hom_{\C}(N,T)$ be the space of $\mathcal{C}_{e}$-linear maps. We consider the following map:
    \begin{align*}
          \tilde{\gamma}_{T}:\lbrace \phi \in \Hom_{\C}(N,T)\mid ||\phi||_{op} \leq 1\} &\longrightarrow \widehat{\Pho}^\times(N \oplus T \oplus B) , \\ 
         \phi &\longmapsto \graph\big(\phi+\sqrt{1-\lVert \phi \rVert^2}\,\phi_0 \big).
    \end{align*}
    We recall from \cite[Section 5.1]{DE26b} that $\Pho^\times(N\oplus T\oplus B)$ consists exactly of graphs of $\mathcal{C}_{e}$-linear holomorphic maps 
    $\phi: N\rightarrow T\oplus B$ of operator norm one. In particular, $\tilde{\gamma}_{T}$ is well-defined (again, using $\mathcal{C}_{e}$ to orient the annihilator photon).

    Now, two annihilator photons $\omega_1 \in \Pho^\times(N\oplus B)$ and $\omega_2 \in \Pho^\times(N \oplus T\oplus B)$ are orthogonal if and only if $\omega_2 =\omega_1$. Thus, $\tilde{\gamma}_{T}$ and $\tau_{B}$ intersect exactly once, namely at $\omega_0=\graph(\phi_0)$. Let us view this intersection in the model $\Gr_2^{or}(U) \times \sphere(\Hom_{\C}(N, U^\bot))$ from Lemma \ref{lem:ModelPho} by using the homeomorphism $F^{-1}$. From this perspective, we are intersecting 
    the 2-disk 
    \[ \lbrace N\rbrace \times D,\; \text{where}\; D = \big\{\phi+\sqrt{1-\lVert \phi \rVert^2}\phi_0 \mid \phi \in \Hom_{\C}(N,B) \big\}\subset    \sphere(\Hom_{\C}(N, U^\bot)),\]
    with the copy of $\sphere^1\times \sphere^2$ given by 
    $\Gr_2^{or}(U) \times \sphere(\Hom_{\C}(N,B)).$
    These two submanifolds evidently intersect transversely at $\omega_0=(N, \phi_0) \in \Gr_2^{or}(U) \times \sphere(\Hom_{\C}(N, U^\bot))$. 
    Note that applying $F$ doesn't change the linking number, which is purely topological. 
\end{proof}

Next, we discuss lifting the geometric structures from \cite{DE26b}.
In particular, we shall use a double cover 
$\hat{\mathscr{M}}_{\sigma}$ of $\mathscr{M}_{\sigma}$, 
which corresponds to lifting the geometric structure from 
$\mathcal{F}_{\sigma}$ to $\widehat{\mathcal{F}}_{\sigma}$.
For concreteness, we describe these fiber bundles  $\sphere^1 \times \sphere^2 \rightarrow \hat{\Fiber}_{\sigma} \rightarrow S$ now. 
Let $(\V,\Phi)$ be a $\sigma$-bundle and $\mathscr{L},T,N,B$ be the Frenet framing of $\V^{\R}$. Then we consider the following fiber bundles:  
\begin{align}\label{eq:DoubleCover}
    \begin{cases}
    \hat{\mathscr{M}}_{\beta} &= \widehat{\Ein}(\mathscr{L} \oplus N \oplus B) \cong \sphere(\mathscr{L}\oplus N)\oplus \sphere(B) ,\\
     \hat{\mathscr{M}}_{\alpha} &= \Gr_{2}^{or}(\mathscr{L}\oplus N)\oplus \sphere(\Hom_{\C}(N,B)). 
     \end{cases}
\end{align}
The bundle $\sphere^1\times \sphere^2 \rightarrow \hat{\overline{\Fiber_\sigma}} \rightarrow \widetilde{S}$ is then the corresponding degree two covering of $\overline{\Fiber_\sigma}$.
As $\widetilde{S}$ is contractible, the bundle $\hat{\overline{\Fiber_\sigma}} \rightarrow \widetilde{S}$ is topologically trivial and hence $\pi_1(\hat{\overline{\Fiber_\sigma}})\cong \Z$.  

For later, we note a description of the complex line bundles whose sphere bundles appear above. 

\begin{remark}\label{Remk:SpecialLineBundle}
There are complex line bundle isomorphisms as follows. 
When $\sigma =\beta$, we have $B \cong \K^3$ and 
when $\sigma =\alpha$, we have $\Hom_{\C}(N,B) \cong \mathcal{T}^3\K^3$. See \cite[Section 4.2, Section 5.6]{DE26b}. 
\end{remark}

We now show that a loop in the developed fiber is simply linked with the lifted limit map. The key idea is to reduce to the model cases considered above by appropriate homotopies. 

\begin{lemma}[Simple linking: fiber and limit map]
\label{lem:linkingfiberAndLimitMap}
Let $(\V,\Phi)$ be a polystable $\sigma$-bundle, which 
is also a Hodge bundle with $\mathrm{div}(\delta)\neq0 $ when $\sigma =\beta$, and $\hat{\xi}^{\sigma}: \mathbb{S}^1 \rightarrow \widehat{\mathcal{F}}_{\sigma}$ be the lifted limit map. 
Fix a point $p \in \tilde{S}$ and let $\gamma:\sphere^1 \rightarrow \widehat{\mathcal{F}}_{\sigma}$ be the inclusion of a generator of $\pi_1(\hat{\overline{{\Fiber}_\sigma}}_{\mid p})$. 
Then $\ln(\gamma, \hat{\xi}^{\sigma}) = \pm 1$. 
\end{lemma}

\begin{proof}
Let $\gamma: \mathbb{S}^1 \rightarrow \hat{\overline{{\Fiber}_\sigma}}_{\mid p}$ be a loop representing a generator of $\pi_1(\hat{\overline{{\Fiber}_\sigma}}_{\mid p})$. By \eqref{eq:DoubleCover}, $\hat{\overline{\Fiber_\sigma}}$ can be viewed as a circle bundle over $\sphere(U)$, where $U = \mathscr{L}\oplus N$. We start by homotoping $\gamma$ inside $\hat{\overline{{\Fiber}_\sigma}}_{\mid p}$ to the model generator $\gamma_B:\mathbb{S}^1 \rightarrow\hat{\overline{{\Fiber}_\sigma}}_{\mid p}$ given by $\widehat{\Ein}(\hat{\mathscr{L}}_{p} \oplus B_p)$ if $\sigma=\beta$ (see Remark \ref{rem:CircleEinstein}) and any lift of 
$\Pho^\times(N_{p} \oplus B_p)$ if $\sigma=\alpha$ (see Lemma \ref{Lem:OrthogonalPhoton}). 
We shall write again $\tau:\mathbb{S}^1\times \mathbb{S}^2 \rightarrow \widehat{\mathcal{F}_{\sigma}}$ for the thickening of the lifted limit set $\hat{\Lambda} = \im(\hat{\xi}^{\sigma})$. By Theorem \ref{thm:BetaAvoidsLimitSet} and \ref{thm:AlphaAvoidsLimitSet} this homotopy avoids the thickening of $\hat{\Lambda}$, and hence does not change the linking number. With $p$ still fixed, we write $(e, T,N,B)$ for an oriented Frenet splitting refining $(\mathscr{L}_p, T_p,N_p, B_p)$. 
We then imagine $e$ as orienting each of $\mathscr{L},T,N,B$, and we endow  the corresponding oriented subspaces with a hat for distinction. 
 \medskip

Since $\dev: \hat{\overline{{\Fiber}_\sigma}} \rightarrow \widehat{\mathcal{F}}_{\sigma}$ is injective on fibers, we can identify $\gamma$ with its image under $\dev$.
Now, consider the map $\tau_{\infty}: \sphere^1\times \sphere^2 \rightarrow \widehat{\mathcal{F}}_{\sigma}$, which will serve as an approximation of $\tau$, given as follows:
\[ \begin{cases} 
\tau_\infty: K_{\widehat{\Ein}(\hat{\mathscr{L}}\oplus T)} \hookrightarrow \widehat{\mathcal{F}}_{\beta} &  \sigma=\beta ,\\
\tau_\infty: K_{\widehat{\Pho}(\hat{N}\oplus T)} \hookrightarrow \widehat{\mathcal{F}}_{\alpha} & \sigma =\alpha.
\end{cases}
\]
By Lemmas \ref{Lem:AlphaFibers+LimitSet} and \ref{Lem:BetaFibers+LimitSet}, we have $\ln(\gamma_B, \tau_{\infty}) =\pm 1$. 

We now construct a homotopy $(\tau_t)_{t\in [0,\infty]}$ between $\tau=\tau_0$ and $\tau_{\infty}$ that avoids $\gamma_B$. We begin with the case $\sigma=\beta$. It suffices to build a homotopy $(\xi_t)_{t\in [0,\infty]}$ between the curves $\xi_0=\xi^\beta$ and $\xi_\infty$, a parametrization of $\widehat{\Ein}(\hat{\mathscr{L}}\oplus T)$, such that the thickenings of the intermediate curves avoid $\gamma$. Given a unit tangent vector $v\in \T_p\widetilde{\Sigma}$, define $x(v)\in \partial\pi_1S$ the endpoint of the $g_{\beta}$-geodesic with initial velocity $v$. The cone $\hat{C}'(v)$ contains $\xi^\beta(x)$ by Corollary \ref{Cor:BetaLimitSetTrapped}. This cone $\hat{C}'(v)$ is preserved by the one-parameter subgroup $\big(\exp(t\Psi_0(v))\big)_{t \in \R}$. Indeed recall that $\hat{C}'(v)=\lbrace [x]\in \RP^6\mid B_v'(x)> 0\rbrace $ where $B_v' = h|_{\mathbb{E}_{2}(v)}-h|_{\mathbb{E}_{0}(v)}-h|_{\mathbb{E}_{-2}(v)}-2h|_{\mathbb{E}_{1}(v)}-2h|_{\mathbb{E}_{-1}(v)}$, $h$ is the harmonic metric, and $\mathbb{E}_{k}(v)$ is the $k$-eigenspace of $\Psi_0$. Now, one directly checks that 
$B_v'(X)>0$ implies $B_v'(\exp(t\Psi_0(v))\cdot X)>0$ for any $t \geq 0$. Indeed, the positive term in $B_v'(X)$ grows proportionally more when applying $\exp(t\Psi_0)$.

Moreover, for every $\ell\in C(v)$, $(\exp(t\Psi_0(v))\cdot \ell)_{t\geq 0}$ converges to $\ell_0\in \widehat{\Ein}(\mathscr{L}\oplus T)$. For $t\geq 0$, we define $\xi_t: \partial \pi_1S\to \mathcal{F}_\sigma$ so that:
\[ \xi_t(x(v))=\exp(t\Psi_0(v))\cdot \xi_0(x(v)).\] 
The dynamics of $\Psi_0(v)$ described above implies that these curves $\xi_t$ converge pointwise to $\xi_\infty$. Hence it remains only to check that the thickenings of these intermediate curves $\hat{\Lambda}_t:= \im(\xi_t)$ are disjoint from $\gamma$. This is a consequence of Theorem \ref{thm:BetaAvoidsLimitSet}, as the thickening of any element of any cone $\hat{C}'(v)$ for $v\in \T\widetilde{\Sigma}$ is disjoint with the developed fiber $\dev( \hat{\overline{{\Fiber}_\sigma}}|_{p})$, and hence in particular is disjoint with $\gamma$.

The proof for $\sigma=\alpha$ is completely analogous, up to replacing the isotropic lines in $\hat{C}'(v)$ by the photons included in the corresponding cones $C(v)$, and using Corollary \ref{Cor:AlphaLimitSetTrapped} and Theorem \ref{thm:AlphaAvoidsLimitSet}.
\medskip

These two homotopies show that $\ln(\gamma, \tau) =\ln(\gamma_B, \tau) =\ln(\gamma_B, \tau_\infty)=\pm 1$ for $\sigma \in \{\alpha,\beta\}$.
\end{proof}

We can now relate the fundamental group of the domain of discontinuity $\hat{\Omega}_{\sigma}$ in $\widehat{\mathcal{F}}_{\sigma}$ to the fundamental group of the moving base of pencil manifold $\hat{\overline{\Fiber}_{\sigma}}$. 

\begin{theorem}\label{thm:DomainPi1}
Let ${\dev}: \hat{\overline{{\Fiber}}_{\sigma}}\rightarrow \widehat{\mathcal{F}}_{\sigma}$ be the developing map associated to a polystable $\sigma$-Hodge bundle $(\V,\Phi)$, for $\sigma \in \{\alpha,\beta\}$, that also satisfies $\mathrm{div}(\delta)\neq0$ when $\sigma =\beta$. Then there is a short exact sequence 
\[ 1\longrightarrow N\longrightarrow \pi_1(\hat{\Omega}^\rho_{\sigma})\longrightarrow\pi_1(\hat{\overline{{\Fiber}_{\sigma}}})\longrightarrow 1, \]
where $N$ is a finite group. 
Moreover, the induced homomorphism $\dev_*:\pi_1(\hat{\overline{{\Fiber}_{\sigma}}})\to \pi_1(\hat{\Omega}^\rho_{\sigma})$ defines a section of the short-exact sequence, which induces an isomorphism $\pi_1(\hat{\Omega}_\rho^\sigma) \cong N\rtimes \pi_1(\hat{\overline{{\Fiber}_{\sigma}}})$.
\end{theorem}

\begin{proof}
Let $\hat{\xi}^{\sigma}: \sphere^1 \rightarrow \widehat{\mathcal{F}}_{\sigma}$ be the lifted limit map. 
We use the linking number, as in Definition \ref{Defn:LinkingG2Flags}, to define the following homomorphism: 
\[ q: \pi_1(\hat{\Omega}^\rho_{\sigma})\to \Z \cong \pi_1(\hat{\overline{{\Fiber}_{\sigma}}}), \qquad  
q([\gamma ])= \ln(\gamma, \hat{\xi}^{\sigma}) .\]
Note that $q$ is surjective by Lemma \ref{lem:linkingfiberAndLimitMap}, as the value $\pm 1\in \Z$ is attained by the loop $[\gamma_B]$.

Up to changing the orientation of the loop $\gamma_B$, which generates $\pi_1(\hat{\overline{{\Fiber}_{\sigma}}})$, one can assume that $q([\gamma_B])=+1$. Thus, the induced homomorphism $\dev_*:\pi_1(\hat{\overline{{\Fiber}_{\sigma}}})\to \pi_1(\hat{\Omega}^{\rho}_{\sigma})$ is a right inverse of $q$. 

Since $\dev: \hat{\overline{{\Fiber}_{\sigma}}} \rightarrow \hat{\Omega}^{\rho}_{\sigma}$ is a finite covering, we know that $\pi_1(\hat{\Omega}^{\rho}_{\sigma})$ is a finite extension of $\pi_1(\hat{\overline{{\Fiber}_{\sigma}}})\cong \Z$. Thus, the subgroup $N :=\ker(q)$ is finite. The desired result follows. 
\end{proof}

\begin{corollary}\label{Cor:FundamentalGroupSES}
There is a short exact sequence: 
\[ 1 \longrightarrow N \longrightarrow \pi_1(\hat{M}^\rho_{\sigma}) \longrightarrow \pi_1(\hat{{\Fiber}}_{\sigma}) \longrightarrow 1\ .\]
\end{corollary}

\begin{proof}  
Consider the following commutative diagram of short-exact sequences:
\[
\begin{tikzcd}
&&1\arrow[d]&1 \arrow[d]& \\
1 \arrow[r] &  N \arrow[r]  & \pi_1(\hat{\Omega}_\sigma^\rho) \arrow[r,rightarrow] \arrow[d] & \pi_1( \hat{\overline{{\Fiber}_{\sigma}}}) \arrow[r] \arrow[d]&1 \\
& & \pi_1(\hat{M}^\rho_{\sigma})  \arrow[d] & \pi_1(\hat{{\Fiber}}_{\sigma})  \arrow[d] \\
&&\pi_1S \arrow[d]\arrow[r, rightarrow]&\pi_1S\arrow[d]& \\
&&1&1& \\
\end{tikzcd}
\]
First note that $N \triangleleft  \pi_1(\hat{M}_\sigma^\rho)$. 
Indeed, $N$ is the normal subgroup of all torsion elements in $\pi_1(\hat{\Omega}_\sigma^\rho)$, and is hence preserved by every automorphism of $\pi_1(\hat{\Omega}_\sigma^\rho)$, and we have $\pi_1(\hat{\Omega}^\rho_\sigma)\triangleleft\pi_1(\hat{M}^\rho_\sigma)$.

We then quotient the left vertical short-exact sequence by $N$, and use that we have a map $ \pi_1(\hat{{\Fiber}}_\sigma)\rightarrow \pi_1(\hat{M}_\sigma^{\rho})$. Indeed, we have a $\rho$-equivariant map $\dev: \hat{\overline{{\Fiber}_\sigma}} \rightarrow \hat{\Omega}_\sigma^{\rho}$ that descends to a map $\dev_0: {\hat{\Fiber}}_\sigma \rightarrow \hat{M}_\sigma^{\rho}$. 
We set $\theta = {\dev}_*$, which descends to $\theta_0: \pi_1(\hat{\mathscr{M}}_{\sigma})\rightarrow \pi_1(\hat{M}_{\sigma}^{\rho})/N$. We obtain:  
\[
\begin{tikzcd}
1\arrow[d]&1 \arrow[d]& \\
  \pi_1(\hat{\Omega}_\sigma^\rho)/N \arrow[r,leftarrow,"\cong", "\theta" below] \arrow[d] & \pi_1(\hat{\overline{{\Fiber}_{\sigma}}}) \arrow[d] \\
 \pi_1(\hat{M}_\sigma^\rho)/N  \arrow[d] & \pi_1(\hat{{\Fiber}}_\sigma) \arrow[d]\arrow[l,"  \theta_0" below] \\
\pi_1S \arrow[d]\arrow[r,leftrightarrow,"\cong"]&\pi_1S\arrow[d]& \\
1&1& \\
\end{tikzcd}
\]
This diagram is commutative, as $\dev$ is $\rho$-equivariant. By the 5-lemma in the category of groups, we conclude the middle map is an isomorphism. Tracing back to the original diagram before the quotient by $N$, we obtain our desired short-exact sequence, where 
\begin{itemize}
    \item The map $N \rightarrow \pi_1(\hat{M}_\sigma^{\rho})$ is the composition of inclusions $N \rightarrow \pi_1(\hat{\Omega}_\sigma^{\rho}) \rightarrow \pi_1(\hat{M}_\sigma^{\rho})$.
    \item The map $\pi_1(\hat{M}_\sigma^{\rho}) \rightarrow \pi_1(\hat{{\Fiber}}_\sigma)$ induces the isomorphism $\pi_1(\hat{M}_\sigma^{\rho})/N \cong \pi_1(\hat{{\Fiber}}_\sigma)$. $\qedhere$
\end{itemize}
\end{proof}

We believe that $N$ is in fact trivial. However, Corollary \ref{Cor:FundamentalGroupSES} is still strong enough for our main application of interest: distinguishing connected components of Anosov representations.

\subsection{Connected components of Anosov representations}
We now show that the fundamental group of the manifold $\hat{\mathscr{M}}_{\sigma}^{\rho}$ is an invariant for the connected component of $\rho$ in the space of $P_{\sigma}$-Anosov representations in $\Gtwosplit$. We describe this fundamental group explicitly, which is a torsion-free extension of $\pi_1S$, and then show that many components of $P_{\sigma}$-Anosov representations are detected. 

Let us begin with an elementary result about the isomorphism classes of extensions of torsion-free groups by finite groups.

\begin{proposition}\label{Prop:GroupTheoryExtension}
    Let $H_1,H_2$ be torsion-free groups and let $N_1,N_2$ be finite groups. Consider groups $G_1$, $G_2$ that admit exact sequences:
$$1\to N_1\to G_1\to H_1\to 1,$$
$$1\to N_2\to G_2\to H_2\to 1.$$
If $G_1$ and $G_2$ are isomorphic, then $H_1$ and $H_2$ are isomorphic.
\end{proposition}

\begin{proof} The subgroup $N_i$ is naturally identified with a normal subgroup of $G_i$ that consists of exactly all torsion elements by the hypotheses. 
Hence, any isomorphism $ G_1 \cong G_2$ induces an isomorphism $H_1\cong G_1/N_1 \cong G_2/N_2\cong H_2$. 
\end{proof}

The following result describes the fundamental group of the  $(\Gtwosplit, \widehat{\mathcal{F}}_{\sigma})$-manifold $\hat{\Fiber}_{\sigma}$. 
\begin{lemma}\label{Lem:pi1Presentation}
Let $(\V,\Phi)$ be a polystable $\sigma$-bundle on a Riemann surface $\Sigma_g=(S_g, J)$ and let
$\hat{\Fiber}_{\sigma}$ be the manifold in \eqref{eq:DoubleCover}. Then 
\begin{enumerate}[label=(\roman*)]
    \item The fundamental group of $\hat{\Fiber}_{\sigma}$ is equal to the following central extension of $\pi_1(S_g)$:
    \[ \pi_1(\hat{{\Fiber}}_\sigma)\cong \Gamma_g^d =\left\langle t,a_1,b_1,\cdots, a_g,b_g\; \bigg | \; t^{d}\prod_{i=1}^g[a_i,b_i], \; [t, a_i], \;[t,b_i]\right\rangle, \]
where $d=\deg(\mathcal{L}_{\sigma})\geq 0$ and 
\[ \mathcal{L}_{\sigma} = \begin{cases} \mathcal{B} & \sigma = \beta, \\
\mathcal{T}^{3}\mathcal{K}^3 & \sigma = \alpha.
\end{cases} \] 
    \item The group $\Gamma_g^d$ is the central extension $1\to \Z\to \Gamma^d_g\to \Gamma_g\to 1$ associated to the fiber bundle $\sphere^1 \rightarrow \sphere(\mathcal{L}_{\sigma}) \rightarrow S_g$.
    \item The first integer homology group of $\hat{\Fiber}_{\sigma}$ is given by $H_1(\hat{\Fiber}_{\sigma}, \Z) \cong \Z/d\Z\oplus \Z^{2g}$. 
\end{enumerate}
\end{lemma}

Note that the condition $d\geq 0$ in (i) follows from the degree constraints in Lemma \ref{Lem:ConnectedComponents}. 
\begin{proof}
Using Remark \ref{Remk:SpecialLineBundle}, there is an inclusion of fiber bundles  $\sphere(\mathcal{L}_{\sigma})\hookrightarrow \hat{\mathscr{M}}_{\sigma}$, which induces an isomorphism of fundamental groups. Thus, it suffices to prove the given description applies to the group $\Gamma^d_{g} : = \pi_1(\sphere(\mathcal{L}_{\sigma}))$. Since the surface $S_g$ is aspherical, observe that the long exact sequence of homotopy groups of the fibration $\sphere(\mathcal{L}_{\sigma})\rightarrow S_g$ provides a short-exact sequence $1 \rightarrow \Z \rightarrow \Gamma^{d}_{g} \rightarrow \Gamma_g\rightarrow 1$. 

Let us consider the restriction of the bundle $\sphere(\mathcal{L}_{\sigma}) \rightarrow S_g$ to $S_g\setminus D$, for a small disk $D\subset S_g$. The standard characterization of the degree of $\mathcal{L}_\sigma$ is that if one trivializes the corresponding circle bundle as $(S_g\setminus D)\times \mathbb{S}^1$, then the total space of $\sphere(\mathcal{L}_{\sigma})$ is obtained by gluing in $D\times \mathbb{S}^1$ by the map 
\[ \phi: \partial D\times \mathbb{S}^1\to \partial(S_g\setminus D)\times \mathbb{S}^1, \qquad (x,s)\mapsto (x, s+dx),\]
where $\sphere^1 \simeq \R/\Z$. Hence, the fundamental group of $\sphere(\mathcal{L}_{\sigma})$ realized after the gluing is the quotient of $\Z\times F_{2g}$, where $F_{2g}$ the free group generated by the $2g$ generators of $\pi_1(S_g \setminus D)$, with the class of the curve $\gamma:\mathbb{S}^1 \to \partial (S_g\setminus D)\times \mathbb{S}^1$, given by $x\mapsto (x, s_0+xd)$, for some point $s_0\in \mathbb{S}^1$. This curve $\gamma$ is represented by $t^d\Pi_{i=1}^g[a_i,b_i]$ where $t$ is the generator of the $\Z$-factor and $(a_i,b_i)_{i=1}^g$ are the standard generators of $\pi_1(S_g\setminus D)$.

The homology group $H_1(\hat{\Fiber}_{\sigma},\Z)$ is the abelianization of $\Gamma_g^d$, obtained by imposing the relation $t^d=1$ on the free abelian group with generators $(t,a_1,b_1,\cdots, a_g,b_g)$. This group is isomorphic to $\Z/d\Z \oplus \Z^{2g}$. 
\end{proof}

Now, we combine the two lemmas. 
\begin{corollary}[A geometric invariant]\label{Cor:pi1Invariant}
Fix representations $\rho_i\in \chi^{Anosov}_{\sigma}(S_g,\Gtwosplit)$ associated to polystable $\sigma$-bundles $(\V, \Phi_i)$ over  Riemann surfaces $\Sigma_i=(S_g,J_i)$, with corresponding $(\Gtwosplit, \widehat{\mathcal{F}}_{\sigma})$-geometric structures $\hat{\Fiber}^{i}_\sigma$. If $\rho_1, \rho_2$ are in the same connected component of $P_{\sigma}$-Anosov representations, then $\pi_1(\hat{{\Fiber}}^1_\sigma) \cong \pi_1(\hat{{\Fiber}}^2_\sigma)$.  
\end{corollary}

\begin{proof}
For $\rho: \pi_1S \rightarrow \Gtwosplit$ that is $P_{\sigma}$-Anosov, denote $\hat{M}_{\sigma}^{\rho}:= \rho(\pi_1S)\backslash \hat{\Omega}_{\sigma}^{\rho}$, where 
$\hat{\Omega}_{\sigma}^{\rho} \subset \widehat{\mathcal{F}}_{\sigma}$ is the lift of the domain of discontinuity $\Omega_{\sigma}^{\rho} \subset \mathcal{F}_{\sigma}$ to the double cover.

Set $G_i = \pi_1( \hat{M}_\sigma^{\rho_i})$ and $H_i = \pi_1(\hat{{\Fiber}}^i_\sigma)$. By Lemma \ref{Lem:pi1Presentation}, $H_i$ is torsion-free. Also $G_1 \cong G_2$ as in \cite[Theorem 2.10]{AMTW25}. 
Hence, we apply Proposition \ref{Prop:GroupTheoryExtension} and conclude the corollary. 
\end{proof}

We can explain more intrinsically the geometric invariant that has been detected.  

\begin{remark}[An intrinsic geometric invariant]
Let $\rho:\pi_1S\rightarrow \Gtwosplit$ be associated to a polystable $\sigma$-bundle, which is Hodge and has $\mathrm{div}(\delta)\neq 0$ when $\sigma= \beta$. Set $\Lambda =\im(\hat{\xi}^{\sigma})$ as the limit set in $\widehat{\mathcal{F}}_{\sigma}$. 

Theorem \ref{thm:DomainPi1} and Lemma \ref{Lem:pi1Presentation} imply that $\pi_1(\hat{\Omega}_{\rho}^{\sigma}) \cong N \rtimes \Z$ and moreover that any primitive infinite order element $\gamma \in \pi_1(\hat{\Omega}_{\rho}^{\sigma})$ satisfies $\ln(\gamma, \Lambda)=\pm 1$ and $\gamma$ descends to an order $d$ element $[\gamma] \in H_1(\hat{M}_{\rho}^{\sigma}, \Z)$, unless $d=0$ in which case $[\gamma]$ is infinite order. 
Here, $d$ is given by Lemma \ref{Lem:pi1Presentation}. 
\end{remark}

We can then certify distinct connected components of Anosov representations by analyzing the possible degrees of $\mathcal{L}_{\sigma}$. To this end, we present the following lemma, which just summarizes results from \cite{CT24} and \cite{DE26b} on the cases $\sigma = \beta$ and $\sigma  =\alpha$, respectively. 

\begin{lemma}[Components of $\sigma$-bundles]\label{Lem:ConnectedComponents}
Let $\mathcal{M}_{\sigma}(\Sigma)$ be the moduli space of polystable $\sigma$-bundles on $\Sigma = (S_g,J)$. 
Let $d$ be the integer $d=\deg(\mathcal{T})$ when $\sigma = \alpha$ and $d=\deg(\mathcal{B})$ when $\sigma=\beta$. 
Then 
\begin{enumerate}[label=(\roman*), itemsep = 1ex]
    \item $\mathcal{M}_{\alpha}(\Sigma)$ has $4g-3$ connected components $\mathcal{M}_{\alpha}^d(\Sigma)$, labeled by $-2g+2\leq d \leq 2g-2$. 
    \item $\mathcal{M}_{\beta}(\Sigma)$ has $6g-5$ connected components $\mathcal{M}^d_{\beta}(\Sigma)$, labeled by $0\leq d \leq 6g-6$.
\end{enumerate}
\end{lemma} 
\begin{proof}
(i) The result follows from \cite[Proposition 5.20]{DE26b}. 

(ii) The result follows mostly from \cite[Proposition 5.6]{CT24}, and the mild extension in \cite[Proposition 4.7]{DE26b}.   
\end{proof}

Finally, we obtain the result of interest: a verification of many connected components of $P_{\sigma}$-Anosov representations using $\sigma$-bundles. For the statement, we also denote by $\mathcal{M}^{*}_\beta(\Sigma)$ the sub-moduli space of Hodge bundles in $\mathcal{M}_\beta^d(\Sigma)$ that also have degree $d >0$. Recall that a $\beta$-bundle $(\V,\Phi)$, as in \eqref{BetaCyclic}, is a Hodge bundle exactly when $\alpha=0$ or $\delta=0$. 

\begin{theorem}[Correspondence of strata]\label{thm:StrataCorrespond}
Let $\mathcal{M}_{\sigma}(\Sigma)$ be the moduli space of polystable $\sigma$-bundles on $\Sigma=(S_g,J)$. Then 
\begin{enumerate}[label = (\roman*), itemsep = 1ex]
    \item $\NAH_{\Sigma}$ induces an injective map between the set of connected components of $\mathcal{M}_{\alpha}(\Sigma)$ and the set of connected components of $\chi_{\alpha}^{Anosov}(S_g,\Gtwosplit)$. 
    \item Every connected component of $\NAH_{\Sigma}(\mathcal{M}_{\beta}(\Sigma) \setminus \mathcal{M}_{\beta}^0(\Sigma))$ contains $P_{\beta}$-Anosov representations. Moreover, the map $\NAH_{\Sigma}:\mathcal{M}_\beta^*(\Sigma)\rightarrow\chi_{\beta}^{Anosov}(S_g,\Gtwosplit)$
    is injective on the level of connected components.
\end{enumerate}
\end{theorem}

\begin{proof}
(i) By Theorem \ref{thm:AlphaAnosov}, the map $\NAH_{\Sigma}:\mathcal{M}_{\alpha}(\Sigma)\rightarrow \chi(S_g,\Gtwosplit)$ has image in $\chi_{\alpha}^{Anosov}(S,\Gtwosplit)$. 
By Lemma \ref{Lem:ConnectedComponents}, the map $\mathcal{M}_{\alpha}(\Sigma) \rightarrow \Z$ by 
$[\mathcal{T}, \beta, \delta] \rightarrow \deg(\mathcal{T}^3\mathcal{K}^3)$ is injective and has image in $\Z_{\geq 0}$. The result then follows by Lemma \ref{Lem:pi1Presentation} and Corollary \ref{Cor:pi1Invariant}. \medskip  

(ii) Similarly, by Theorem \ref{Thm:BetaAnosov}, we have a well-defined map $\NAH_{\Sigma}:\mathcal{M}_{\beta}^{*}(\Sigma) \rightarrow \chi_{\beta}^{Anosov}(S_g,\Gtwosplit)$. Indeed, polystable $\beta$-bundles for which $\deg(\mathcal{B})>0$ have differential $\delta$ with at least one zero. 
By Corollary \ref{Cor:pi1Invariant} and Lemma \ref{Lem:ConnectedComponents}, the injectivity of this map on the level of connected components follows. 

Now, by \cite[Proposition 5.8]{DE26b}, every connected component of $\mathcal{M}_{\beta}(\Sigma) \setminus \mathcal{M}_{\beta}^0(\Sigma)$ contains a Hodge bundle with either $\alpha =0$ or $\delta =0$. Hence, every connected component of the image of $\mathcal{M}_{\beta}(\Sigma) \setminus \mathcal{M}_{\beta}^0(\Sigma)$ under $\NAH_{\Sigma}$ contains $P_{\beta}$-Anosov representations. \end{proof}

\begin{corollary}[Connected components of Anosov representations]
The space $\chi_{\alpha}^{Anosov}(S_g,\Gtwosplit)$ has at least $4g-3$ connected components and $\chi_{\beta}^{Anosov}(S_g,\Gtwosplit)$ has at least $6g-6$ connected components. 
\end{corollary}

\appendix

\section{Connected Components of \texorpdfstring{$\chi(S,\Gtwosplit)$}{chi(S,G2')}}
\label{Apendix:CC}

 In this appendix we explain how to use the results from \cite{AABCGGPO19} to count the connected components of the character variety $\chi(S,\Gtwosplit)$ for a closed surface of genus $g\geq 2$. In particular, $\chi(S,\Gtwosplit)$ has exactly three connected components. This fact, which was apparently already known by the community but never written, was communicated to us by Brian Collier. 

To label the components, we recall that the maximal compact subgroup of $\Gtwosplit$ is $K \cong \SO(4,\R)$. 
Now, any principal $\Gtwosplit$-bundle reduces structure uniquely up to isomorphism to $K$. Hence, the flat rank seven vector bundle $E_{\rho}$ associated to a representation $\rho:\pi_1S \rightarrow \Gtwosplit$ decomposes as $E=E^3_{\rho} \oplus E^{4}_{\rho}$. Each of the bundles $E^3_{\rho}, E^{4}_{\rho}$ are orientable. 
Thus, the unique topological invariant of $\rho$ is the common second Stiefel-Whitney class $w_2(E^3_{\rho})=w_2(E^4_{\rho})$, which shall aid in the following classification. 

\begin{theorem}[Components of the $\Gtwosplit$-character variety]
 \label{thm:AppendixCountConnectedComponents}
 The character variety $\chi(S,\Gtwosplit)$ has exactly three connected components: 
 \begin{itemize}
     \item the Hitchin component $\mathcal{C}_{\text{Hit}}$, 
     \item a component $\mathcal{C}_0$ containing representations such that $w_2(E^4_{\rho})=0$, and containing representations factoring through $K$,
     \item a component $\mathcal{C}_1$ containing representations such that $w_2(E^4_{\rho})\neq 0$, and  containing representations factoring through $K$.
 \end{itemize} 
 \end{theorem}

We complete the proof of this theorem at the end of this appendix, as a simple consequence of the technical Lemma \ref{Lem:G2Minima}.

\begin{remark}
    Theorem \ref{thm:AppendixCountConnectedComponents} goes in the direction of the general conjecture (see \cite{BCGBO24}) that for a simple real Lie group $G$ with semisimple maximal compact subgroup $K$, the connected components of the $G$-character variety admit the following count: 
    \[ |\pi_0(\chi(S,G))| = |\{\text{topological invariants in }K\}|\; + \; |\{G-\text{higher Teichmüller spaces}\}|.\]
    For $G=\Gtwosplit$, we get $3=2+1$ from $2=|\pi_1(K)|$ via $K \cong \SO(4)$, and $1$ from the Hitchin component.
\end{remark}

Let us briefly resume the setting of \cite{AABCGGPO19}. If we fix a Riemann surface structure $\Sigma$ on the surface $S$, the non-abelian Hodge correspondence provides a homeomorphism between the character variety $\chi(S,\Gtwosplit)$ and the moduli space $\mathcal{M}_{\Sigma}(\Gtwosplit)$ of $\Gtwosplit$-polystable Higgs bundles $(\mathcal{E},\Phi)$ on $\Sigma$. Given a representative Higgs bundle $(\V,\Phi)$, one can associate the following number:
\[ E(\mathcal{E},\Phi)=\int_\Sigma \lVert \Phi\rVert_h^2. \]
Here, $h$ is any harmonic metric solving Hitchin's equations. Since $E$ is independent of the gauge representative, $E$ defines a function $E: \mathcal{M}_{\Sigma}(\Gtwosplit) \rightarrow \R$, first introduced in Hitchin's seminal work \cite{Hit87}, and used subsequently to count connected components of surface group character varieties in many other contexts. We refer to the survey \cite{Bra24} for further details and history.\medskip

Hitchin proved the function $E$ is proper and its critical points, on the smooth locus, are exactly the fixed points of the $\C^*$-action on the space of Higgs bundles, given for $\lambda\in \C^*$, by $\lambda\cdot(\mathcal{E},\Phi)=(\mathcal{E},\lambda\Phi)$.

\begin{remark}
A helpful fact here is that $\Gtwosplit$-stability is weaker than $\SO(3,4)$-stability, meaning more Higgs bundles are $\Gtwosplit$-stable than $\SO(3,4)$-stable. In particular, there are few $\Gtwosplit$-strictly polystable bundles to consider as potential minima of $E$.
\end{remark}

On the other hand, polystability in complex Lie groups is more rigid. 
\begin{remark}
By \cite{GPGMiR12}, $\Gtwo^\C$-polystability and $\SL(7,\C)$-polystability agree. 
\end{remark}

The key lemma for the count of connected component is the following classification for local minima of $E$. In other words, there is only one non-trivial minimum of the function $E$. 
\begin{lemma}[Classifying local minima]\label{Lem:G2Minima}
Let $(\mathcal{E},\Phi)$ be a $\Gtwosplit$-polystable Higgs bundle on $\Sigma$ that is a local minimum of $E$. Then either it is the unique $\C^*$-fixed point in the Hitchin component or $\Phi=0$.   
\end{lemma}

\begin{proof}
Suppose that $(\V,\Phi)$ is a local minimum of $E$. We shall distinguish the cases that $(\mathcal{E},\Phi)$ is $\Gtwosplit$-stable (see \cite[{Definition 3.2}]{GPGMiR12}) and the case when it is \emph{strictly $\Gtwosplit$-polystable}, i.e., polystable but not stable. \medskip 

\textbf{Case 1: $(\mathcal{E},\Phi)$ is $\Gtwosplit$-stable.} 
Let us recall facts from \cite[Appendix A]{AABCGGPO19}. Since $(\mathcal{E},\Phi)$ is a critical point of $E$, it is a $\C^*$-fixed point. Thus, for each $\lambda\in \C^*$ there is a gauge transformation $g_\lambda$ of $\mathcal{E}$ such that $\Ad_{g_\lambda}(\Phi)=e^{i\lambda}\Phi$. By \cite[Proposition A.12]{AABCGGPO19}, we can avoid questions of uniqueness of $g_{\lambda}$ and conclude that the following infinitesimal gauge transformation is always well-defined: 
\[\psi=\frac{\mathrm{d}}{\mathrm{d}\lambda}\mid_{\lambda=0}g_\lambda. \]

Let $\g_2' = \frakk \oplus \frakp$ be a Cartan decomposition (written as $\mathfrak{h}+\mathfrak{m}$ in \cite{AABCGGPO19}). Taking the complexified splitting, we obtain $\g_2^{\C} = \frakk^{\C} \oplus \frakp^{\C}$. The infinitesimal gauge transformation 
$\psi$ is a section $\psi \in \Omega^0(\Sigma, \mathcal{E}^{\R}[\frakk])$ that satisfies $\nabla^h \psi=0$, where $\nabla^h$ is the Chern connection associated to the harmonic metric on $\mathcal{E}[\mathfrak{g}_2']$.

As explained in \cite{AABCGGPO19}, for any point $p \in \Sigma$, the element $\psi_p \in \End(\V_p)$ has the same eigenvalues. Since $\psi$ lies in a compact subalgebra, its eigenvalues are all imaginary. 
For the rest of the proof we write $\g = \g_2'$ and $\g^{\C} = \g_2^{\C}$ so that the subscripts `2' do not cause confusion with the additional decorations that we now introduce. By studying the eigenspaces of $\ad_{\psi}$, we obtain a corresponding weight space decomposition of the Lie algebra bundle $\V[\g^\C]$ as follows. We first split according to the complexified Cartan decomposition: 
\begin{equation}
    \V[\g^\C] = \V[\frakk^\C] \oplus \V[\frakp^{\C}],
\end{equation}
then further decompose into weight spaces
\begin{equation}\label{eq:Grading}
    \V[\g^\C]_{\ell} = \V[\frakk^\C]_{\ell} \oplus \V[\frakp^{\C}]_{\ell},
\end{equation}
where the subscript $\ell$ denotes the $i\ell$ eigenspace of $\psi$.  \medskip 

Observe that by definition of $\psi$, we have $[\psi, \Phi] = i \Phi$. 
Thus, $\Phi \in \Omega^0(\Sigma, \V[\frakp^{\C}]_1)$. Hence, the adjoint of $\Phi$, once restricted to $\V[\frakk^{\C}]$, becomes a map $\ad_{\Phi}: \V[\frakk^{\C}] \rightarrow \V[\frakp^{\C}]\otimes \K$, which then decomposes into a sum of holomorphic vector bundle maps of the following type:
\begin{equation}\label{eq:AppendixHypercohomology}
    \phi_\ell: \mathcal{E}[\mathfrak{k}^\C]_{\ell}\to \mathcal{E}[\mathfrak{p}^\C]_{\ell+1}\otimes \mathcal{K}.
\end{equation}

We now recall the Lie-theoretic characterization of smooth local minima of $E$.
\begin{proposition}[{\cite[{Proposition 3.4 and Proposition A.2}]{AABCGGPO19}}]
\label{prop:AppendixPConditionSmooth}
    A $\Gtwosplit$-Higgs bundle $(\mathcal{E},\Phi)$-that is $\Gtwo^{\C}$-stable is a local minimum of $E$ if and only if for every $\ell>0$, the map $\phi_\ell$ from \eqref{eq:AppendixHypercohomology} is an isomorphism.
\end{proposition}
By applying Proposition \ref{prop:AppendixPConditionSmooth}, we will see that Case 1 becomes nothing more than a (tedious) exercise in $\g_2'$ Lie theory. \medskip

As $i\psi$ is a regular semisimple element with real eigenvalues, we can find a Cartan subalgebra $\mathfrak{a}^{\C} \subset \mathfrak{p}^{\C}$ containing $i\psi$. 
We distinguish $3$ different cases depending on the regularity of $\tau = \boldsymbol{\mu}(i\psi)$, the Cartan projection of $i\psi$. Note that $\tau=0$ is impossible or else $\Phi =0$, which contradicts stability. 
\medskip 

We can choose simple roots $\alpha,\beta$ such that $\alpha(\tau),\beta(\tau ) \geq0 $, with at least one of $\alpha(\tau), \beta(\tau)$ being positive. 
Note that the positive roots $\{\alpha,\beta, \alpha+\beta, \alpha+2\beta, \alpha+3\beta, 2\alpha+3\beta\}$ assemble into a poset according to their evaluation on $\tau$. 
We make some general observations before considering cases:
\begin{enumerate}[label=(\alph*)]
    \item $\V[\langle e_{2\alpha+3\beta}\rangle] \subset \V[\frakp^{\C}]$. Indeed, $2\alpha+3\beta$ is necessarily maximal in the poset, which implies $ \V[\langle e_{2\alpha+3\beta}\rangle]\subset \ker(\ad_{\Phi})$, so that the desired containment follows from  \eqref{eq:AppendixHypercohomology}. 
    \item $\V[\g^{\C}]_0 \subset \V[\frakk^{\C}]_0$. This is a consequence of the Cartan decomposition and the fact that $\V[\g^{\C}]_0$ must preserve the grading \eqref{eq:Grading}. Hence, $\V[\frak{a}^{\C}]_0 \subset \V[\frak{k}^{\C}]_0$.  
    \item For any root $\gamma$, we have $\V[\langle e_{\gamma}, e_{-\gamma}\rangle ]$ completely contained in $\V[\frakk^{\C}]$ or $\V[\frakp^{\C}]$, by the previous point, and the Cartan decomposition. 
\end{enumerate}
\medskip 

\textbf{Case 1.1: Regular case.} Let us assume that $\tau$ is \emph{fully regular}, meaning that $\alpha(\tau),\beta(\tau) > 0$. Up to the action of the Weyl group, we suppose $\alpha =r^*-s^*$ and $\beta = s^*$, so  
that $\tau$ obtains the form 
\[ \tau = \diag(r+s, r, s,0,-s,-r,-r-s),\]
where $r > s > 0$. First, note that since $\g$ is rank two, the eigenspace $\mathcal{E}[\mathfrak{g}^\C]_0$, a Cartan subalgebra, is of dimension two. Note that $\mathcal{E}[\mathfrak{g}^\C]_0=\mathcal{E}[\mathfrak{k}^\C]_0$  by (b). 

Suppose, for contradiction, that $\V[\langle e_{\alpha+3\beta}\rangle] \subset \V[\frakp^{\C}]$. Hence, $\V[\langle e_{-\alpha-3\beta}\rangle] \subset \V[\frakp^{\C}]$ by (c). Using the Cartan decomposition, one finds $\V[\langle e_{\alpha}\rangle ] \subset \V[\frakk^{\C}]$ and then $\V[\langle e_{-\alpha}\rangle ] \subset \V[\frakk^{\C}]$. One then finds that $\V[\langle e_{\beta}\rangle ] \subset \V[\frakp^{\C}]$. This then implies the whole Cartan decomposition. 
Namely, 
\[\frakk^{\C} = \langle e_{\alpha+2\beta}, e_{\alpha}, e_{-\alpha}, e_{-\alpha-2\beta}\rangle \oplus \mathfrak{a}^{\C}.\]
Now, we know by Proposition \ref{prop:AppendixPConditionSmooth} that $\ad_{\Phi}$ is an isomorphism on the two positive $\frakk^{\C}$-weight spaces:  
\[ \phi_{(\alpha+2\beta)(\tau)}:\V[\langle e_{\alpha+2\beta}\rangle]\stackrel{\cong}{\longrightarrow} \V[\frakp^{\C}]_{(\alpha+2\beta)(\tau)+1}, \qquad  
\phi_{\alpha(\tau)}:\V[\langle e_{\alpha} \rangle]\stackrel{\cong}{\longrightarrow} \V[\frakp^{\C}]_{\alpha(\tau)+1}.\]
The first and second largest positive weight space in $\mathfrak{k}^{\C}$ must be mapped to the largest and second largest positive weight space in $\mathfrak{p}^{\C}$ by $\Phi$. Hence $(\alpha+2\beta)(\tau)+1=(2\alpha+3\beta)(\tau)$ and $\alpha(\tau)+1=(\alpha+3\beta)(\tau)$, so $\alpha(\tau)=2/3$ and $\beta(\tau)=1/3$. Hence, we learn $\Phi\in \mathcal{E}[\mathfrak{g}^\C]_1=\mathcal{E}[e_{\alpha+\beta}]$. However, then we obtain a contradiction, as $\phi_{\frac{2}{3}}: \V[\langle e_{\alpha}\rangle]\rightarrow \V[\langle e_{\alpha+3\beta}\rangle]$ is supposed to be an isomorphism, which is impossible. Hence, this case does not occur. See Figure \ref{Fig:WeightsCase1.1}.

Thus, $\V[\langle e_{\alpha+3\beta}\rangle] \subset \V[\frakk^{\C}]$. Following similar reasoning, we find $\V[\langle e_{\beta} \rangle\oplus \langle e_{-\beta} \rangle ]\subset \V[\frakp^{\C}]$. This leads to $\V[e_{\alpha+2\beta}] \subset \V[\frakp^{\C}]$ and $\V[\langle e_{\alpha+\beta}\rangle ] \subset \V[ \frakk^{\C}]$. Now, we know the weights $\alpha+\beta = r$ and $ \alpha+3\beta = r+2s$ are positive. By Proposition \ref{prop:AppendixPConditionSmooth}, we conclude that $r= s=1$ is mandatory,  $\V[\frakp^{\C}]_1 = \V[\langle e_{\alpha}, e_{\beta}\rangle ]$, and 
$\Phi$ has non-vanishing $e_{\alpha}$ and $e_{\beta}$-components. This completely determines the Higgs bundle $(\mathcal{E},\Phi)$, which is the $\C^*$-fixed point in the Hitchin component in $\mathcal{M}_{\Sigma}(\Gtwosplit)$. 
\begin{figure}[ht]
    \includegraphics[width=.60\textwidth]{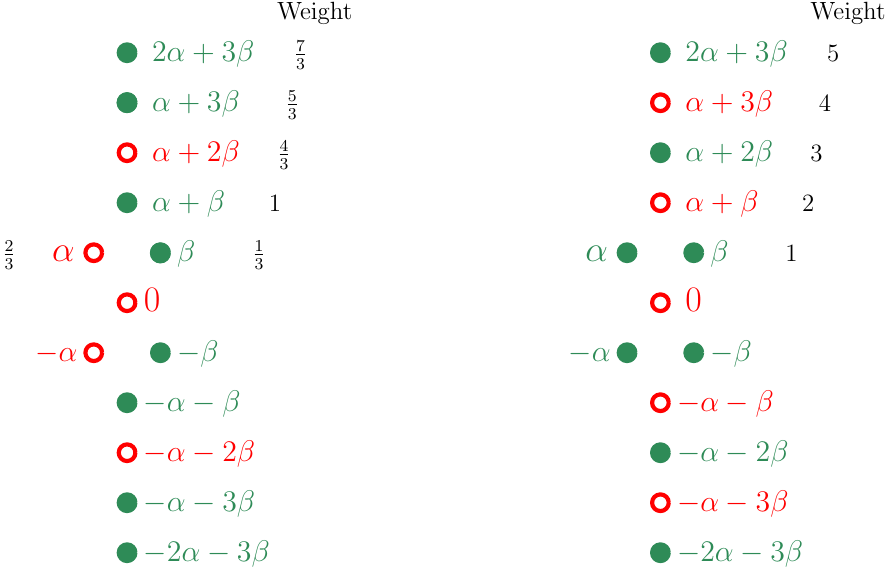}
    \caption{The weights and reality conditions in Case 1.1. On the left is the impossible configuration found if $e_{\alpha+3\beta} \in \p^{\C}$ and on the right is the Fuchsian-Hitchin configuration. The green filled nodes are in $\p^{\C}$ and the red empty nodes in $\frakk^{\C}$.}
    \label{Fig:WeightsCase1.1}
\end{figure}

\textbf{Case 1.2: $\alpha$-case.}  Let us assume that $\beta(\tau)= 0$ and $\alpha(\tau)\neq 0$. In this case, the eigenvalues of the adjoint of $i\psi$ are $(-2x,-x,0,x,2x)$, with corresponding eigenspaces of dimensions $(1,4,4,4,1)$. 
See Figure \ref{fig:Weights1.2} for this case.

\begin{figure}[ht]
    \centering
    \includegraphics[width=0.35\linewidth]{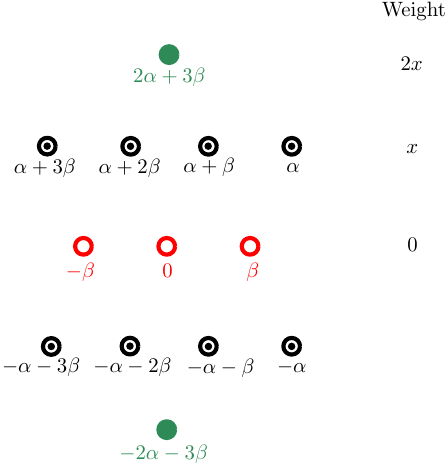}
    \caption{The weight lattice in Case 1.2. The green filled nodes are in $\p^{\C}$, red empty nodes in $\frakk^{\C}$ and the black dotted nodes are weight spaces of unknown type. }
    \label{fig:Weights1.2}
\end{figure}

Moreover, we know that $1$ is an eigenvalue of $i\psi$, which then presents two sub-cases. \medskip 

\textbf{Case 1.2(i): $\mathbf{2x=1}$.} In this case, the eigenvalues of $\ad_{i\psi}$ are 
$(-1,-\frac{1}{2},0,\frac{1}{2}, 1)$. 
We can easily discard this case. 
First, $\V[\g^{\C}]_{1} \subset \V[\frakp^{\C}]$ by (a) and 
$\dim(\V[\frakk^{\C}]_0)=\dim(\V[\g^{\C}])_0=4$ by (b). 
Hence, $\dim(\V[\frakk^{\C}]_{\frac{1}{2}}) > 0$. 
However, $\phi_{\frac{1}{2}}$ cannot be an isomorphism 
since $\Phi \in \V[\g^{\C}]_{1}=\V[\langle e_{2\alpha+3\beta}\rangle ]$. Thus, 1.2(i) does not transpire.

\textbf{Case 1.2(ii): $\mathbf{x=1}$.} 
Hence, $\tau = \diag(1,1,0,0,0,-1,-1) \in \overline{\mathfrak{a}}^+$ in this case and the eigenvalues of $\ad_{i\psi}$ are $(2,1,0,-1,-2)$. 

In particular, $e_{\beta} \in \g_0^{\C} = \frakk^{\C}_0$ by (b). 

Here, we know that $\V[\g^{\C}]_2=\V[\frakp^{\C}]_2$. For  $\phi_1: \V[\frakk^{\C}]_1 \rightarrow \V[\frakp^{\C}]_2$ to be an isomorphism, we must have that
$\dim(\V[\frakk^{\C}]_1)  = 1$. 
Thus, $\dim(\V[\frakp^{\C}]_1)=3$.

The adjoint action of $e_\beta$ preserves $\mathfrak{g}_1\simeq \R^4$ and is nilpotent of order $4$, which means that it does not preserve any non-trivial splitting of $\R^4$. On the other hand, $e_{\beta}\in \mathfrak{k}^{\C}$ is also supposed to preserve the $(3+1)$ splitting $(\mathfrak{k}^{\C}\cap\mathfrak{g}_1)\oplus (\mathfrak{p}^{\C}\cap\mathfrak{g}_1)$. Hence, Case 1.2(ii) does not transpire. \medskip 

We conclude that Case 1.2 does not occur: no such local minima exist. 
 \medskip 

\textbf{Case 1.3: $\beta$-case.}  Let us assume $\alpha(\tau)=0$ and $\beta(\tau)\neq 0$. In this case, the eigenvalues of the form $\ad_{i\psi}$ are of $(-3x,-2x,-x,0,x,2x,3x)$, where $x=\beta(\tau)>0$. Here, we know that the Cartan projection $\tau$ obtains the form $\tau = \diag(2x,x,x,0,-x,-x,-2x)$. Hence, the weight spaces of $\ad_{i\psi}$ are explicitly known: they are given by the $\beta$-height splitting. 
However, we do not yet know the exact weights yet. 
See Figure \ref{fig:Weights1.3}. 

\begin{figure}[ht]
    \centering
    \includegraphics[width=0.24\linewidth]{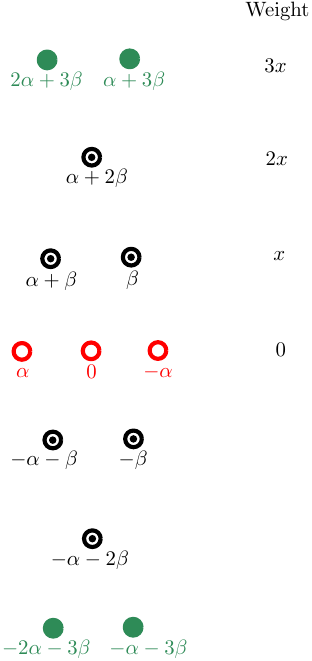}
    \caption{The weight lattice in Case 1.3. The green filled nodes are in $\p^{\C}$, red empty nodes in $\frakk^{\C}$, and black dotted nodes are weight spaces of unknown type.}
    \label{fig:Weights1.3}
\end{figure}

Since we know $1$ is a weight, there are only three options. We proceed as follows: \medskip 

\textbf{1.3(i) Suppose $3x=1$.} This case is easily dispelled by the same reasoning as in Case 1.2(i). \medskip 

\textbf{1.3(ii) Suppose $2x=1$.} The corresponding eigenspaces of $\ad_{i\psi}$ in $\g^{\C}$ are $(-\frac{3}{2}, -1, -\frac{1}{2}, 0, \frac{1}{2}, 1, \frac{3}{2})$, with  dimensions $(2,1,2,4,2,1,2)$. In this case, we must see $\V[\g^{\C}]_{\frac{3}{2}}\subset \V[\frakp^{\C}]$.

Thus, we must have 
$\phi_{\frac{1}{2}}: \V[\frakk^{\C}]_{\frac{1}{2}} \rightarrow \V[\frakp^{\C}]_{\frac{3}{2}}$ an isomorphism by Proposition \ref{prop:AppendixPConditionSmooth}. 
Hence, $\dim(\V[\frakk^{\C}]_{\frac{1}{2}})= 2$, which then implies $\dim(\V[\frakk^{\C}])\geq 8$ by (b) and (c), a contradiction. (Recall that $\dim_{\C}(\frakk^{\C})=6$)
\medskip 

\textbf{1.3(iii) Suppose $x=1$.} Hence, $\ad_{i\psi}$ has eigenvalues $(3,2,1,0,-1,-2,-3)$. 
By (b), we have $\V[\g^{\C}]_3=\V[\frakp^{\C}]_3$. 
Hence, we have $\phi_2: \V[\frakk^{\C}]_2\rightarrow \V[\frakp^{\C}]_3$ an isomorphism by Proposition \ref{prop:AppendixPConditionSmooth}, which is impossible by dimensions. \medskip

We conclude that Case 1.3 does not occur either. \medskip

So far, we have found that the only $\Gtwosplit$-stable minimum of $E$ occurs in the $\Gtwosplit$-Hitchin component. All that remains is to consider the non-stable case. We will show the only such minima have $\Phi=0$. We follow the standard approach of deforming the complex structure to show that any candidate non-trivial minima of $E$ do not occur. \medskip

\textbf{Case 2: Strictly polystable.} Now let us assume that $(\mathcal{E},\Phi)$ is not $\Gtwosplit$-stable, but is polystable. This means that this bundle reduces to an $L$-polystable Higgs bundle $(\mathcal{E}_0,\Phi_0)$ for some Levi subgroup $L\subset \Gtwosplit$. In particular, we see $L=L_\alpha$ or $L= L_\beta$, where the Levi subgroups are associated to the two roots of $\Gtwosplit$. Each of these Levi subgroups is isomorphic to $\GL(2,\R)$, as we recall momentarily.

We now describe the general structure of $L$-polystable Higgs bundles $(\V_0,\Phi_0)$. In particular, we see two possibilities: $\Phi_0 =0$ or $(\V_0,\Phi_0)$ is obtained from an $\SO(2)$-Higgs bundle $(\mathcal{T},0)$, given by a degree zero holomorphic line bundle $\mathcal{T}$, and a $\GL(2,\R)$-stable bundle of the form $\mathcal{E}_0=\mathcal{L}\oplus \mathcal{L}^{-1}$, with 
\[ \Phi_0=\begin{pmatrix}
0 & 0 \\
\eta & 0
\end{pmatrix}, \]
where $\eta \in H^0(\mathcal{K} \mathcal{L}^{-2})$. That is, the latter Higgs bundles obtain the form $(\mathcal{T},0) \star(\V_0,\Phi_0)$ as in \cite[Section 5.2]{BCGBO24}. Observe that stability implies $\deg(\mathcal{L})> 0$, while existence of $\eta$ implies $\deg(\mathcal{L})\leq g-1$. 

We need only consider the nontrivial case that $\Phi_0\neq 0$ below, and we shall specify the concrete form this general description takes when $L=L_{\beta}$ and when $L=L_{\alpha}$. Below, the $\Gtwosplit$-structure of these Higgs bundles is exactly the same as  for $\beta$-bundles and $\alpha$-bundles as in \cite{DE26b}. \medskip 

\textbf{Case 2($\boldsymbol{\beta})$: $L=L_\beta$.} We briefly recall the structure of $L_{\beta}$. One can realize $L_{\beta}$ as the mutual stabilizer in $\Gtwosplit$ of a pair of transverse isotropic lines $\ell, \ell' \in \Ein^{2,3}$. This subgroup $L_{\beta}$ preserves a 2+2+1+1+1 splitting of $\R^{3,4}$. The representation of $L_{\beta}$ on each of the isotropic two-planes, which are dual, is standard. Hence, a polystable  
$L_{\beta}$-Higgs bundle $(\V,\Phi)$ obtains the form: 
\begin{equation}\label{eq:Appendixalpha}
    \begin{tikzcd}
     \mathcal{T}^2 & \mathcal{T}\mathcal{L} \arrow[r, "\eta"] & \mathcal{T}\mathcal{L}^{-1}&\mathcal{O} &\mathcal{T}^{-1}\mathcal{L}\arrow[r, "\eta"]&\mathcal{T}^{-1}\mathcal{L}^{-1}&\mathcal{T}^{-2}
    \end{tikzcd}.
\end{equation} 

In fact, the Levi subgroup $L_{\beta}$ setwise stabilizes a timelike in $\R^{3,4}$. This implies $L_{\beta}$ is a Lie subgroup of a particular $\Z_2$-extension of $\SL(3,\R)$ in $\Gtwosplit$ that we denote by $\widehat{\SL}(3,\R)$ (see \cite[Corollary 3.11]{CT24}). To avoid discussing $\widehat{\SL}(3,\R)$-Higgs bundles, we can directly perform the analogous deformation that one would apply if we were truly in $\SL(3,\R)$.

To define the deformation, fix a nonzero smooth section $\theta \in \Omega^{0,1}(\Sigma, \mathcal{T}^{-1}\mathcal{L}^{-1})$. Since $\deg(\mathcal{L}) > 0$, we can demand $[\theta]$ is nonzero in $H^{0,1}(\Sigma,\mathcal{T}^{-1}\mathcal{L}^{-1})\simeq H^0(\Sigma, \mathcal{L}\mathcal{T}\K)^*$, identified via Serre duality.  
We consider the bundle $(\mathcal{E}',\phi)$ where $\mathcal{E}'$ has the smooth structure of $\mathcal{E}$ and the modified holomorphic structure $\delbar_{\V'}:=\overline{\partial}_\mathcal{E}+s\Theta$, where $s\in\R^*$ is a real number to be chosen and $\Theta \in \Omega^{0,1}(\V[\frakk^{\C}])$ is given by 

\[ \Theta = \begin{pmatrix} &&&&&&\\
                            &&&&&& \\
                            \theta&&&&&& \\
                            &i\sqrt{2}\theta&&&&& \\
                            &&0&&& \\
                            &&&-i\sqrt{2}\theta&&&\\
                            &&&&-\theta&&
    
\end{pmatrix}\]
The deformation $\Theta$ is modeled on the root vector $e_{-\alpha-\beta} \in \g_2^{\C}$, and indeed $\Theta \in \Omega^{0,1}(\V[\frakk^{\C}])$. 
Since $\Theta$ commutes with $\Phi$, the latter being modeled on $e_{-\alpha}$ fiberwise, we find that $\Phi \in H^0(\End(\V') \otimes \K)$ and moreover $\Phi \in H^0(\V'[\mathfrak{p}^{\C}]\otimes \K)$. Thus, $(\V',\Phi)$ is also a $\Gtwo'$-Higgs bundle. Moreover, the $\Gtwosplit$-Higgs bundle $(\V',\Phi)$ is semistable for $s$ small enough as semistability is an open condition. 
For any $t \in \R^*$, consider the gauge transformation $g_t \in \Omega^0(\Sigma, \End(\V))$ given by:
\[ g_{t} := \diag(t^3,t^2,t,1,t^{-1},t^{-2},t^{-3}). \]
Note the following two key properties:  
\[ \Ad_{g_t}\cdot \Phi = t^{-1}\Phi, \qquad \Ad_{g_t} \cdot (\delbar_{\V} +\Theta )= \delbar_{\V} +t^{-2}\Theta. \] 
Hence, in the moduli space $\mathcal{M}_{\Sigma}(\Gtwosplit)$, we see $\lim_{t \rightarrow \infty} [(\mathcal{E}',t\Phi)]=[(\mathcal{E},\Phi)]$. By \cite[Lemma A.20]{AABCGGPO19}, we conclude that the original Higgs bundle $(\V,\Phi)$ is not a local minimum for $E$ unless $\Phi=0$. Here, we may invoke this lemma as we already know that for all $t\in \R$ the bundle $(\mathcal{E}',t\Phi)$ is not $\mathcal{S}$-equivalent to $(\mathcal{E}',t\Phi)$, when $\Phi \neq 0$, since $\theta$ is nonzero in cohomology. \medskip  

\textbf{Case 2($\boldsymbol{\alpha})$: $L=L_\alpha$.} We now recall the structure of $L_{\alpha}$. The Levi subgroup is realized as $L_{\alpha} = \Stab_{\Gtwosplit}(\omega ) \cap \Stab_{\Gtwosplit}(\omega')$ for a pair of transverse annihilator photons $\omega, \omega' \in \Pho^\times$. The subgroup $L_{\alpha}$ preserves a 3+2+2 decomposition of $\R^{3,4}$, namely 
$\omega\times_{\imoct} \omega' \oplus \omega \oplus \omega'$. Here, $\omega,\omega'$ are dual isotropic 2-planes and $\omega\times_{\imoct} \omega' \cong \R^{1,2}$. Correspondingly, an $L_{\alpha}$-Higgs bundle obtains the form: 
\begin{equation}
    \begin{tikzcd}\label{eq:AlphaCandidate}
     \mathcal{T}\mathcal{L} \arrow[r, "\eta"] & \mathcal{T}\mathcal{L}^{-1}&\mathcal{L}^2\arrow[r, "-i\sqrt{2}\eta "]&\mathcal{O} \arrow[r, "-i\sqrt{2}\eta "]&\mathcal{L}^{-2}&\mathcal{T}^{-1}\mathcal{L}\arrow[r, "\eta"]&\mathcal{T}^{-1}\mathcal{L}^{-1}
    \end{tikzcd} .
\end{equation}  

Now, the subgroup $L_{\alpha}$ is also a Lie subgroup of $\SO(2,2)$, which is realized by $H=\Stab_{\Gtwosplit}(\omega\oplus \omega')$. One difficulty of this case is that the above Higgs bundle \emph{is} a local minimum of $E$ in the smaller space $\mathcal{M}_{\Sigma}(\SO(2,2))$ \cite[Proposition 3.12]{AABCGGPO19}. Thus, we must deform outside of $\SO(2,2)$. 

To address the candidate Higgs bundle \eqref{eq:AlphaCandidate}, we explicitly construct a path of Higgs bundles for which the functional $E$ decreases, as in Case 2($\beta$).
To define the deformation of Higgs bundles, we once again fix the Higgs field and deform only the complex structure. 
This time, we select $\theta \in \Omega^{0,1}(\Sigma, \mathcal{T}^{-1}\mathcal{L}^{-3})$. Again by Serre duality, we see that $\theta$ can be chosen to be non-trivial in cohomology.  
Consider the bundle $(\mathcal{E}',\Phi)$ where $\mathcal{E}'$ is the smooth bundle $\mathcal{E}$ equipped with the holomorphic structure $\overline{\partial}_\mathcal{E}+s\Theta$, where $s \in \R^*$ is a constant to be chosen and $\Theta \in \Omega^{0,1}(\V[\frakk^{\C}])$ is

\[ \Theta = \begin{pmatrix} &&&&&&\\
                            &&&&&& \\
                            &&&&&& \\
                            &&&&&& \\
                            \theta&&&&&& \\
                            &0&&&&&\\
                            &&-\theta&&&&
    
\end{pmatrix}.\]
The deformation $\Theta$ is modeled on the root vector $e_{-\alpha-3\beta} \in \g$, and satisfies $\Theta \in \Omega^{0,1}(\V[\frakk^{\C}])$. 
Since $\Theta$ commutes with $\Phi$, the latter modeled on $e_{-\beta}$, we have $\Phi \in H^0(\End(\V') \otimes \K)$. Thus, $(\V,\Phi')$ is also a $\Gtwo'$-Higgs bundle, which is semistable for $s$ sufficiently small.\footnote{Note that as an $\SO(2,2)$-Higgs bundle, $\Theta$ does \emph{not} correspond to a section of $\Omega^{0,1}(\V[\frakk^{\C}])$.} 

For any $t \in \R^*$, we consider the gauge transformation $g_t \in \Omega^0(\Sigma, \End(\V))$ given by 
\[ g_{t} := \diag(t^3,t^2,t,1,t^{-1},t^{-2},t^{-3}) .\] 
Note that 
\[ \Ad_{g_t}\cdot \Phi = t^{-1}\Phi, \qquad  \Ad_{g_t} \cdot (\delbar_{\V} +\Theta )= \delbar_{\V} +t^{-4}\Theta.\]
Thus, in the moduli space $\mathcal{M}_{\Sigma}(\Gtwosplit)$, we see $\lim_{t \rightarrow \infty} [(\mathcal{E}',t\Phi)]=[(\mathcal{E},\Phi)]$. 
If $\Phi \neq 0$, then since $(\V',\Phi)$ is not $\mathcal{S}$-equivalent to $(\V,\Phi)$ as in the previous case, we conclude again by \cite[Lemma A.20]{AABCGGPO19} that the original Higgs bundle $(\V,\Phi)$ is not a local minimum for $E$. This completes Case 2($\alpha$) and the proof of Lemma \ref{Lem:G2Minima}. \end{proof}

We now finally verify the connected component count using the key Lemma \ref{Lem:G2Minima}. 

\begin{proof}[{Proof of Theorem \ref{thm:AppendixCountConnectedComponents}}]
Since $E$ is proper, each connected component of the character variety contains a local minimum. Now, we know that $|\pi_0(\chi(S,\Gtwosplit))|\geq 1 + |\pi_0(K)| =3$ because of topological invariants associated to $K =\SO(4,\R)$ and due to the Hitchin component. The components $\mathcal{C}_0$ and $\mathcal{C}_1$ contain representations $\rho_i:\pi_1S \rightarrow K \hookrightarrow  \Gtwosplit$ associated to Higgs bundles $(\V_i, 0)$ with zero Higgs field. These Higgs bundles are necessarily the global minima of the function $E|_{\mathcal{C}_i}$. We know also that the Fuchsian-Hitchin case in $\mathcal{C}_{\Hit}$ corresponds to the unique $\C^*$-fixed point in $\Gtwosplit$-Hitchin section of $\Sigma$, which is also the global minimum of $E$ on $\mathcal{C}_{\Hit}$. 
On the other hand, Lemma \ref{Lem:G2Minima} asserts that these three global minima are the only local minima of $E$ on $\chi(S,\Gtwosplit)$, meaning $|\pi_0(\chi(S,\Gtwosplit))|=3$. 
\end{proof}

Finally, we can sharpen Theorem \ref{thm:StrataCorrespond}.

\begin{theorem}[Components of Anosov representations]
Let $S$ be a closed surface of genus $g\geq 2$ and 
write $\chi(S,\Gtwosplit) = \mathcal{C}_0 \sqcup \mathcal{C}_1 \sqcup \mathcal{C}_{\Hit}$ as in Theorem \ref{thm:AppendixCountConnectedComponents}. Then we have the following: 
\begin{enumerate}[label=(\roman*),itemsep = 1ex]
    \item The set of $P_{\alpha}$-Anosov representations in $\mathcal{C}_0$ and $\mathcal{C}_1$ has at least $2g-2$ connected components each.
    \item The set of $P_{\beta}$-Anosov representations in $\mathcal{C}_0$ and $\mathcal{C}_1$ has at least $3g-4$ and $3g-3$ connected components, respectively. 
\end{enumerate} 
Moreover, for $\alpha$ and $\beta$-bundles, the associated representation is in $\mathcal{C}_0$ if and only if $d$ is even, where $d=\deg(\mathcal{T})$ or $d=\deg(\mathcal{B})$, respectively.
\end{theorem}

Note that in the count we must exclude the case $d=6g-6$ for $\beta$-bundles and $d=2g-2$ for $\alpha$-bundles, as the corresponding representations are in the Hitchin component.

\begin{proof}
By Lemma \ref{Lem:ConnectedComponents}, Theorem \ref{thm:StrataCorrespond} and Theorem \ref{thm:AppendixCountConnectedComponents}, it suffices to prove the `moreover' claim. Recall the Frenet frame splitting $\V^{\R} \cong (\mathscr{L}\oplus N) \oplus (T \oplus B)$ for $\sigma$-bundles from \S\ref{Subsec:SigmaBundles}. We use that each of $T,N,B$ is a real plane bundle that also has the structure of a complex line bundle. Also, as a real vector bundle, $E^4_{\rho} \cong T\oplus B$. Recall that for $E \rightarrow S$ a real rank four vector bundle, equipped with the additional structure of a complex vector bundle, the second Stiefel-Whitney class $w_2(E) \in H^2(S, \Z_2) \cong \Z_2$ is equal to $\deg(E)\bmod 2$, where $\deg$ is the degree of $E$, viewed as a complex vector bundle. Hence, we may conclude by describing the complex structures on $T$, $B$ in each case. \medskip 

For $\beta$-bundles, we have $\mathcal{B} \cong B$ and $T \cong \mathcal{K}^{-1}$. Thus, $\deg(E^4_{\rho}) \bmod 2 =\deg(\mathcal{B}) \bmod 2$. \medskip 

For $\alpha$-bundles, we have $T \cong \mathcal{T}^{-1}$ and 
$B \cong \mathcal{T}^{2}\K$. Hence, $\deg(E^4_{\rho}) \bmod 2 = \deg(\mathcal{T}) \bmod 2$. 
\end{proof}

\bibliographystyle{alpha}
\small{\bibliography{bib}}

\end{document}